\documentclass[11pt,letterpaper, reqno]{amsart}
\usepackage{amsmath, amsthm, amsfonts, amssymb, color}
\usepackage{mathrsfs, bbm}

\allowdisplaybreaks

\newtheorem{thm}{Theorem}[section]
\newtheorem{cor}[thm]{Corollary}
\newtheorem{lem}[thm]{Lemma}
\newtheorem{remark}[thm]{Remark}
\newtheorem{prp}[thm]{Proposition}

\newtheorem{defn}[thm]{Definition}
\newcommand{\scr}[1]{\mathscr #1}
\definecolor{wco}{rgb}{0.5,0.2,0.3}

\numberwithin{equation}{section} \theoremstyle{remark}

\def\1{{\mathbbm 1}}

\def\wt{\widetilde}

\def\R{\mathbb R}

\def\E{\mathbb E}

\def\sX{\mathcal X}
\def\sL{\mathcal L}
\def\L{\mathcal L}

\def\wh{\widehat}

\def\<{\langle}
\def\>{\rangle}

\def\pf{\noindent{\bf Proof.} }

 \def\beq{\begin{equation}}

\def\B{\scr B}

 \def\P{\mathbb P} 
  \def\ee{\varepsilon}

\def\sL{\mathcal L}

\begin{document}

\title{Boundary Harnack principle for non-local operators on metric measure spaces }


\author{Zhen-Qing Chen}
\address{Department of Mathematics, University of Washington, Seattle, WA 98195, USA}
\curraddr{}
\email{zqchen@uw.edu}
\thanks{}

\author{Jie-Ming Wang}
\address{School of Mathematics and Statistics,  Beijing Institute of Technology,
Beijing 100081, P. R. China}
\curraddr{}
\email{wangjm@bit.edu.cn}
\thanks{Research partially supported  by
 NNSFC   Grant   11731009.}



\date{}

\dedicatory{}

\begin{abstract}
In this paper,  a necessary and sufficient condition is obtained for the scale invariant boundary Harnack inequality (BHP in abbreviation) for
 a large class of Hunt processes on metric measure spaces that are in weak duality with another  Hunt  process.
We next consider
 a discontinuous subordinate Brownian motion with Gaussian component $X_t=W_{S_t}$ in $\R^d$
for which the   L\'evy density of the subordinator $S$  satisfies  some mild comparability condition.
We   show that   the scale invariant BHP  holds for the subordinate Brownian motion $X$ in any Lipschitz domain
 satisfying the interior cone condition with common angle $\theta\in (\cos^{-1}(1/\sqrt d), \pi)$,  but fails
  in any truncated circular cone  with angle $\theta \leq \cos^{-1}(1/\sqrt d)$,
    a Lipschitz domain  whose Lipschitz constant   is larger than or equal to $1/\sqrt{d-1}.$
\end{abstract}

\maketitle

\section{Introduction}

Let $d\geq 1$.  For $x\in \R^d$ and $r>0$, denote by $B(x, r)$   the open ball in $\R^d$ centered at $x$ with radius $r$.
We say a boundary Harnack principle (BHP in abbreviation) holds
 for an operator $\L$  on $\R^d$ in  a domain $D\subset \R^d$
if  for any nonnegative   functions $u$ and $v$
that are $\L$-harmonic in $D\cap B(z_0, r)$ and vanish   on $B(z_0, r)\cap D^c$  with  $z_0\in \partial D$
 in suitable sense,
there is a constant $c(r)>0$ independent of $u$ and $v$ so that
$$
\frac{u(x)}{v(x)} \leq c\frac{u(y)}{v(y)} \quad \hbox{for any }
x, y\in  B(z_0, r/2)\cap  D.
$$
 If the constant $c(r)$ can be chosen to be independent of $r\in (0, r_0]$ for some $r_0>0$,
we say the scale invariant BHP holds for $\L$ in $D$;
see \eqref{e:1.6} for a precise statement.
 BHP, when it holds, has many far reaching implications.
 Especially,
 the scale invariant BHP  is crucial in the study of Martin boundary, Fatou theorem and Green function estimates of $\sL$ in open sets.

 Scale invariant BHP  for the Laplace operator $\Delta$ (or equivalently for Brownian motion)
 in Lipschitz domains in $\R^d$  is  independently established  by
  Ancona \cite{A},  Dahlberg \cite{Da} and Wu  \cite{Wu}.
  This BHP  result
   is    later extended
 to  non-tangential accessible (NTA) domains by Jerison and Kenig \cite{JK},
 to more general  elliptic operators in Lipschitz domains  by Caffarelli, Fabes, Mortola and Salsa
 \cite{CFMS} and Fabes, Garofalo, Marin-Malave
and Salsa \cite{FGMS},
and to  H\"older domains  for elliptic operators by Ba\~nuelos, Bass and Burdzy \cite{BBB} and  Bass and Burdzy  \cite{ BB2, BB4}.
The scale invariant BHP for the Laplace operator in uniform domains  is   established by Aikawa \cite{Ai}.
There is  now
 a very substantial literature on BHP for differential operators.
 For diffusions on geodesic spaces,  the scale invariant BHP on inner uniform domains have been established in \cite{GSC, LS, Lie} under a  two-sided  Guassian or sub-Gaussian heat kernel bounds assumption.
In \cite{BM}, it is shown that  for symmetric diffusions on geodesic spaces,  the scale invariant BHP on inner uniform domains holds under the assumption of scale invariant elliptic Harnack inequality and the existence of regular Green function. By a recent result of \cite{BCM}, the assumption of the existence of regular Green function can be dropped.

\smallskip

The study of BHP for non-local operators (or equivalently, discontinuous Markov processes)
started relatively recent and many  progresses have been made in recent years.
  Bogdan \cite{Bo} showed the scale invariant BHP holds for $\Delta^{\alpha/2}:=-(-\Delta)^{\alpha/2}$
(or equivalently, for rotationally symmetric $\alpha$-stable processes)
 with $\alpha\in (0,2)$ in bounded Lipschitz domains.
This result is later extended to $\kappa$-fat open sets by Song and Wu \cite{SW} and then to arbitrary open sets
by Bogdan, Kulczycki and Kwasnicki in \cite{BKK}.
The scale invariant BHP has been shown to hold for   pure jump subordinate Brownian motions in
 $\kappa$-fat  (or corkscrew)
open sets and for
a  class of pure jump symmetric L\'evy process in open sets  \cite{KSV, KSV5}.
BHP has also been established for some non-homogenous non-local operators that are generators
of some discontinuous processes. They include  censored stable processes \cite{BBC}, truncated stable processes (i.e. isotropic stable processes   with jump size larger than $1$ removed) \cite{KS},   a class of discontinuous Markov processes whose generators are fractional Laplacians under a lower order perturbation \cite{CRY}, and
   stable subordinate of Brownian motion  on Sierpi\'nski  carpet \cite{Sto} and
on infinite Sierpi\'nski gasket  \cite{KK}.
  A non-scale-invariant BHP for stable-like operator in non-divergence form is given in \cite{ROS}.

\smallskip

Despite the significant progress for
 differential elliptic operators and purely non-local operators,
 results on BHP
 for non-local operators with diffusive part, or equivalently, for
discontinuous processes with Gaussian components,   are still  quite limited at present, especially on non-smooth domains.
Chen, Kim, Song and Vondra\v{c}ek \cite{CKSV} established a scale invariant boundary Harnack principle with explicit boundary decay rate  in $C^{1,1}$
  open sets
  for $\L^b :=\Delta- b(-\Delta)^{\alpha/2}$ that is uniform in the constant $b\in (0, M]$.
  This result is  later extended in \cite{KSV2} to
\begin{equation}\label{e:1.1a}
\L^\phi:= -\phi (-\Delta)
\end{equation}
	in $C^{1,1}$  open sets,
	where $\phi$ is a complete Bernstein function on $[0, \infty)$
  that has a positive drift $b>0$ 	  (see  Remark \ref{R:3.5}(i) for its definition)
	satisfying $\phi (0)=0$ and some additional conditions.
	 	Associate with such non-local operator $\L^\phi$ is a L\'evy process that is the independent sum of a Brownian motion and a  purely discontinuous  subordinate Brownian motion.
	The scale invariant boundary Harnack principle with explicit boundary decay rate  in $C^{1,1}$ open sets for a large class of diffusions with jumps, which may not be L\'evy processes,
	is obtained recently by  Chen and Wang \cite{ChW}.
A non-scale-invariant BHP  for jump-type Feller
processes (which can possibly have diffusive components)  in weak dual having strong Feller property
on metric measure state spaces is obtained by Bogdan,  Kumagai and  Kwasnicki \cite{BKK2}, under a comparability condition  of
the jump kernel and a Urysohn-type condition on the domain of the generator of the process.
 It was an open problem till now if the scale invariant BHP holds for
$\Delta + \Delta^{\alpha/2}$ in Lipschitz domains or not.

\smallskip

One of the goals of this paper is to answer this open problem.
We show that
the scale invariant BHP holds for $\Delta + \Delta^{\alpha/2}$  in Lipschitz domains in $\R^d$
satisfying the interior cone condition with common angle $\theta\in (\cos^{-1}(1/\sqrt d), \pi)$,
but fails  in any truncated circular cone  with angle $\theta \leq \cos^{-1}(1/\sqrt d)$,
    a Lipschitz domain  whose Lipschitz constant   is larger than or equal to $1/\sqrt{d-1}.$
In fact, we have established this result not only for $\Delta + \Delta^{\alpha/2}$ but also for more general
non-local operator $\L^\phi$ of \eqref{e:1.1a}. See Theorem \ref{T:1.10} below for a precise statement.
Our approach is largely probabilistic but also uses some analytic techniques.
To show the above result,  we investigate and establish, under some mild assumptions,
 the necessary and sufficient conditions of the scale invariant BHP for discontinuous Hunt processes on  open sets of  locally compact separable metric spaces. This result is then applied to get the  geometric conditions on Lipschitz domains
for which the scale invariant BHP holds for  $\L^\phi$.

\smallskip

We now carefully state the setting and the main results of this paper.
 Let $(\sX, \B, d)$ be a locally compact separable metric space with $\sigma$-algebra $\B$   and  $m$ be a Radon measure  on $\sX$ that has full support.
Let $(X_t, \zeta, \mathcal{F}_t, \P_x)$ be an  irreducible Hunt  process taking values in  $\sX$
in the sense of \cite{BG}.  Here we say  the Hunt process $X$ is  irreducible on $\sX$ if for any non-empty  open set $U\subset \sX,$
$ \P_x (T_U<\infty)>0 $ for every   $x\in \sX,$
where $T_U:=\inf\{t>0: X_t\in U\}.$
The $\alpha$-potential operator $G^\alpha$ of $X$ is defined by
\begin{equation}\label{e:1.1}
G^\alpha f(x):=\E_x \int_0^\infty e^{-\alpha t} f(X_t)\,dt, \quad \alpha> 0, x\in \sX
\end{equation}
for each bounded measurable function $f.$
 We assume that $X$ is  in weak duality to another
   Hunt process $\wh X$ with respect to  the Radon measure $m$ on $\sX$   in the following sense.
For each  nonnegative Borel  functions $f, g$ in $\sX$ and for each $\alpha>0,$
\begin{equation}\label{e:1.1'}
\int_\sX  G^\alpha f(x)\cdot g(x) m(dx)=\int_\sX  f(x) \wh G^\alpha g(x) m(dx),
\end{equation}
where $\wh G^\alpha$ is the $\alpha$-potential operator of $\wh X$.
Denote by $G^\alpha(x, \cdot)$ and $\wh G^\alpha(x, \cdot)$ the $\alpha$-potential measures of the processes $X$ and $\wh X$, respectively.

 \smallskip
 It is known that any Hunt process admits a L\'evy system that describes how the process jumps.
Suppose that the Hunt process $X$ has the L\'evy system $(N(x, dy), dt).$ That is, for any  nonnegative function $f$ on $\R_+\times
\sX\times \sX$ vanishing along the diagonal of $\sX \times \sX$, for any stopping time $T$
with respect to the minimal admissible  augmented  filtration generated by $X$ and $x\in \sX$,
\begin{equation}\label{e:2.1'}
\E_x  \Big[ \sum_{s\leq T}f(s, X_{s-}, X_s); X_{s-}\neq X_s \Big]
 =  \E_x \Big[ \int_0^T\int_\sX f(s,X_s,y) N(X_{s-}, dy)\,ds \Big].
\end{equation}

\smallskip

  We   consider the following
 assumptions on the  irreducible Hunt process $X$ in this paper.
The first three are on the Hunt process $X$, while the fourth one is a condition on the open subset $D$ of $\sX$.

 \begin{enumerate}
\item[{\bf (A1)}] For all $x\in \sX$ and $\alpha>0,$  $G^\alpha(x, \cdot)$ and $\wh G^\alpha(x, \cdot)$ are absolutely continuous with respect to $m(dy).$

\item[{\bf  (A2)}] Every semi-polar set of $X$ is polar.

\smallskip

\item[{\bf   (A3)}]
There exist $c>1$ and $r_0>0$ such that for any $x_0\in \sX, 0<r<r_0$ and $ x\in B(x_0, r)$,
$$c^{-1}N(x_0, dy)\leq N(x, dy) \leq cN(x_0, dy)   \quad \hbox{on } B(x_0, 2r)^c.
$$

\item[{\bf (A4)}]  $D$ is  an open subset of $\sX$ so that
  there exist $c>1$ and $\delta_0>0$ such that for any $z_0\in \partial D$ and $r\in (0, \delta_0),$
$$
\E_x \left[ \tau_{D\cap B(z_0, 4r)} \right] \leq c \E_x \left[ \tau_{D\cap B(z_0, 2r)}\right] \quad \hbox{for } x\in D\cap B(z_0, r),
$$
where  $\tau_B:=\inf\{t \geq 0: X_t \notin B\}$ denotes the first exiting time from $B$ by  the process $X$.

 \end{enumerate}

  Under assumption {\bf (A1)}, for each $\alpha>0,$ there exists a unique $\alpha$-potential kernel $G^\alpha(x, y)$ defined on $\sX\times \sX$ such that for  each nonnegative Borel  function $f,$
$$
G^\alpha f(x)=\int_\sX G^\alpha(x, y) f(y) m(dy), \quad \wh G^\alpha f(x)=\int_\sX G^\alpha(y, x) f(y) m(dy),
$$
and  $x\mapsto G^\alpha(x, y)$ is $\alpha$-excessive with respect to $X,$  $y\mapsto G^\alpha(x, y)$ is $\alpha$-excessive with respect to $\wh X;$ see \cite[Theorem 13.2]{ChungW},   or \cite[Theorem 1]{KW},  or the remarks after \cite[Proposition VI.1.3]{BG}.
 Assumption  {\bf  (A2)} is known as Hunt's hypothesis (H) for the process $X$.
 It implies that the process $X_t$ never hits irregular points.
  Assumption {\bf  (A3)} is a comparability condition on the jumping kernel  $N(x, dy)$ of $X.$
     Note that when $X$ is a Hunt process having continuous sample paths,
then its L\'evy system is trivial, that is, $N(x, dy)=0$ for every $x\in \sX$.
In this case, condition {\bf  (A3)} is automatically satisfied. In this paper, we are mostly concerned with
discontinuous Hunt processes.

\smallskip

 For an open set $U\subset \sX,$ denote by $X^U$ and $\wh X^U$ the part processes of $X$ and $\wh X$ killed upon exiting $U$ respectively.
  These part processes are Hunt processes on $U$; see, for example,  \cite[Exercise 3.3.7]{CF}.
For each $\alpha>0,$  denote by $G^\alpha_U$ and $\wh G^\alpha_U$ the $\alpha$-potential operators  of $X^U$ and  $\wh X^U,$ which are similarly defined as in \eqref{e:1.1}  with $X^U$ and $\wh X^U$ in place of $X.$
  By \cite[Remark 13.26]{ChungW},  the semigroups  of $X^U$ and  $\wh X^U$ are in weak duality relative to $m(dx).$
  Hence for each $\alpha>0,$ $G^\alpha_U$ and $\wh G^\alpha_U$ are in weak duality relative to $m(dx).$
As $X^U$ and $\wh X^U$ are subprocesses of $X$ and $\wh X,$ it is easy to see that for each $\alpha>0,$  $G^\alpha_U(x, dy)$ and $\wh G^\alpha_U(x, dy)$ are absolutely continuous with respect to $m(dy)$ under Assumption  (A1).
 Hence,  for each $\alpha>0,$ there exists a unique $\alpha$-potential kernel $G^\alpha_U(x, y)$ of $X^U$  on $U\times U$ in this case.
 Since $G^\alpha_U(x, y)$ is non-decreasing in $\alpha>0$ for $(x, y)\in U\times U,$ we define the potential kernel
$G_U(x, y):=\lim_{\alpha\rightarrow 0}G^\alpha_U(x, y)$ for $(x, y)\in U\times U,$ which exists  uniquely but may be infinite under Assumption  (A1).
 We set $G_U(x, y)=0$ for $(x, y)\notin U\times U.$
By monotone convergence theorem, for each nonnegative Borel function $f$ on $U,$
$$G_U f(x)=\int_U G_U(x, y) f(y) m(dy)=\E_x \int_0^\infty f(X^U_t) \,dt=\E_x \int_0^{\tau_U} f(X_t) dt.$$
Since the increasing limit of excessive function is also an excessive function, for each  $y\in U,$ $x\mapsto G_U(x, y)$ is excessive with respect to $X^U,$ and for each  $x\in U,$  $y\mapsto G_U(x, y)$ is excessive with respect to $\hat X^U$.

\begin{defn}\label{De:1}
   For an open set $U$ of $\sX$,  we say $G_U(x, y)$ is
the Green function of $X$ in $U$ if it is the potential kernel of $X^U.$ That is, the following conditions hold:
\begin{enumerate}
 \item[\rm (i)]
 For any Borel bounded function  $f\geq 0$ on $U$,
\begin{equation}\label{e:1.3'}
\E_x \int_0^{\tau_U} f(X_s) ds = \int_U G_U (x, y) f(y) m(dy), \quad \hbox{for} \quad x\in U.
\end{equation}

\item[\rm (ii)]
For each fixed $y\in U,$ $x\mapsto G_U(x, y)$ is excessive with respect to $X^U.$ For each fixed $x\in U,$  $y\mapsto G_U(x, y)$ is excessive with respect to $\wh X^U.$

 \end{enumerate}
\end{defn}

\smallskip

For a Borel set $A\subset \sX$,  define $T_A:=\inf\{t>0: X_t\in A\}.$

\begin{defn}\label{D:1} \rm
A point $y$ is said to be regular of $X$ for a Borel set $A$ if $\P^y(T_A=0)=1.$
 Denote by $A^r$ all the regular points of $A$.
\end{defn}

The following is the first main result of this paper.
 See Definition \ref{D:2.1} below for the terminology of harmonic and regular harmonic functions.

\begin{thm}\label{T:2.1}
Let   $D$  be  an open subset of $\sX$  so that
 $D^c$ has a non-empty interior.
\begin{enumerate}

 \item[\rm (i)] Suppose that Assumptions (A1)-(A4) hold.
 Suppose that  there exist $\delta= \delta(X, D)>0$ and $C=C(X, D)\geq 1$   such that for any  $z_0\in \partial D,$
 $r\in (0, \delta)$ and $x_1, x_2\in D_r(z_0)$,
\begin{equation}\label{e:1.4a}
 G_{D_{4r}(z_0)}(x_1, y)
 \leq C\dfrac{\E_{x_1} [\tau_{D_{2r}(z_0)}]}{\E_{x_2}[ \tau_{D_{2r}(z_0)}]}G_{D_{4r}(z_0)}(x_2, y)
 \end{equation}
 for   $ y\in D_{4r}(z_0)\setminus  D_{3r}(z_0),$
where $D_r(z_0):=D\cap B(z_0, r).$
 Then    the  scale invariant boundary Harnack principle holds for $X$ in $D$;
 that is,
 there exist positive constants $\delta_1=\delta_1(X, D)$ and $C_1   =C_1(X, D)\geq 1 $ such that for any
 $z_0\in\partial D$,  $r\in (0, \delta_1)$ and  any two nonnegative  functions $h_1$ and $h_2$ that are
regular harmonic with respect to $X$  in  $D\cap B(z_0, 2r)$ and  vanish
on $(D^c)^r\cap B(z_0, 2r),$
\begin{equation}\label{e:1.6}
h_1(x)h_2(y)\leq C_1h_1(y)h_2(x)  \quad \hbox{for} \: x,y\in D\cap B(z_0,  r/2).
\end{equation}

\item[\rm (ii)]   Suppose that Assumptions (A1) and (A3) hold.
 Suppose that
     there is some $r_1 >0 $ so that
\begin{equation}\label{e:1.5}
N(x, \sX\setminus B(x,  r_1))>0  \quad \hbox{ for every } x\in \sX.
\end{equation}
  If the  scale invariant boundary Harnack principle \eqref{e:1.6}  holds,
then \eqref{e:1.4a} holds.

 \item[\rm (iii)]  In particular,  suppose that Assumptions (A1)-(A4) and \eqref{e:1.5} hold. The  scale invariant boundary Harnack principle     holds for $X$ in an open set $D$  if and only if  \eqref{e:1.4a} holds.
  \end{enumerate}
\end{thm}

\begin{remark} \label{R:1.4} \rm
\begin{enumerate}

\item[(i)]
 Note that in Theorem \ref{T:2.1}(i), we  do not need to assume $X$ is discontinuous.  However,
the condition that $X$ is discontinuous is essential for  Theorem \ref{T:2.1}(ii).
It is known that the comparability condition \eqref{e:1.4a} is not necessary for the scale invariant BHP for Brownian motion
on Euclidean spaces; see Remark \ref{R:1.5} below.

 \item[(ii)]
  Assumption (A2) is not needed in Theorem \ref{T:2.1}(i) if we consider the scale invariant BHP   only  for those
nonnegative  regular harmonic functions that vanish
 on $D^c\cap B(z_0, 2r)$ in Theorem \ref{T:2.1}(i).
The same remark applies to other results in this paper that involves the regular points of $D^c$, for example,
Proposition \ref{P:2.3}.
On the other hand,  if $h$ is a nonnegative   harmonic function with respect to $X$ in $D\cap B(z_0, 2r)$  vanishing continuously
 on  $(D^c)^r\cap B(z_0, 2r)$ with $z_0\in \partial D$,
  then under  Assumption (A2), by a similar argument as that for \cite[Lemma 4.2]{CKSV},
 $h$ is a regular harmonic function with respect to $X$ in $D\cap  B(z_0, s) $ that vanishes on  $(D^c)^r\cap B(z_0, s)$ for every $s\in (0, 2r)$.
 This is the reason we formulated Theorem \ref{T:2.1}(i) in this way so its conclusion
 is applicable to nonnegative   harmonic functions with respect to $X$ in $D\cap  B(z_0, 2r)$  that vanish continuously
 on $(D^c)^r\cap  B(z_0, 2r)$.
 In Theorem \ref{T:1.10}, the BHP for subordinate Brownian motions   is stated for  nonnegative   harmonic functions  in $D\cap B(z_0, 2r)$  that vanish continuously  on $(D^c)^r \cap  B(z_0, 2r)$ when $D$ is a Lipschitz domain (for which $(D^c)^r = D^c$).

\item[(iii)]  When $\sX=\R^d$,  assumption \eqref{e:1.5} holds for every   discontinuous L\'evy process
$X$ on $\R^d$.
 By Theorem \ref{T:2.1}(iii), under assumption \eqref{e:1.5},
the condition \eqref{e:1.4a} is   necessary and sufficient for the scale invariant BHP to hold for a discontinuous Hunt process $X$ in any open subset $D \subset \sX$ under Assumptions (A1)-(A4).
To see that the condition \eqref{e:1.4a} is  necessary for the scale invariant BHP for the discontinuous Hunt process
$X$ in  any  open set  $D$ ,
 let   $u(x):=\P_x(  X_{\tau_{D_{2r}(z_0)}} \in   B(z_0, 4r)^c )$
   and $v(x):=\E_x \left[ G_{D_{4r}(z_0)}(X_{\tau_{D_{2r}(z_0)}} , y)  \wedge n \right] $  with  $z_0\in \partial D$
  and $y\in D_{4r}(z_0)\setminus D_{3r}(z_0),$  where $D_r(z_0)=D\cap B(z_0, r).$
  The functions   $u$ and $v$ are  nonnegative    regular  harmonic   in $D\cap B(z_0, 2r)$ and   vanish
  on    $   D^c   \cap B(z_0, 2r)$.
   By using the L\'evy  system formula and Assumption (A3) (see Lemma \ref{L:2.2} below),
  we can show that for any $x_1, x_2\in D\cap B(z_0, r),$ the ratio  $u(x_1)/u(x_2)$ is comparable to $\E_{x_1} \tau_{D_{2r}(z_0)}/ \E_{x_2} \tau_{D_{2r}(z_0)}$ uniformly in  small $r$.
   Passing $n\rightarrow\infty$,  this shows that \eqref{e:1.4a} is a necessary condition for the scale invariant BHP to hold on $D$
   for the discontinuous Hunt process $X$.

 \item[(iv)] In fact, we   will show  under the assumption of Theorem \ref{T:2.1}(i)
 that for any non-trivial, positive and bounded function $h $ that is
 $X$-regular harmonic in  $D\cap B(z_0, 2r)$ and vanishes  on  points of  $(D^c)^r\cap B(z_0, 2r),$
$$ \dfrac{h(x )}{h(y)}\leq c\dfrac{\E_{x } \left[ \tau_{D_{2r}(z_0)} \right] }{\E_{y} \left[ \tau_{D_{2r}(z_0)} \right] }
 \quad \hbox{for } x,y\in D\cap B(z_0, r/2).
$$
See \eqref{e:2.19} below. This says that the boundary decay rate for $h$ at $\partial D\cap B(z_0,  r)$ is $\E_x \tau_{D_{2r}(z_0)}$.

 \end{enumerate}
\end{remark}

 Assumption (A4) in fact is a mild condition. We will show in  Theorem \ref{T:2.8}
  below that the Assumption (A4) holds for any L\'evy process
 in $\R^d$ and for any open subset $D\subset  \R^d.$
  Let $Y$ be a  L\'evy process in $\R^d$ with L\'evy triple $(A, b, \nu)$,
 where  $A=(a_{ij})_{1\leq i, j\leq d}$ is a symmetric non-negative definite constant matrix,
$b\in \R^d$ is a constant vector, and $\nu$ is a non-negative Borel measure on $\R^d\setminus \{0\}$
satisfying $\int_{\R^d} (|z|^2 \wedge 1) \nu(dz).$
Note that $\nu(x, B)=\nu(B-x)$ for each Borel set $B$ in $\R^d.$
The assumption (A3) can be restated as follows.

\smallskip

\noindent {\bf (A3')}  Either $\nu =0$ or  $\nu(\R^d\setminus \{0\})> 0$  and there exist $c>1$ and $r_0>0$ so  that for any $0<r<r_0$,
 \begin{equation}\label{e:2.4}
c^{-1}\nu(dz-x)\leq \nu(dz) \leq c\nu(dz-x)
\quad \hbox{on  } B(0, 2r)^c \hbox{ for any } x\in B(0, r).
 \end{equation}

The following follows readily from Theorem \ref{T:2.1} and Theorem \ref{T:2.8} below.

\begin{cor}\label{C2}
 Let $Y$ be a  discontinuous L\'evy process in $\R^d$ with $\nu(\R^d) > 0$.  Let $D$ be an open set. Suppose (A1)-(A2) and (A3') hold.
 Then the  scale invariant boundary Harnack principle holds for $Y$ in an open set $D$  if and only if \eqref{e:1.4a} holds.
\end{cor}

\smallskip

 Corollary \ref{C2} holds for a large class of discontinuous L\'evy processes.
 A well known result of Hartman and Wintner \cite{HW} asserts
that    a L\'evy process $Y$ has a bounded continuous density function $p(t, y-x)$ with respect to the Lebesgue measure
if its   L\'evy exponent $  \psi$ has the property that
\begin{equation}\label{e:HW}
\lim_{|\xi| \to \infty} \frac{{\rm Re}   \psi  (\xi)}{\ln (1+|\xi|)} =\infty.
\end{equation}
In this case,
for all $x\in \R^d$ and $\alpha>0,$ the resolvent operator $G^\alpha(x, \cdot)$ and the co-resolvent operator $\wh G^\alpha(x, \cdot)$ are absolutely continuous with respect to the Lebesgue measure with the potential kernel $\wh G^\alpha(x, y)=G^\alpha(y, x).$ That is, the Assumption (A1) holds for
 any L\'evy process $Y$ that satisfies condition \eqref{e:HW}.
 It is well known that any subordinate Brownian motion in $\R^d$ has a continuous transition density function  with respect to the Lebesgue measure.
 For the Assumption (A2),  Hunt's hypotheses (H)  holds for a large class of  L\'evy processes, including
 symmetric L\'evy processes and L\'evy processes whose L\'evy exponents satisfy sector condition.
We refer the reader to \cite{HS} for  a survey on   Hunt's  hypotheses (H).

\smallskip

 Next, we apply Theorem \ref{T:2.1} to study scale invariant BHP   for  subordinate Brownian motions
 or equivalently, for non-local operators $\L^\phi$ of \eqref{e:1.1a},
 on non-smooth  domains of  $\R^d$.
   Recall that a subordinator is an increasing L\'evy process starting from $0$,  which can be characterized through its Laplace exponent $\phi:$
$\E[e^{-\lambda S_t}]=e^{-t\phi(\lambda)}$, $\lambda\geq 0.$ The Laplace exponent of a subordinator is a Berstein function and  has the representation
\begin{equation}\label{e:1.12}
 \phi  (\lambda)=b\lambda +\int_0^\infty (1-e^{-\lambda t})\,\mu(dt),
\end{equation}
where  $b\geq 0$ and $\mu$ is a  measure (called the L\'evy measure of $\phi$)
 on $[0, \infty)$ satisfying $\int_0^\infty (1\wedge t)\mu(dt)<\infty$.
   With a little abuse of notation, we still use $X$ to denote the subordinate Brownian motion  associated with
   the generator  $\L^\phi:=-\phi(-\Delta)$;  that is,
   $X_t=W_{S_t},$
where $W$ is a standard Brownian motion on $\R^d$  and $S$ is a  subordinator
  with the Laplace exponent
    $\phi(\lambda)$  independent of $W$.
 For an open set $D$ of $\R^d$,  denote by $G_D(x, y)$
the Green function of $X$ in $D$ in the sense of \eqref{e:1.3'} with $m(dy)=dy$.
It is known that $G_D(x, y)$ is a positive continuous function on $D\times D$ off its diagonal.

\smallskip

  As stated in the paragraph after  Corollary \ref{C2},
any subordinate Brownian motion in $\R^d$ satisfies Assumptions {\bf (A1)} and {\bf (A2)}.
   We will show   (in the paragraph containing  \eqref{e:mu1})
 that
   if there is a density $\mu(t)$ of the L\'evy measure $\mu (dt)$ of the subordinator $S$ with respect to the Lebesgue measure on $[0, \infty)$  satisfying the condition \eqref{e:mu2} below,
 then  the condition {\bf (A3')} holds.
  We get the following theorem from Corollary \ref{C2}.

\begin{thm}\label{T0}
 Suppose
   that  $\phi$ is the Laplace exponent  of a discontinuous subordinator $S$  whose L\'evy measure $\mu (dt)$ has a density $\mu(t)$
     with respect to the Lebesgue measure on $[0, \infty)$
 which has the property that there is a constant $c>0$ so that
\begin{equation}\label{e:mu2}
\mu (t) \leq  c \mu (2t)  \    \hbox{ for } t\in (0, 8)  \quad \hbox{and} \quad
\mu(t)\leq c\mu(t+1) \    \hbox{ for } t\in (1, \infty).
\end{equation}
Then the  scale invariant boundary Harnack principle holds for $\L^\phi $ in an open set $D$ in $\R^d$ if and only if \eqref{e:1.4a} holds.
\end{thm}

\begin{remark}\label{R:1.5} \rm

\begin{enumerate}
  \item[(i)]    Examples of Laplace exponent $\phi$   of subordinator that satisfy the conditions of Theorem \ref{T0} include:
\begin{itemize}
\item[(a)]  $\phi (\lambda) =b \lambda+  \lambda^\alpha$ for any $\alpha \in (0, 1)$ and $b\geq 0$;
\item[(b)]  $\phi (\lambda) =b \lambda+  \lambda^\alpha + \lambda^\beta $ for any $\alpha, \beta \in (0, 1)$  and $b\geq 0$;
\item[(c)]  $\phi (\lambda) =b \lambda+ (\lambda+m)^\alpha - m^\alpha$  for any $\alpha \in (0, 1)$,  $m>0$  and $b\geq 0$;
\item[(d)]  $\phi (\lambda) =b \lambda+ \log (\lambda+1) $ for any $b\geq 0$.
\end{itemize}
\smallskip
All these functions are complete Bernstein functions.  Note that if the Laplace exponent $\phi$ is a complete Bernstein function,
then the second inequality in \eqref{e:mu2} is automatically satisfied; see Remark \ref{R:3.5}.

\smallskip
\item[(ii)]   Theorem \ref{T0}  is in sharp contrast with the Brownian motion case.  It  is known that the comparability condition
    \eqref{e:1.4a}  does not hold  in general on bounded  Lipschitz domains for the Laplacian operator or Brownian motion.
In fact, it is shown in Atar, Athreya and Chen \cite{AAC}  that
 the   Green function and the mean exit time  of Brownian motion
 are comparable in any  bounded Lipschitz domain  that satisfies a uniform  interior cone condition with
  the common cone angle $\theta\in (\cos^{-1}(1/\sqrt d), \pi)$,    but
   are not comparable  in a truncated cone with angle  $ \theta\in (0,   \cos^{-1}(1/\sqrt d)).$
  On the other hand, the scale invariant BHP for Brownian motion  holds  on any bounded Lipschitz domains.
	 Our results show that the scale invariant BHP is unstable   on non-smooth domains  for non-local operators that have diffusive components
	in the sense that among the family of non-local operators $\sL^{a, b}:=a \Delta +b \Delta^{\alpha/2}$ with $0<\alpha<2$
	and non-negative $a$ and $b$ so that  $a+b>0$, the Euclidean domains over which the scale invariant BHP for $\sL^{a, b}$ holds has  a phase transition at $a=0$ and at $b=0$.   See Theorem \ref{T:1.10} below for further information on this.
\end{enumerate}
\end{remark}

\smallskip

In the following,
we will  give geometric conditions for  Lipschitz open sets  in $\R^d$ with $ d\geq 2$
on which  the scale invariant BHP for $\L^\phi$ holds  when  the Laplace exponent $\phi$ has a positive drift $b$ and a non-trivial  L\'evy measure $\mu$.

\begin{defn}\label{D1} \rm
   \begin{enumerate}
    \item[\rm (i)]
An open set $D$ in $\R^d$  is said to be Lipschitz  if there are positive constants
$R_0 $ and  $\Lambda_0 $ so that for every $z\in\partial D,$ there is  an orthonormal coordinate system
$CS_z: y=(y_1, \cdots, y_{d-1}, y_d)=:(\wt {y}, y_d)\in \R^{d-1}\times \R$ with its origin at $z$
and    a Lipschitz function
$\varphi=\varphi_z: \R^{d-1}\rightarrow \R$ with Lipschitz constant $\Lambda_0$
 such that
$$
B(z, R_0)\cap D= \left\{y=(\wt {y}, y_d)\in B(0,R_0) \hbox{ in } CS_z: y_d>\varphi(\wt {y}) \right\}.
$$
The pair $(R_0, \Lambda_0)$ is called the characteristics of the Lipschitz open set $D.$

   \item[\rm (ii)]  We say an open set $D\subset \R^d$ is $C^1$
   if it is Lipschitz with the local boundary functions $\varphi_z$ in the above definition being $C^1$
   and $\nabla \varphi_z$ having uniform modulus of continuity.
    \end{enumerate}
\end{defn}

For   $\theta\in (0, \pi)$ and $r>0,$ define
\begin{equation}\label{e:cone}
\Gamma_\theta(r):=\{x=(x_1, \ldots, x_d) \in \R^d: |x|<r \hbox{ and } x_d>|x|\cos\theta\},
\end{equation}
 which is a   truncated  circular cone with apex $2\theta$ and axis along the positive  $x_d$-direction.
 We    call $\theta$ the angle of the cone $\Gamma_\theta (r)$, and
    will denote  $\Gamma_\theta(1)$ by $\Gamma_\theta$.

\begin{defn} \rm
We say that a domain $D$ satisfies the  interior cone
condition with common angle $\theta\in (0, \pi)$  if
there is $a_0>0$ such that
for every point $z\in \partial D,$
there is a cone $\Gamma_\theta(z, a_0)\subset D$ with vertex at $z$ that is conjugate to $\Gamma_\theta(a_0);$
that is, $\Gamma_\theta(z, a_0)$ is the cone with vertex at $z$ that is obtained from $\Gamma_\theta(a_0)$ through parallel translation and rotation.
\end{defn}

It is easy to see that a Lipschitz domain naturally satisfies the interior cone condition   for some common angle $\theta >0$.

\begin{thm}\label{T:1.10}
 Let $d\geq 2$.
 Suppose that $\phi$ is the Laplace exponent of a discontinuous subordinator that has a positive drift $b>0$ and satisfies the condition
 of Theorem \ref{T0}.
Let  $D$ be a Lipschitz domain in $\R^d$ with characteristics $(R_0, \Lambda_0)$ that satisfies the interior cone condition with common angle $\theta\in (\cos^{-1}(1/\sqrt d), \pi)$.
Then there exists $ C= C(d, \phi, R_0,  \Lambda_0,  \theta)>0$  such that for any $z_0\in\partial D, r\in (0, R_0/2)$ and any two  nonnegative harmonic functions $h_1$ and $h_2$ with respect to  the subordinate Brownian motion  $X$ with generator   $\L^\phi$ in
 $D\cap B(z_0, 2r)$  vanishing continuously on $D^c\cap B(z_0, 2r),$
$$
\dfrac{h_1(x)}{h_1(y)}\leq C \dfrac{h_2(x)}{h_2(y)} \quad \hbox{for any }
x,y\in D\cap B(z_0, r).
$$
If $D$ is a truncated cone $\Gamma_\theta$ with  angle $\theta\in (0, \cos^{-1}(1/\sqrt d)],$
then the  above scale invariant BHP fails to hold  for  $\sL^\phi$.
\end{thm}

 Theorem \ref{T:1.10} will be established through Theorem \ref{T4} and Theorem \ref{T5}.
As noted earlier,
a Lipschitz domain naturally satisfies the interior cone condition with a common angle.
 Theorem \ref{T:1.10} states that the scale invariant BHP holds for $\L^\phi$ in  Lipschitz domains in $\R^d$ satisfying the interior cone condition with common angle $\theta\in (\cos^{-1}(1/\sqrt d), \pi)$. This includes as special case the Lipschitz domains with the Lipschitz constant  strictly less than $1/\sqrt{d-1}$.
  On the other hand, any cone with an obtuse angle is a Lipschitz domain
 that satisfies interior cone
 condition with common angle $\theta\in (\cos^{-1}(1/\sqrt d), \pi)$ but whose Lipschitz constant can be arbitrarily large.
  Theorem  \ref{T:1.10}
   is sharp in the sense that the  scale invariant BHP
    fails in any truncated cone  with angle $\theta \leq \cos^{-1}(1/\sqrt d)$,
   which is
  a Lipschitz domain  whose Lipschitz constant   is larger than or equal to $1/\sqrt{d-1}.$
  Theorem \ref{T:1.10} in particular also implies that the scale invariant boundary Harnack inequality
 holds for  $\sL^\phi$ in any $C^1$ domains; see Corollary \ref{C1}.
 The scale invariant BHP  also holds for any piecewise   $C^1$ smooth domains for which there is some
$\theta > \cos^{-1}(1/\sqrt d)$ so that the (truncated) interior $\theta$-cone condition is satisfied at each non-smooth points on the boundary.

  Theorem \ref{T:1.10} is also in stark contrast with the case of  the fractional Laplacian  in addition to the case of the  Laplacian.
  It is known that the scale invariant  BHP holds for the fractional Laplacian $\Delta^{\alpha/2}$ with $\alpha \in (0, 2)$ on
  any   open set  (see \cite{BKK}).

A probabilistic ``box" method  is developed by Bass and Burdzy  \cite{BB1}  to give a new proof of scale invariant BHP for Laplacian in Lipschitz domains.
It is further developed  to establish BHP in twisted H\"older domains \cite{BB2}, and later adapted by Aikawa \cite{Ai} analytically for uniform domains.
In this paper,  we use the box method  in Lemma \ref{L:2.4} to establish the comparability of
   the Green function and the mean exit time
    for discontinuous
 subordinate Brownian motions having Gaussian components
  in Lipschitz domains that satisfy
   the uniformly interior cone condition with common angle $\theta\in (\cos^{-1}(1/\sqrt d), \pi)$.
The counterexample that the scale invariant BHP fails in  a truncated cone with angle $\theta\in (0, \cos^{-1}(1/\sqrt d)]$   relies on the key observation  that
 would the condition  \eqref{e:1.4a}  hold  in a truncated cone,  the Green function and the mean exit time of Brownian motion in the truncated cone would be comparable along the axis approaching the vertex of the cone. The latter is shown to fail in \cite{AAC}   if the angle $\theta$ of the cone  is strictly smaller than
$\cos^{-1}(1/\sqrt d).$
 We further show in this paper that this also fails when  $\theta=\cos^{-1}(1/\sqrt d).$

\smallskip

The remainder of this paper is organized as follows.  In Section 2,
 we give a representation of non-negative   regular harmonic functions of  a Hunt process $X$ in bounded open sets in $\sX$ that vanish on a portion of the boundary in terms of Green function. We use it to prove Theorem \ref{T:2.1}.
We then show in Theorem \ref{T:2.8} that the Assumption (A4) holds for any L\'evy process
 in $\R^d$ and for any open subset $D\subset  \R^d$, from which and Theorem \ref{T:2.1} we immediately get Corollary \ref{C2}.
  In Section 3,  we  establish the scale invariant BHP  for discontinuous subordinate Brownian motions having
  Gaussian components on Lipschitz domains that satisfy the interior cone condition with common angle $\theta\in (\cos^{-1}(1/\sqrt d), \pi)$ in Theorem \ref{T4}.
  As a consequence, we  deduce in Corollary \ref{C1}   that the scale invariant BHP holds for any   $C^1$ domains with uniform modulo of continuity. This in particular extends the scale invariant BHP results \cite{CKSV, KSV2}
  that are established using a test function method developed in \cite{BBC}.
   The counterexample that the scale invariant BHP fails in  a truncated cone with angle $\theta\in (0, \cos^{-1}(1/\sqrt d)]$ is proved in Theorem \ref{T5}.
     The conclusion of Theorem \ref{T:1.10} is a combination of Theorems \ref{T4} and \ref{T5}.
   In the Appendix,  we provide the details of the proof for Pruitt's result of  Theorem \ref{T:0.3}  for   L\'evy processes in $\R^d$.

\section{Scale invariant BHP on metric measure spaces}

In this section, we shall establish the scale invariant BHP of  Theorem \ref{T:2.1}
 in a general framework for discontinuous Hunt processes on   locally compact separable metric spaces
 and present a proof of Theorem \ref{T:2.8}  that the Assumption (A4) holds for any L\'evy process
 in $\R^d$ and for any open subset $D\subset  \R^d$.

Let $(\sX, d, m)$ be a locally compact separable metric space with a Radon measure $m.$
Let $(X_t, \zeta, \mathcal{F}_t, \P_x)$ be a  discontinuous Hunt process  with state space  $\sX$.
 We assume that $X$ is  in weak duality to another Hunt process $\wh X$ with respect to a  $\sigma$-finite measure $m$ on $\sX$ that has full support in the  sense of \eqref{e:1.1'}.
 The duality assumption is assumed here as we need to use  some potential theoretic results from
   \cite{BG, ChungW}   in our proof of Theorem \ref{T:2.1}.

 Throughout this section,
we always assume $X$ is a Hunt process having the L\'evy system $(N(x, dy), dt).$
 For a locally compact open set $U$ of $\sX$ with $\overline U\subsetneq \sX$,  denote by $G_U(x, y)$
the Green function of $X$ in $U.$
  As mentioned in the Introduction, the Green function $G_U(x, y)$ exists uniquely for each  $(x, y)\in U\times U$ under Assumption (A1).
 For an open set $U\subset \sX,$ denote by $X^U$ the part process of $X$ killed upon exiting $U$ .
 For each open subset $B\subset \sX,$  $\tau_B:=\inf\{t \geq 0: X_t \notin B\}$ denotes the first exiting time from $B$ by  the process $X$.

\begin{defn}\label{D:2.1} \rm
Suppose $U$ is an open set in $\sX.$ A real-valued function $u$ defined on $\sX$ is said to be
\begin{itemize}
\item[\rm (i)]
 {\it  harmonic in   $U$ with respect to $X$}
if for every open set $B$  whose closure is a compact subset of $U,$
   $\E_x |u(X_{\tau_B})|<\infty$ and
 $u(x)=\E_x u(X_{\tau_B})$ for  each  $x\in B.$
  In particular, we say $u$  is   regular harmonic in   $U$ with respect to $X$
if    $\E_x |u(X_{\tau_U})|<\infty$ and
$u(x)=\E_x u(X_{\tau_U})$  for each  $x\in U.$

\item[\rm (ii)]
 {\it $\alpha$-excessive with respect to $X^U$ with $\alpha\geq 0$}
if $u$ is non-negative,
$ u(x)\geq e^{-\alpha t}\E_x u(X^U_t)$ and $u(x)=\lim_{t\downarrow 0} e^{-\alpha t}\E_x  \left[ u(X^U_t)\right]$
 for every $t>0$ and   $x\in U.$
 If $\alpha=0,$ we say $u$ is an excessive function with respect to $X^U.$

\end{itemize}
\end{defn}

\begin{lem}\label{L:2.2}
 Suppose Assumption (A3) holds.
Let $D$ be an open set in $\sX$ and $r_0$ be the constant in assumption (A3).
There exists $C\geq 1$ such that for any $z_0\in \partial D, r\in (0, r_0)$ and $x\in D\cap B(z_0, r),$
\begin{eqnarray*}
C^{-1}  N(z_0, B^c(z_0, 2r)) \, \E_x \tau_{D\cap B(z_0,r)}
&  \leq & \P_x(X_{\tau_{D\cap B(z_0, r)}}\in B^c(z_0, 2r))\\
& \leq & C N(z_0, B^c(z_0, 2r)) \, \E_x \tau_{D\cap B(z_0, r)}.
\end{eqnarray*}

\end{lem}

\pf Fix $z_0\in \partial D.$ For the simplicity of notation, for each $r>0,$ let $D_r(z_0):=D\cap B(z_0, r)$ and $B_r(z_0):=B(z_0, r).$
 By the L\'evy system formula \eqref{e:2.1'},
\begin{eqnarray*}
&&\P_x(X_{\tau_{D_r(z_0)}}\in B^c_{2r}(z_0))
=\E_x\sum_{s\leq \tau_{D_r(z_0)}}{\mathbbm 1}_{\{X_{s-}\in D_r(z_0), X_s\in B^c_{2r}(z_0)\}}\\
&=&\E_x\int_0^{\tau_{D_r(z_0)}}\int_{B^c_{2r}(z_0)}N(X_{s}, dz)\,ds
\asymp  N(z_0, B^c_{2r}(z_0))\, \E_x \tau_{D_r(z_0)} ,
\end{eqnarray*}
where the comparability in the last line is due to Assumption (A3).
\qed

\smallskip

\begin{lem}\label{L:2.4n}
 Suppose  Assumption (A1)  holds.
 Let  $D$ be an open subset of  $\sX$ so that   $D^c$ has a non-empty interior.
Then there is a strictly positive function $g_0$ on $D$ such that  $g_0\in L^1(D; m)$ and
 $0<G_D g_0 (x) \leq 1$  for every $x\in D.$
 Consequently, for every $x\in D$, $G_D(x, y)<\infty$ for $m$-a.e. $y\in D$.
\end{lem}

\pf Denote the transition semigroup of $X$ by $\{P_t; t\geq 0\}$.
 Recall that a measure $\mu$ on $(\sX, \B)$ is excessive for   the Hunt process $X$  if $\mu P_t\leq \mu$ for evert $t>0.$
 From the weak duality assumption, the reference measure $m(dx)$  is clearly an excessive measure of $X$.
 Denote by $m_D(dx)$ the restriction of $m(dx)$  on $D.$  As the part process $X^D$
 of $X$ killed upon leaving $D$ is in weak duality of $\widehat X^D$ with respect to $m_D$, the measure $m_D$ is an excessive
 measure for  the Hunt process $X^D$.
 Fix a positive function $g$ in $D$ such that $\int_D g(y) m (dy)<\infty.$
 Define
 $$m_d:=m_D(\cdot \cap \{x: G_D g(x)<\infty\}), \quad m_c:=m_D(\cdot\cap \{x: G_D g(x)=\infty\}).$$
 Then $m_D=m_d+m_c.$
 By \cite[Theorem 4.3]{FM}, the measures $m_d$ and $m_c$ are  excessive measures with respect to $X^D$ in $D$.
In particular,  $m_c$ is sub-invariant with respect to the resolvents $\{G^\alpha_D: \alpha>0\}$ in the sense that
$m_c\circ (\alpha G^\alpha_D)\leq m_c$ for every $\alpha >0$.
Recall that  $G^\alpha_D:=\int_0^\infty e^{-\alpha t} P^D_t dt$, where $\{P^D_t; t\geq 0\}$
is the transition semigroup of the part process $X^D.$

In what follows, we show that $m_c(D)=0$ by contradiction. Suppose  $m_c(D)>0$.
 Since $G_D g(x)=\infty$ holds $m_c$-a.e. on $D$,  it follows from the equivalence of (iii) and (v) in \cite[Proposition 2.5]{BCRo} applied to
  the measure $m_c$ on $D$ that  there exists a Borel set $N\subset D$ with $m_c(N)=0$ so that
 $ \alpha G^\alpha_D 1(x)=1$   for every $ x\in D\setminus N$ and every   $ \alpha>0$.
Observe that
 $\alpha G^\alpha_D 1(x) = {\mathbb E}_x \int_0^{\tau_D} \alpha e^{-\alpha t} dt = 1-\E_x e^{-\alpha \tau_D}$
 for every $ x\in D.$
 Thus we have
 \begin{equation}\label{e:exit}
 {\mathbb P}_x (\tau_D=\infty)=1
 \quad \hbox{for every } x\in D\setminus N.
 \end{equation}
 On the other hand,  since ${\bar D}^c$ is non-empty,  we have by the irreducibility of $X$ that
$\P_x (\tau_D<\infty)>0$ for every $x\in D.$
 This contradicts to
\eqref{e:exit} as  $m(D\setminus N)\geq m_c(D\setminus N)>0$.
This contradiction shows that $m_c= 0$ and hence $m=m_d$, which proves that
$G_D g < \infty$ $m$-a.e. on $D.$

 For each integer $n\geq 1$, let $A_n:= \{x\in D: G_D g (x)\leq n\}$. Then $m ( D\setminus \cup_{n=1}^\infty A_n)=0$.
 Since for each $x\in D,$ $G_D(x, dy)  = G_D(x, y) m_D (dy)$ is  absolutely continuous with respect to $m(dy)$ under Assumption (A1), we have
\begin{equation}\label{e:0.2}
 G_D(x,  D\setminus \cup_{n=1}^\infty A_n )    =0  \quad \hbox{for every } x\in D.
\end{equation}
Note that for $x\in D$ and every $t>0$,
$P^D_t G_D g(x) = \E_x \int_t^{\infty} g(X^D_s )ds$,
which increases to $G_Dg(x) =\E_x \int_0^\infty g(X^D_s )ds$ as $t \downarrow 0$.
So  $G_D g$ is an excessive function of $X^D$.
By the strong Markov property of $X^D$ and \cite[Theorem II.2.12]{BG},
for every $x\in D$,
\begin{eqnarray*}
G_D ( {\mathbbm 1}_{A_n} g) (x) &=& \E_x \int_{T_{A_n}}^\infty ({\mathbbm 1}_{A_n} g)(X^D_s )ds =
\E_x   [ G_D  ({\mathbbm 1}_{A_n}g)   (X^D_{T_{A_n}})]\\
& \leq & \E_x   [ (G_D  g)  (X^D_{T_{A_n}})]  \leq n,
\end{eqnarray*}
where $T_{A_n}:=\inf\{t>0: X_t\in A_n\}.$
Define $g_0:= \sum_{n=1}^\infty {\mathbbm 1}_{A_n}(n2^n)^{-1}  g + {\mathbbm 1}_{D\setminus \cup_{k=1}^\infty A_k}$.
Then $g_0>0$ on $D$ with $g_0 \in L^1(D; m)$ and
$$
G_D g_0 (x)  = \sum_{n=1}^\infty \frac{1}{n 2^n} G_D ({\mathbbm 1}_{A_n} g) (x)\leq  \sum_{n=1}^\infty \frac{n}{n 2^n}=1
\quad \hbox{for every }  x\in D,
$$
  where the first equality holds due to \eqref{e:0.2}.
On the other hand,
since $g_0>0$ on $D$, $\int_0^{\tau_D} g_0(X_s) ds >0$ ${\mathbb P}_x$-a.s and hence
$G_D g_0(x)={\mathbb E}_x \int_0^{\tau_D} g_0(X_s) ds >0$ for every $x\in D$.
This  establishes that $0<G_D g_0 (x) \leq 1$ for every $x\in D$.  Since $G_D g_0(x)= \int_D G_D(x, y) g_0(y) m(dy)$,
we have for each $x\in D$, $G_D(x, y)<\infty$ for $m$-a.e. $y\in D$.
\qed

\smallskip

  Lemma \ref{L:2.4n} shows that the semigroup of $X^D$  is transient  in the sense of
  \cite[Definition on page 86]{ChungW} that there exists a strictly positive function $g_0$ on $D$ such that  $0<G_D g_0 (x) <\infty$  for every $x\in D.$

\smallskip

\begin{prp}\label{P:2.3}
Suppose Assumptions (A1)-(A2)  hold.
 Let $D$ be an open subset of  $\sX$ so that   $D^c$ has a non-empty interior.  Let  $U$ and $V$ be relatively compact  open sets in $\sX$ with $  \overline{V}  \subseteq U$ and $D\cap V\neq \varnothing.$
  Suppose $h$ is a nonnegative    regular  harmonic function with respect to $X$ in $D\cap U$ and   vanishes
  on $(D^c)^r\cap U$.
 Then
 there exist a sequence of Radon measures $\{\mu_n; n\geq 1\}$ on $D\setminus V$
  so that     $G_{D\cap U}\mu_n(x):= \int_{D\setminus V} G_{D\cap U}(x, y) \mu_n (dy)$   is uniformly bounded on $D\cap V$ and increasing   in $n\geq 1$,  and for $x\in   D\cap V,$
  \begin{eqnarray}\label{e:2.5a}
h(x)& = &  \lim_{n\to \infty} \int_{D\setminus  V} G_{D\cap U}(x, y)\mu_n(dy) \nonumber\\
&&  + \int_{D\cap U} G_{D\cap U}(x,y)\int_{U^c}h(z) N(y, dz)\,m(dy).
\end{eqnarray}

\end{prp}

\pf Let $h$ be a nonnegative regular harmonic function with respect to $X$ in $D\cap U$ vanishing  on
 $(D^c)^r\cap U.$
Then for each $x\in D\cap U,$
$$\begin{aligned}
h(x)&=\E_x h(X_{\tau_{D\cap U}})\\
&=\E_x  \left[h(X_{\tau_{D\cap U}}), X_{\tau_{D\cap U-}}=X_{\tau_{D\cap U}} \right]
+\E_x \left[ h(X_{\tau_{D\cap U}}), X_{\tau_{D\cap U}-}\neq X_{\tau_{D\cap U}} \right] .
\end{aligned}$$
Let
\begin{equation}\label{e:2.2}
h_1(x):=
\begin{cases}
h(x)-\E_x  \left[h(X_{\tau_{D\cap U}}); X_{\tau_{D\cap U}-}\neq X_{\tau_{D\cap U}} \right]
\quad & \hbox{for } x\in D\cap U, \\
0  & \hbox{for } x\in \sX \setminus (D\cap U) .
\end{cases}
\end{equation}
For any open set $B$ with $\overline B\Subset   D\cap U$ and $x\in B$, by the strong Markov property of $X$,
\begin{eqnarray}\label{e:2.3'}
 && \E_x  \left[ h_1(X_{\tau_{B}}); \tau_{B}<\tau_{D\cap U} \right]  \nonumber \\
&=& \E_x  \left[h (X_{\tau_{B}});  \tau_{B} < \tau_{D\cap U} \right]
-\E_x\big[ \E_{X_{\tau_{B}}}
\big[ h(X_{\tau_{D\cap U}}); X_{\tau_{D\cap U}-}\neq X_{\tau_{D\cap U}} \big]; \tau_{B}< \tau_{D\cap U}\big]  \nonumber \\
&=&\E_x  \left[h (X_{\tau_{B}});  \tau_{B} < \tau_{D\cap U} \right]
- \E_x \left[  h(X_{\tau_{D\cap U}});    \tau_{B}< \tau_{D\cap U} \hbox{ and }  X_{\tau_{D\cap U}-}\neq X_{\tau_{D\cap U}}   \right]
 \nonumber \\
&=& \E_x  \left[h (X_{\tau_{B}})\right] -  \E_x  \left[h (X_{\tau_{B}}) ;  \tau_{B} = \tau_{D\cap U} \right] \nonumber\\
&&-\E_x \left[  h(X_{\tau_{D\cap U}});    \tau_{B}< \tau_{D\cap U} \hbox{ and }  X_{\tau_{D\cap U}-}\neq X_{\tau_{D\cap U}}   \right] \nonumber \\
&=&  h(x) -  \E_x \left[  h(X_{\tau_{D\cap U}});    X_{\tau_{D\cap U}-}\neq X_{\tau_{D\cap U}}   \right]  \nonumber \\
&=&  h_1(x),
\end{eqnarray}
where  the second to the last equality is due to the fact that  $X_{\tau_{D\cap U}-}\neq X_{\tau_{D\cap U}}  $ on $\{  \tau_{B}= \tau_{D\cap U} \}$.
This shows that
$h_1$ is a non-negative  harmonic function in $D\cap U$ with respect to   the part process $X^{D\cap U}$.

In the following, we show  that $h_1$ is  regular  harmonic  in $D\cap V$ with respect to   the part process $X^{D\cap U}.$
In fact, by a similar argument of \eqref{e:2.3'}, for $x\in D\cap V,$
\begin{eqnarray}\label{e:2.4n}
 &&  \E_x  \left[ h_1(X_{\tau_{D\cap V}}); \tau_{D\cap V}<\tau_{D\cap U} \right]   \nonumber \\
 &=&\E_x  \left[h (X_{\tau_{D\cap V}});  \tau_{D\cap V} < \tau_{D\cap U} \right]   \nonumber \\
 &&-\E_x \left[\E_{X_{\tau_{D\cap V}}}
\left[ h(X_{\tau_{D\cap U}}); X_{\tau_{D\cap U}-}\neq X_{\tau_{D\cap U}} \right]; \tau_{D\cap V}< \tau_{D\cap U}\right]  \nonumber \\
&=&\E_x  \left[h (X_{\tau_{D\cap V}})\right] -  \E_x  \left[h (X_{\tau_{D\cap V}}) ;  \tau_{D\cap V} = \tau_{D\cap U} \right]   \nonumber \\
&& -\E_x \left[  h(X_{\tau_{D\cap U}});    \tau_{D\cap V}< \tau_{D\cap U} \hbox{ and }  X_{\tau_{D\cap U}-}\neq X_{\tau_{D\cap U}}  \right].
\end{eqnarray}

 It is well known  (see e.g. \cite[Proposition (3.3) on p.80]{BG})   that  $ D^c\setminus (D^c)^r$ is semi-polar.
Hence it is polar by  Assumption (A2); that is,
$\P_y(X_t  \mbox{ ever hits } D^c\setminus (D^c)^r)=0$ for all $y.$
Since $h$  is zero at those points on $(\partial D)^r\cap V\subset (D^c)^r\cap U$, then
\begin{equation}\label{e:2.8n}
\E_x  \left[h (X_{\tau_{D\cap V}}) ;  \tau_{D\cap V} = \tau_{D\cap U}, X_{\tau_{D\cap V}}\in \partial D\cap V \right]=0, \quad x\in D\cap V.
\end{equation}
Hence, for $x\in D\cap V,$
\begin{eqnarray}\label{e:2.5n}
&&  \E_x  \left[h (X_{\tau_{D\cap V}}) ;  \tau_{D\cap V} = \tau_{D\cap U} \right] \nonumber \\
&=& \E_x  \left[h (X_{\tau_{D \cap V}}) ;  \tau_{D\cap V}= \tau_{D\cap U}, X_{\tau_{D\cap V}-}= X_{\tau_{D\cap V}} \right]
 \nonumber \\
&& \quad+\E_x  \left[h (X_{\tau_{D \cap V}}) ;  \tau_{D\cap V}= \tau_{D\cap U}, X_{\tau_{D\cap V}-}\neq X_{\tau_{D\cap V}} \right]
 \nonumber \\
&=& \E_x  \left[h (X_{\tau_{D \cap V}}) ;  \tau_{D\cap V}= \tau_{D\cap U}, X_{\tau_{D\cap V}}\in \partial D\cap V \right]
 \nonumber \\
&& \quad+\E_x  \left[h (X_{\tau_{D \cap V}}) ;  \tau_{D\cap V}= \tau_{D\cap U}, X_{\tau_{D\cap V}-}\neq X_{\tau_{D\cap V}} \right]  \nonumber \\
&=& \E_x  \left[h (X_{\tau_{D\cap U}}) ;  \tau_{D\cap V} = \tau_{D\cap U}, X_{\tau_{D\cap U}-}\neq X_{\tau_{D\cap U}} \right].
\end{eqnarray}
As by the strong Markov property, $h(x)=\E_x  \left[h (X_{\tau_{D\cap V}})\right]$ for $x\in D\cap V.$
Thus by \eqref{e:2.4n} together with \eqref{e:2.5n}, we have for every $x\in D\cap V$,
\begin{equation}\label{e:2.6}
\begin{aligned}
&\quad \E_x  \left[ h_1(X_{\tau_{D\cap V}}); \tau_{D\cap V}<\tau_{D\cap U} \right]\\
&= h(x) - \E_x \left[  h(X_{\tau_{D\cap U}});    X_{\tau_{D\cap U}-}\neq X_{\tau_{D\cap U}}   \right]
=h_1(x).
\end{aligned}
\end{equation}
Hence,  $h_1$ is  regular  harmonic  in $D\cap V$ with respect to   the part process $X^{D\cap U}.$

  Note that $h_1$ is a nonnegative harmonic function of $X^{D\cap U}$. So it is  an excessive function of $X^{D\cap U}$ (see e.g. \cite[Proposition 3.4]{KW}). By applying Lemma \ref{L:2.4n} in $D\cap U$ in place of $D$  and \cite[Proposition 10 in Section 3.2]{ChungW},  there exist a sequence of bounded positive functions $\{g_n\}_{n\geq 1}$ on $D\cap U$ such that
\begin{equation}\label{e:2.10}
h_1(x)=\lim_{n\rightarrow\infty}\uparrow G_{D\cap U} g_n(x)=\lim_{n\rightarrow\infty}\uparrow G_{D\cap U} g_n(x), \quad x\in D\cap U.
\end{equation}
Consequently, it follows from \eqref{e:2.6} and \eqref{e:2.10},  the monotone convergence theorem and Fubini's theorem that  for every  $x\in D\cap V$,
\begin{eqnarray*}
h_1(x)  &=& \E_x h_1(X^{D\cap U}_{\tau_{D\cap V}})\\
&= &  \E_x \left[   \lim_{n\to \infty}   G_{D\cap U} g_n (X^{D\cap U}_{\tau_{D\cap V}}) \right]  =  \lim_{n\to \infty}   \E_x \left[     G_{D\cap U} g_n (X^{D\cap U}_{\tau_{D\cap V}}) \right]   \\
&=& \lim_{n\to \infty}   \int_{D\cap U} \E_x      G_{D\cap U} (X^{D\cap U}_{\tau_{D\cap V}}, y)g_n(y)\,m(dy)\\
&=& \lim_{n\to \infty}   \int_{D\cap U} \E_x      G_{D\cap U} (X^{D\cap U}_{T_{D\setminus V}}, y)g_n(y)\,m(dy),
\end{eqnarray*}
where $T_{D\setminus V}:=\inf\{t>0: X_t\in D\setminus V\}.$
Recall that $X$ is   in weak duality to   the Hunt process $\wh X$ with respect to  the Radon measure $m$ on $\sX.$
Let $\wh X^{D\cap U}$ be the part process of $\wh X$ killed upon exiting $D\cap U$. For a Borel set $A\subset \sX$,  define $\wh T_A:=\inf\{t>0: \wh X_t\in A\}.$ Denote by $\wh \P_y$ and $\wh\E_y$  the probability and expectation of $\wh X$ starting from $y.$
Define $\wh P^{D\cap U}_{D\setminus V}(y, dz):=\wh\P_y(\wh X^{D\cap U}_{\wh T_{D\setminus V}}\in dz).$
By virtue of the  formula in \cite[Theorem 1.16 in Section VI]{BG}, we have
$$\E_x G_{D\cap U} (X^{D\cap U}_{T_{D\setminus V}}, y)=\wh\E_y G_{D\cap U} (x, \wh X^{D\cap U}_{\wh T_{D\setminus V}}), \quad x, y\in D\cap U.$$
Hence,  for every  $x\in D\cap V$,
\begin{eqnarray*}
h_1(x)
&=& \lim_{n\to \infty}   \int_{D\cap U} \wh \E_y  G_{D\cap U} (x, \wh X^{D\cap U}_{\wh T_{D\setminus V}})g_n(y)\,m(dy) \\
&=& \lim_{n\to \infty}   \int_{D\cap U} \int_{D\setminus V}  G_{D\cap U} (x, z) \wh P^{D\cap U}_{D\setminus V}( y, dz)g_n(y)\,m(dy) ,
\end{eqnarray*}
where in the last equality, we used the fact that $\wh P^{D\cap U}_{D\setminus V}( y, dz)$ is supported on $D\setminus V$  (see e.g. \cite[Theorem 2 in Section 3.4]{ChungW}).
For each $n\geq 1,$ let $\mu_n(dz):=\int_{D\cap U}\wh P^{D\cap U}_{D\setminus V}(y, dz)g_n(y)\,m(dy).$ Then $\{\mu_n(dz)\}_{n\geq 1}$ is a sequence of Radon measures on $D\setminus V$ such that
\begin{equation}\label{e:2.12}
h_1(x)=\lim_{n\to \infty}\int_{D\setminus V} G_{D\cap U}(x, z)\mu_n(dz), \quad x\in D\cap V.
\end{equation}

By the L\'evy system formula \eqref{e:2.1'}, Assumption (A2) and that $h$ vanishes on $(D^c)^r\cap U,$ we have for $x\in D\cap V,$
\begin{eqnarray*}
&& \E_x  \left[h(X_{\tau_{D\cap U}}); X_{\tau_{D\cap U}-}\neq X_{\tau_{D\cap U}} \right] \\
&=&\E_x\sum_{s\leq \tau_{D\cap U}} h(X_s){\mathbbm 1}_{\{X_{s-}\neq X_s; X_s\in U^c\}}\\
&=&\E_x\int_0^{\tau_{D\cap U}}\int_{U^c} h(z)N(X_s, dz)\,ds\\
&=& \int_{D\cap U} G_{D\cap U}(x,y)\int_{U^c}h(z) N(y, dz)\,m(dy).
\end{eqnarray*}
This together with  \eqref{e:2.2} and \eqref{e:2.12} yields the desired conclusion \eqref{e:2.5a}.
 \qed

\smallskip

\noindent{\bf Proof of Theorem \ref{T:2.1}.}
Let $D$ be an open subset of  $\sX$ so that   $D^c$ has a non-empty interior.
Let $z_0$ be a boundary point  in $\partial D.$
For the simplicity of notation, for each $r>0,$ we denote by  $D_r(z_0):=D\cap B(z_0, r)$ and $B_r(z_0):=B(z_0, r).$

\smallskip

(i)   We   show that the condition \eqref{e:1.4a} is  sufficient for the scale invariant BHP for $X$ in $D$
 under assumptions (A1)-(A4).
Let   $\delta, r_0$ and $\delta_0$ be the constants in the condition \eqref{e:1.4a}, (A3) and (A4).
Let $\delta_1:=(\delta\wedge \delta_0\wedge r_0)/4.$
Let $r\in (0,  \delta_1), z_0\in \partial D$ and
let $h$ be a nonnegative regular harmonic function in $D_{4r}(z_0)$ vanishing on   $(D^c)^r\cap B_{4r}(z_0).$
By taking $U=B_{4r}(z_0)$ and $V=B_{3r}(z_0)$  in Proposition \ref{P:2.3},
there is a sequence of  Radon measures $\{\mu_n; n\geq 1\}$ on $D_{4r}(z_0)\setminus D_{3r}(z_0)$
so that for  $x\in D_{3r}(z_0),$
  \begin{eqnarray*}
  h(x)&=&  \lim_{n\to \infty} \int_{D_{4r}(z_0)\setminus D_{3r}(z_0)} G_{D_{4r}(z_0)} (x, y) \mu_n (dy)\\
  && +\int_{D_{4r}(z_0)} G_{D_{4r}(z_0)}(x,y)\big[\int_{B(z_0, 4r)^c}h(z) N(y, dz)\big]\,m(dy)\\
  &=& :h_1(x)+h_2(x).
  \end{eqnarray*}
 By condition \eqref{e:1.4a}, there exists $c_1$ such that for any $r\in (0, \delta_1)$ and  $x_1, x_2\in D_r(z_0),$
\begin{eqnarray}\label{e:2.13'}
h_1(x_1)
&=&  \lim_{n\to \infty} \int_{D_{4r}(z_0)\setminus D_{3r}(z_0)} G_{D_{4r}(z_0)} (x_1, y) \mu_n (dy)
\nonumber \\
&\leq & c_1 \dfrac{\E_{x_1} \tau_{D_{2r}(z_0)}}{\E_{x_2}\tau_{D_{2r}(z_0)}}
 \lim_{n\to \infty} \int_{D_{4r}(z_0)\setminus D_{3r}(z_0)} G_{D_{4r}(z_0)} (x_2, y) \mu_n (dy)
\nonumber \\
&=& c_1 \dfrac{\E_{x_1} \tau_{D_{2r}(z_0)}}{\E_{x_2}\tau_{D_{2r}(z_0)}} h_1 (x_2).
 \end{eqnarray}
 Note  that  for every  $x\in D_r(z_0)$,
\begin{eqnarray*}
h_2(x)
&=&\int_{D_{4r}(z_0)} G_{D_{4r}(z_0)}(x,y)\int_{B(z_0, 4r)^c}h(z) N(y, dz)\,m(dy)\\
&=&\int_{D_{4r}(z_0)\setminus D_{3r}(z_0)} G_{D_{4r}(z_0)}(x,y)\int_{B(z_0, 4r)^c}h(z) N(y, dz)\,m(dy)\\
&&\quad+\int_{D_{3r}(z_0)} G_{D_{4r}(z_0)}(x,y) \int_{B(z_0, 4r)^c}h(z) N(y, dz)\,m(dy)\\
&=:&h_{2,1}(x)+h_{2,2}(x) .
\end{eqnarray*}
By the condition \eqref{e:1.4a}, we have for any $r\in (0, \delta_1)$ and  $x_1, x_2\in D_r(z_0),$
\begin{eqnarray}\label{e:2.14'}
&& h_{2,1}(x_1)  \nonumber \\
& \leq &c_1\dfrac{\E_{x_1} \tau_{D_{2r}(z_0)}}{\E_{x_2} \tau_{D_{2r}(z_0)}}\int_{D_{4r}(z_0)\setminus D_{3r}(z_0)}G_{D_{4r}(z_0)}(x_2, y)\int_{B(z_0, 4r)^c}  h(z) N(y, dz)\,m(dy)  \nonumber \\
&  \leq &c_1\dfrac{\E_{x_1} \tau_{D_{2r}(z_0)}}{\E_{x_2} \tau_{D_{2r}(z_0)}}h_{2, 1}(x_2)\leq c_1\dfrac{\E_{x_1} \tau_{D_{2r}(z_0)}}{\E_{x_2} \tau_{D_{2r}(z_0)}}h(x_2).
\end{eqnarray}
By Assumption (A3), there exists $c_2>0$ such that for any $r\in (0, \delta_1)$ and $x_1\in D_r(z_0),$
\begin{eqnarray} \label{e:2.15a}
h_{2,2}(x_1)&=&   \int_{D_{3r}(z_0)}  G_{D_{4r}(z_0)}(x_1, y)\,m(dy) \int_{B(z_0, 4r)^c}h(z)N(y, dz) \nonumber \\
&\leq & c_2 \E_{x_1}\tau_{D_{4r}(z_0)} \int_{B(z_0, 4r)^c} h(z) N(z_0, dz).
\end{eqnarray}
On the other hand, it follows from Assumption (A3) that there exists $c_3>0$ such that for each $x_2\in D_r(z_0),$
\begin{eqnarray}\label{e:2.16a}
h_{2,2}(x_2)&=& \int_{D_{2r}(z_0)}  G_{D_{4r}(z_0)}(x_2, y)\,m(dy) \int_{B(z_0, 4r)^c} h(z) N(y, dz)\nonumber\\
&\geq & c_3\int_{D_{2r}(z_0)}  G_{D_{2r}(z_0)}(x_2, y)\,m(dy)  \int_{B(z_0, 4r)^c} h(z) N(z_0, dz)\nonumber\\
&=& c_3 \E_{x_2}\tau_{D_{2r}(z_0)} \int_{B(z_0, 4r)^c} h(z) N(z_0, dz).
\end{eqnarray}
 Note that
by the assumption (A4), there exists a constant $c_4$ such that for any $r\in (0, \delta_1),$
$\E_x \tau_{D_{4r}(z_0)}\leq c_4\E_x \tau_{D_{2r}(z_0)}$ for $x\in D_r(z_0).$
Hence, it follows from \eqref{e:2.15a} and \eqref{e:2.16a} that for any $r\in (0, \delta_1)$ and $x_1, x_2\in D_r(z_0),$
\begin{equation}\label{e:2.18}
h_{2,2}(x_1)
\leq \dfrac{c_2c_4}{c_3}  \dfrac{\E_{x_1} \tau_{D_{2r}(z_0)}}{\E_{x_2} \tau_{D_{2r}(z_0)}}h_{2,2}(x_2)\leq \dfrac{c_2c_4}{c_3}  \dfrac{\E_{x_1} \tau_{D_{2r}(z_0)}}{\E_{x_2} \tau_{D_{2r}(z_0)}}h(x_2).
\end{equation}
Note that $h=h_1+h_2= h_1+h_{2, 1}+h_{2, 2}.$ By  \eqref{e:2.13'}, \eqref{e:2.14'} and \eqref{e:2.18},  there exists a constant
$c_5>0 $ so that for any $r\in (0, \delta_1)$ ,
\begin{equation} \label{e:2.19}
h(x_1)\leq c_5h(x_2)\dfrac{\E_{x_1} \tau_{D_{2r}(z_0)}}{\E_{x_2} \tau_{D_{2r}(z_0)}}
\quad \hbox{ for }  x_1, x_2\in D_r(z_0).
\end{equation}
Hence,  for any
 $z_0\in\partial D$,  $r\in (0, \delta_1)$ and  any two nonnegative  functions $h_1$ and $h_2$ that are
regular harmonic with respect to $X$  in  $D\cap B(z_0, 4r)$ and  vanish
on $(D^c)^r\cap B(z_0, 4r),$
$$
h_1(x)h_2(y)\leq c_5^2 h_2(x)h_1(y)  \quad \hbox{for }  x,y\in D\cap B(z_0, r).
$$
  Thus,  the desired conclusion is obtained with $r/2$ in place of  $r$.

\smallskip

(ii)  Suppose that Assumptions (A1) and (A3) and  \eqref{e:1.5} hold.
 We will prove that the condition \eqref{e:1.4a} is a necessary  condition for the scale invariant BHP for
 the discontinuous Hunt process $X$ in $D.$ Let $r_0$ and $r_1$ be the constants in the conditions (A3) and \eqref{e:1.5} respectively.

  Let $r\in (0, (r_0\wedge r_1)/4).$ It is easy to see that $x\mapsto \P_x (X_{\tau_{D_{2r}(z_0)}}\notin B (z_0, 4r) )$ is a nonnegative regular harmonic function of $X$ in $D_{2r}(z_0)$ and vanishes on
$ D^c\cap B(z_0, 2r).$ In the following, we will show that for each $y\in D_{4r}(z_0)\setminus D_{3r}(z_0), x\mapsto G_{D_{4r}(z_0)}(x, y)$ is a regular harmonic function of $X$ in $D_{2r}(z_0)$ and vanishes on
$ D^c\cap B(z_0, 2r).$

Fix $x\in D_{2r}(z_0)$.  By the strong Markov property of $X$,  for each nonnegative Borel function $f$ vanishing off
$D_{4r}(z_0)\setminus \overline{D_{2r}(z_0}),$
$$G_{D_{4r}(z_0)}f(x)= \E_x \int_{  \tau_{D_{2r}(z_0)}}^{\tau_{D_{4r}(z_0)}}  f(X_s )ds=\E_{x} G_{D_{4r}(z_0)}
f( X_{ \tau_{D_{2r}(z_0)}} ).
$$
Hence,  for $m$-a.e. $y$   in  $D_{4r}(z_0)\setminus \overline{D_{2r}(z_0)},$
\begin{equation}\label{e:2.15}
G_{D_{4r}(z_0)}(x, y)=\E_{x} G_{D_{4r}(z_0)}(X_{ \tau_{D_{2r}(z_0)}} , y).
 \end{equation}
We  claim that the equality in \eqref{e:2.15} holds for every $y$ in $D_{4r}(z_0)\setminus \overline{D_{2r}(z_0)}.$
 Indeed, it
 follows  from Definition \ref{De:1} that $y\mapsto G_{D_{4r}(z_0)}(x, y)$ and $y\mapsto \E_{x} G_{D_{4r}(z_0)}( X_{ \tau_{D_{2r}(z_0)}}, y)$ are excessive functions of the part process  $\wh X^{D_{4r}(z_0)}$ in $D_{4r}(z_0).$ We first prove that the restriction of the two functions  in $D_{4r}(z_0)\setminus \overline{D_{2r}(z_0)}$ are  excessive functions of the part process  $\wh X^{D_{4r}(z_0)\setminus \overline{D_{2r}(z_0})}.$  For each Borel set $B$ in $\sX,$ we denote by $\wh P^B_t$ the semigroup of the part process  $\wh X^B$ and $\wh \tau_B$ the first exit time from $B$ by the process $\wh X.$ In fact, we have for each $t>0$ and $y\in D_{4r}(z_0)\setminus \overline{D_{2r}(z_0}),$
\begin{equation}\label{e:2.16}
 \wh P^{D_{4r}(z_0)\setminus \overline{D_{2r}(z_0)}}_t  G_{D_{4r}(z_0)}(x, \cdot)(y)
\leq  \wh P^{D_{4r}(z_0)}_t  G_{D_{4r}(z_0)}(x, \cdot)(y)\leq G_{D_{4r}(z_0)}(x, y).
\end{equation}
Note that by the right continuity of   $\wh X_t$,  for each $y\in D_{4r}(z_0)\setminus \overline{D_{2r}(z_0}),$  $\P_y(\wh \tau_{D_{4r}(z_0)\setminus \overline{D_{2r}(z_0})}>0)=1.$
Since $G_{D_{4r}(z_0)}(x, \wh X_t)$ is right continuous in $t$ almost surely (see e.g. \cite[Theorem II.2.12]{BG}),
\begin{eqnarray}\label{e:2.17n}
 && \lim_{t\rightarrow 0+}\wh P^{D_{4r}(z_0)\setminus \overline{D_{2r}(z_0})}_t  G_{D_{4r}(z_0)}(x, \cdot)(y)  \nonumber \\
& =& \lim_{t\rightarrow 0+}\wh\E_y [G_{D_{4r}(z_0)}(x, \wh X_t); \wh\tau_{D_{4r}(z_0)\setminus \overline{D_{2r}(z_0})}>t]
\nonumber \\
&=&  G_{D_{4r}(z_0)}(x, y).
\end{eqnarray}
Consequently, by \eqref{e:2.16} and \eqref{e:2.17n}, $y\mapsto G_{D_{4r}(z_0)}(x, y)$   for $y\in D_{4r}(z_0)\setminus \overline{D_{2r}(z_0})$ is an excessive function of the part process
$\wh X^{D_{4r}(z_0)\setminus \overline{D_{2r}(z_0)}}.$
By a similar argument,  $y\mapsto \E_{x} G_{D_{4r}(z_0)}( X_{ \tau_{D_{2r}(z_0)}} , y)$  for $y\in D_{4r}(z_0)\setminus \overline{D_{2r}(z_0})$ is also an excessive function of the part process of $\wh X^{D_{4r}(z_0)\setminus \overline{D_{2r}(z_0})}.$
Hence, by  Assumption (A1) and \cite[Proposition 1.3 in Chapter VI]{BG}, \eqref{e:2.15} holds
for every $y\in  D_{4r}(z_0)\setminus \overline{D_{2r}(z_0}).$
 However, we can not  conclude that $x\mapsto G_{D_{4r}(z_0)}(x, y)$ is regular harmonic in  $D_{2r}(z_0)$ as we do not know the finiteness
of $G_{D_{4r}(z_0)}(x, y)$. Nevertheless,
 for every $y\in  D_{4r}(z_0)\setminus \overline{D_{2r}(z_0})$, by the monotone convergence theorem,
\begin{equation}\label{e:2.21}
G_{D_{4r}(z_0)}(x, y)
=\lim_{n\to \infty} \E_{x} \big[ G_{D_{4r}(z_0)}( X_{ \tau_{D_{2r}(z_0)}} , y) \wedge n\big]
 \end{equation}
for every $x\in D_{2r}(z_0).$

Suppose the scale invariant BHP in (i) of Theorem \ref{T:2.1} holds for $X$ in $D$.
As for each $r\in (0, (r_0\wedge r_1)/4)$,  $y\in D_{4r}(z_0)\setminus D_{3r}(z_0)$ and $n\geq 1$,
 $x\mapsto \P_x (X_{\tau_{D_{2r}(z_0)}}\notin B (z_0, 4r) )$ and $x\mapsto
  \E_{x} \left[ G_{D_{4r}(z_0)}( X_{\tau_{D_{2r}(z_0)}} , y) \wedge n\right]$
     are nonnegative regular harmonic functions in $D_{2r}(z_0)$ and vanish  on
$D^c \cap B(z_0, 2r),$   it follows from \eqref{e:1.6} that there exist positive constants $\delta=\delta(X, D)$ and
$c=c(X, D)$  such that for any $r\in (0, \delta)$,  $x_1, x_2\in D_r(z_0),$ $y\in D_{4r}(z_0)\setminus D_{3r}(z_0)$ and $n\geq 1$,
\begin{eqnarray}\label{e:2.22}
 && \E_{x_1} \Big[ G_{D_{4r}(z_0)}( X_{\tau_{D_{2r}(z_0)}} , y) \wedge n\Big]
 \, \P_{x_2}(X_{\tau_{D_{2r}(z_0)}}\in B^c_{4r}(z_0))  \nonumber \\
& \leq&  c  \, \E_{x_2} \Big[ G_{D_{4r}(z_0)}( X_{\tau_{D_{2r}(z_0)}}, y) \wedge n  \Big]
\, \P_{x_1}(X_{\tau_{D_{2r}(z_0)}}\in B^c_{4r}(z_0)).
\end{eqnarray}
 Passing $n\to \infty$, we have by \eqref{e:2.21} that
\begin{eqnarray}\label{e:2.7a}
    && G_{D_{4r}(z_0)}(x_1, y)     \, \P_{x_2}(X_{\tau_{D_{2r}(z_0)}}\in B^c_{4r}(z_0))\nonumber\\
 &\leq  &c   \, G_{D_{4r}(z_0)}(x_2 , y)  \, \P_{x_1}(X_{\tau_{D_{2r}(z_0)}}\in B^c_{4r}(z_0)).
 \end{eqnarray}
Note that by  Assumption (A3) and  Lemma \ref{L:2.2},
\begin{equation}\label{e:2.7b}
 \P_x(X_{\tau_{D_{2r}(z_0)}}\in B^c_{4r}(z_0))\asymp \E_x\tau_{D_{2r}(z_0)}\cdot N(z_0, B^c_{4r}(z_0) )
 \end{equation}
for $x\in D_r(z_0),$ where the comparison constant is independent of  $r\in (0, (r_0\wedge r_1)/4)$.
Hence, \eqref{e:1.4a} follows from \eqref{e:2.7a} and \eqref{e:2.7b}.
 This completes the proof.
\qed

\smallskip

In the remainder of this section, we show that   Assumption {\bf (A4)} holds for any  open subset $U\subset \R^d$ and for any  L\'evy process  in $\R^d$
 with $d\geq 1$.   Let $Y$ be a  L\'evy process in $\R^d$. Denote by $\psi$ the L\'evy exponent of $Y$; that is,
$
\E \exp(i(u, Y_t))=\exp (t\psi(u))$ for $u\in \R^d$ and $t>0.$
It is well known that there exists a symmetric non-negative definite constant matrix
$A=(a_{ij})_{1\leq i, j\leq d}$, a constant vector  $b\in \R^d$, and a non-negative Borel measure $\nu$ on $\R^d\setminus \{0\}$
satisfying $\int_{\R^d} (|z|^2 \wedge 1) \nu(dz)$ so that for $u\in \R^d$,
\begin{equation}\label{e:2.25}
\psi(u )=-\dfrac{1}{2} (u, Au)+i(b, u)+\int_{\R^d} \left(e^{i(u, z)}-1-i(u, z){\mathbbm 1}_{\{|z|\leq 1\}}\right)\nu(dz).
\end{equation}
The measure $\nu$ is called the L\'evy measure of $Y$ and $(A, b, \nu)$ is called the L\'evy triple of $Y$.
 One can deduce that the infinitesimal generator of $Y$ is
 \begin{eqnarray*}
\L f(x)&=& \dfrac{1}{2}\sum_{i, j=1}^d a_{ij} \dfrac{\partial^2}{\partial x_i x_j} f(x)\\
&& +b\cdot \nabla f(x)+\int_{\R^d} \left(f(x+z)-f(x)-\nabla f(x) \cdot z{\mathbbm 1}_{\{|z|\leq 1\}}\right)\nu(dz)
\end{eqnarray*}
 for $f\in C_c^2(\R^d).$ For $r>0$, we define
$\tilde b_r:=b  - {\mathbbm 1}_{\{r \leq 1 \}} (r) \int_{r<|z|\leq 1} z \nu(dz)
   + {\mathbbm 1}_{\{r >1\}}(r)  \int_{1<|z|\leq r} z \nu(dz)$
and
 \begin{equation}\label{e:Phi}
  \Phi(r):= \frac1{r^2} \sum_{i=1}^d a_{ii}+\int_{\R^d} \big(\frac{|z|^2}{r^2}\wedge 1\big) \nu(dz)
  +\frac{|\tilde b_r|}{r}.
  \end{equation}

\begin{lem}\label{L:2.10}
Let $U$ be an open set.
There exists $C=C(d)>0$  such that for any $r>0$ and $z_0\in \partial U,$
$$\P_x(Y_{\tau_{U_r(z_0)}}\in   U_{2r}(z_0)\setminus U_r(z_0) )\leq C \,  \Phi(r)  \E_x\tau_{U_r(z_0)}   \quad \hbox{for } x\in U_{r/2}(z_0),$$
where $U_{r}(z_0):=U\cap B(z_0, r).$
\end{lem}

\pf Let $\varphi\in C_c^\infty(\R^d)$ be such that $0\leq \varphi\leq 1,$  $\varphi(y)=0$
if $|y|\leq 1/2$ and $\varphi(y)=1$ if $y\in B(0, 2)\setminus B(0, 1).$ Let $z_0\in \partial U.$
For $r>0,$ define $\varphi_{z_0,r}(y)=\varphi((y-z_0)/r).$
We have
\begin{eqnarray*}
&&\int_{\R^d} \left(\varphi_{z_0,r}(y+z)-\varphi_{z_0,r}(y)-\nabla \varphi_{z_0,r} (y) \cdot z{\mathbbm 1}_{\{|z|\leq 1\}}\right) \nu(dz)\\
&=&\int_{|z|\leq  r } \left(\varphi_{z_0,r}(y+z)-\varphi_{z_0,r}(y)-\nabla \varphi_{z_0,r} (y) \cdot z \right) \nu(dz)\\
&-&{\mathbbm 1}_{(0, 1]}(r)\cdot \int_{r<|z|\leq 1} \nabla \varphi_{z_0,r} (y) \cdot z \nu(dz)\\
&+&{\mathbbm 1}_{(1, +\infty)}(r) \int_{1<|z|\leq r} \nabla \varphi_{z_0,r} (y) \cdot z \nu(dz)
+\int_{|z|> r } (\varphi_{z_0,r}(y+z)-\varphi_{z_0,r}(y)) \nu(dz).
\end{eqnarray*}
Note that
$\dfrac{\partial^2}{\partial y_iy_j} \varphi_{z_0,r}(y)=r^{-2}\dfrac{\partial^2 \varphi }{\partial y_iy_j}  ((y-z_0)/r) $  and $\nabla   \varphi_{z_0,r} (y) = r^{-1} \nabla \varphi ((y-z_0)/r)$.
Hence,
\begin{eqnarray*}
&& |\L \varphi_{z_0,r}(y)|\\
&\leq &\big|\sum_{i, j=1}^d a_{ij} \partial^2_{ij} \varphi_{z_0,r}(y)\big|+\int_{|z|\leq  r } \big|\varphi_{z_0,r}(y+z)-\varphi_{z_0,r}(y)-\nabla \varphi_{z_0,r} (y) \cdot z\big| \nu(dz)\\
&& +|\tilde b_r\cdot \nabla \varphi_{z_0,r}(y)|
+\int_{|z|> r } |\varphi_{z_0,r}(y+z)-\varphi_{z_0,r}(y)| \nu(dz)\\
&\leq &r^{-2} \| D^2 \varphi\|_\infty \Big( \sum_{i, j=1}^d |a_{ij}|+ \int_{|z|\leq  r } |z|^2\nu(dz)\Big)
 +\| \nabla \varphi\|_\infty r^{-1}|\tilde b_r|  \\
&& + \| \varphi \|_\infty \nu(\{z: |z|> r \}) \\
&\leq & \big( d \| D^2 \varphi\|_\infty + \| \nabla \varphi \|_\infty + \| \varphi\|_\infty  \big) \Phi ( r ) =: c_0 \, \Phi ( r ).
\end{eqnarray*}
 Here we used the fact that for symmetric non-negative definitive matrix $(a_{ij})$, $|a_{ij}| \leq \sqrt{a_{ii}a_{jj} } \leq (a_{ii}+a_{jj})/2$
and so  $\sum_{i,j=1}^d |a_{ij}| \leq d \sum_{i=1}^d a_{ii}$.
Thus  by Ito's formula, for $x\in U_{r/2}(z_0),$
\begin{eqnarray*}
&& \P_x(Y_{\tau_{U_r(z_0)}}\in U_{2r}(z_0)\setminus U_r(z_0)) \\
 &\le&  \E_x \varphi_{z_0 ,r}(Y_{\tau_{U_r(z_0)}})
 \, \leq \,  \E_x\int_0^{\tau_{U_r(z_0)}} |\L \varphi_{z_0,r}(Y_s)|\,ds
 \, \leq \, c_0 \Phi( r ) \E_x\tau_{U_r(z_0)}.
 \end{eqnarray*}
\qed

\begin{lem}  \label{L:1.1}
For every $r>0$,
\begin{equation}\label{e:5.5}
   \Phi (r) /16 \leq \Phi (2r) \leq 3 \Phi (r).
   \end{equation}
\end{lem}

\pf  Note that for  $r>0$,
$$
\Phi(2r)=\frac1{4r^2} \sum_{i=1}^d a_{ii}+\int_{\R^d} \Big(\frac{|z|^2}{4r^2}\wedge 1\Big) \nu(dz) +\frac1{2r} |\tilde b_{2r}|.
$$
Thus  we have
\begin{eqnarray*}
\Phi(2r)
&\leq &
\frac1{4r^2} \sum_{i=1}^d a_{ii}+\int_{\R^d} \Big(\frac{|z|^2}{4r^2}\wedge 1\Big) \nu(dz)
+\frac1{2r} |\tilde b_r|
 +   \frac1{2r}  \int_{ r< |z|\leq 2r} |z| \nu(dz)  \\
 &\leq & \frac1{4r^2} \sum_{i=1}^d a_{ii}+\int_{\R^d} \Big(\frac{|z|^2}{4r^2}\wedge 1\Big) \nu(dz)
 +\frac1{2r} |\tilde b_r|
  +   \int_{ r< |z|\leq 2r}  \frac{|z |^2}{2r^2} \nu(dz) \\
 &\leq & \frac1{4r^2} \sum_{i=1}^d a_{ii}+ 3\int_{\R^d} \Big(\frac{|z|^2}{4r^2}\wedge 1\Big) \nu(dz)
  +\frac1{2r} |\tilde b_r|
  \leq  3 \Phi (r).
  \end{eqnarray*}
Similarly, we have
\begin{eqnarray*}
 \Phi(2r)
&\geq & \frac1{4r^2} \sum_{i=1}^d a_{ii}+\int_{\R^d} \Big(\frac{|z|^2}{4r^2}\wedge 1\Big) \nu(dz) +\frac1{16r} |\tilde b_{2r}|\\
 &\geq &   \frac1{4r^2} \sum_{i=1}^d a_{ii}+\int_{\R^d} \Big(\frac{|z|^2}{4r^2}\wedge 1\Big) \nu(dz)
 +\frac1{16r} |\tilde b_r|-   \frac1{16r}  \int_{ r< |z|\leq 2r} |z| \nu(dz)  \\
  &\geq & \frac1{4r^2} \sum_{i=1}^d a_{ii}+  \int_{\R^d} \Big(\frac{|z|^2}{4r^2}\wedge 1\Big) \nu(dz)
  +\frac1{16r}|\tilde b_r|
  -  \frac1{8} \int_{ r< |z|\leq 2r}  \frac{|z |^2}{r^2} \nu(dz) \\
 &\geq & \frac1{4r^2} \sum_{i=1}^d a_{ii}+\frac18 \int_{\R^d} \Big(\frac{|z|^2}{r^2}\wedge 1\Big) \nu(dz)
  +\frac1{16r} |\tilde b_r|
  \geq    \Phi (r)/16.
  \end{eqnarray*}
  This establishes the inequality \eqref{e:5.5}.
\qed

\begin{remark}\label{R:1.2} \rm  The proof of Lemma \ref{L:1.1}  can be easily modified to show that for any $c_1>1$, there is a constant
$c_2>1$ so that
\begin{equation}\label{e:1.2}
   c_2^{-1} \Phi (r)   \leq \Phi (c_1r) \leq c_2 \Phi (r) \quad \hbox{for every } r>0.
   \end{equation}
\end{remark}

\smallskip

 The following well known result is due to
  Pruitt \cite[Theorem 1 and Section 3]{P}.
  However, the proof there
  is only given for one-dimensional L\'evy processes having no Gaussian components.
  For reader's convenience, we provide a detailed proof of Theorem \ref{T:0.3}
  for any L\'evy process in $\R^d$ in the Appendix.
  See  \cite[Theorem 4.7]{Sch} for a result  for  Feller processes  in terms of their symbols.

\begin{thm}\label{T:0.3}
Let $Y$ be a L\'evy process in $\R^d.$ There exists $C=C(d)>1$ such that for every $x\in\R^d$ and $r>0,$
 \begin{equation}\label{e:0.5}
 \dfrac{C^{-1}}{\Phi(r)}\leq\E_x \tau_{B(x, r)}\leq \dfrac{C}{ \Phi(r)}.
 \end{equation}
\end{thm}

\begin{thm}\label{T:2.8}
Assumption {\bf  (A4)} holds true   for  any L\'evy process $Y$ in $\R^d$
 and for any open subset $D\subset \R^d$.
 That is,  for  any   open set  $D\subset \R^d$,  there exists $C=C(d)>1$ such that for any $z_0\in \partial D$ and $r>0,$
\begin{equation}\label{e:3.17}
\E_x\tau_{D\cap B(z_0, 4r)}\leq C \, \E_x\tau_{D\cap B(z_0, 2r)}, \quad x\in D\cap B(z_0, r).
\end{equation}
  \end{thm}

 \pf  Let $z_0\in \partial D.$ For notational  simplicity, for each $r>0,$ let $D_r(z_0):=D\cap B(z_0, r).$
By the strong Markov property of $Y,$ for  $x\in D_r(z_0)$ and $ y\in D_{4r}(z_0)$,
$$
 G_{D_{4r}(z_0)}(x, y)=G_{D_{2r}(z_0)}(x, y)+\E_x  [G_{D_{4r}(z_0)}(Y_{\tau_{D_{2r}(z_0)}}, y); Y_{\tau_{D_{2r}(z_0)}}\in D_{4r}(z_0)]  .
 $$
 Integrating $y$ over $D_{4r}(z_0)$, we have by Lemmas \ref{L:2.10} and \ref{L:1.1} and Theorem  \ref{T:0.3} that  for $x\in D_r(z_0),$
\begin{eqnarray*}
\E_x\tau_{D_{4r}(z_0)}&=& \E_x \tau_{D_{2r}(z_0)}+ \E_x  [\E_{Y_{\tau_{D_{2r}(z_0)}}} \tau_{D_{4r}(z_0)} ; Y_{\tau_{D_{2r}(z_0)}}\in D_{4r}(z_0)] \\
 &=&\E_x \tau_{D_{2r}(z_0)}+ \E_x  \big[ \sup_{z\in D_{4r}(z_0)} \E_{z} \tau_{D_{4r}(z_0)} ; Y_{\tau_{D_{2r}(z_0)}}\in D_{4r}(z_0) \big] \\
&\leq & \E_x \tau_{D_{2r}(z_0)}+\P_x(Y_{\tau_{D_{2r}(z_0)}}\in D_{4r}(z_0))  \, \sup_{z\in \R^d} \E_{z} \tau_{B(z, 8r)}  \\
&\leq & \E_x \tau_{D_{2r}(z_0)}+ c_1\Phi(2r) \E_x \tau_{D_{2r}(z_0)}  \,  \dfrac{c_2}{\Phi(8r)} \\
&\leq &\E_x \tau_{D_{2r}(z_0)}+   c_1 \,  c_3 \,  \E_x \tau_{D_{2r}(z_0)} \\
 &\leq & (1+c_1c_3) \E_x \tau_{D_{2r}(z_0)},
\end{eqnarray*}
where $c_k=c_k(d), k=1, 2, 3.$
This establishes \eqref{e:3.17}.
\qed

\section{Scale invariant BHP in Lipschitz domain}

 In this section,  we shall establish  the scale invariant BHP  result for  discontinuous subordinate Brownian motions having Gaussian components in Lipschitz domains in $\R^d$ with $d\geq 2$  as stated  in Theorem \ref{T:1.10}.  As a consequence, we  deduce in Corollary \ref{C1}   that the scale invariant BHP holds in  any   $C^1$ domains with uniform modulo of continuity  for these subordinate Brownian motions.

Throughout this section, unless specified,  $X$ is a  discontinuous
 subordinate Brownian motion having Gaussian component in $\R^d$  with $d\geq 2$. More specifically,
let $W$ be a Brownian motion in $\R^d$ with Laplacian $\Delta$ as its generator,
and $S$ a  discontinuous subordinator that has a positive drift $b>0$ and is independent of $W$.
Define $X_t:=W_{S_t}$ and  the Laplace exponent of $S$ is
\begin{equation}\label{e:3.1}
\phi(\lambda)=b\lambda +\int_0^\infty (1-e^{-\lambda t})\,\mu(dt),
\end{equation}
where $b> 0$ and $\mu$ is a non-trival measure
 on $[0, \infty)$ satisfying $\int_0^\infty (1\wedge t)\mu(dt)<\infty$.
 In the following, for notational simplicity,     we  assume   without loss of generality that  the drift $b=1.$

 The subordinate Brownian motion $X$ is a symmetric L\'evy process with L\'evy exponent $ \phi  (\lambda^2)$.
  That is,
 $$\E_x [e^{i\xi\cdot (X_t-X_0)}]=e^{-t \phi(|\xi|^2)}, \quad x, \xi\in\R^d.$$
It is a Feller process having strong Feller property
and the infinitesimal generator is
$$
 \sL^{\phi}  f (x) =  \Delta f (x) + \int_{\R^d \setminus \{0\}} (f(x+z) -f(x) - \nabla f(x) \cdot z {\mathbbm 1}_{\{|z| \leq 1\}})
\nu (dz),
$$
for $f\in C^2_b (\R^d)$, where $\nu (dz)= j(|z|) dz$ is the L\'evy measure of the L\'evy process $X$ with
\begin{equation}\label{e:3.2}
j(r) = \int_{(0, \infty)} (4\pi t)^{-d/2} e^{-r^2/(4t)} \mu (dt), \quad r>0.
\end{equation}
The L\'evy system of $X$ is $(N(x, dy), t)$ with
\begin{equation}\label{e:3.3}
N(x, dy)= \nu (dy -x) = j(|y-x|)dy.
\end{equation}

The potential measure of the subordinator $S_t$ is defined to be
$$U(A):=\E \int_0^\infty  {\mathbbm 1}_{\{S_t\in A\}}\,dt.$$
 By a result of Reveu (see \cite[Proposition 1.7]{Ber2}),
  $U(dx)$ is absolutely continuous with respect to the Lebesgue measure on $[0, \infty)$,
  has a strictly positive bounded continuous density function $u (x)$ on $[0, \infty)$ with $ u(0+)=1$.
 In fact,
\begin{equation}\label{e:3.1a}
  u(x) =   \P (\hbox{there is some $t\geq 0$ so that } S_t=x)  \quad \hbox{for }  x\geq 0.
 \end{equation}

For each open set $D,$ let $W^D$ be the part process of Brownian motion $W$ killed upon leaving $D$ and  $Z^{D}_t:=W^{D}(S_t)$.
The process $Z^D$ is called  the subordinate killed Brownian motion in $D$.
 We will use $\zeta$ to denote the life time of the process $Z^D_t.$
 Let $U^{Z^D}(\cdot, \cdot)$ be the occupation density function of $Z^D.$ It follows from  \cite[(4.3)]{SV2} that
 \begin{equation}\label{e:3.5'}
 U^{Z^D}(x, y)=\int_0^\infty p^W_D(t, x, y) u(t)dt,
 \end{equation}
 where $p^W_D$ is the transition density function of the part process $W^D$ of Brownian motion $W$ killed upon exiting  $D.$
  Denote by $G^W_D$ the Green function of $W$ killed upon $D$.

\begin{prp}\label{P:3.1}
Let $D$ be a bounded   Lipschitz  domain with characteristics $(R_0, \Lambda_0)$ in $\R^d$. There exist positive constants $C_k=C_k(d, \phi, R_0, \Lambda_0, {\rm diam}(D))>1, k=1, 2$ such that
$$C_1^{-1}G^W_D(x, y)\leq U^{Z^D}(x, y)\leq G^W_D(x, y) \quad {\rm for}\: x,  y\in D\  \hbox{ and }   d\geq 3$$
 and for $x, y\in D$ with $|x-y|>R_0/2$ and $d=2,$
$$C_2^{-1}G^W_D(x, y)\leq U^{Z^D}(x, y)\leq G^W_D(x, y).
$$
\end{prp}

\pf Note that   by \eqref{e:3.1a},
the potential density function $u$ of $S_t$ is
bounded by $1$.
So by  \eqref{e:3.5'},
$$
 U^{Z^D}(x, y)\leq \int_0^\infty p^W_D(t, x, y) dt= G^W_D(x, y), \quad {\rm for}\: x,  y\in D.
 $$

  Set  $r(x, y):=\delta_D(x)\vee\delta_D(y)\vee|x-y|.$ Let $x_0$ be a fixed point in $D.$
  For $x, y\in D$,  define  $A_{x, y}=x_0$ if $r(x, y)\geq R_0/32$,  and when $r:=r(x, y)\leq R_0/32,$ $A_{x, y}$ is a point in $D$ such that $B(A_{x, y}, \kappa_0 r)\subset D\cap B(x, 3r)\cap B(y, 3r)$ with $\kappa_0:= {1}/{(2\sqrt{1+\Lambda_0^2})}.$
Let   $g_D(x):=G^W_D(x_0, x)\wedge 1$ for $x\in D.$
 Since the scale invariant BHP holds for $\Delta$ in Lipschitz domains,
    by a similar argument as that for  \cite[Proposition 2.2]{AAC}, we have when $d\geq 3,$
 \begin{equation}\label{e:3.1'}
 c_1^{-1}\dfrac{g_D(x) g_D(y)}{g_D^2(A_{x, y})}|x-y|^{2-d}
 \leq  G^W_D(x, y)
 \leq  c_1\dfrac{g_D(x) g_D(y)}{g_D^2(A_{x, y})}|x-y|^{2-d}
  \end{equation}
  and when $d=2,$
 \begin{equation}\label{e:3.1''}
 c_2^{-1}\dfrac{g_D(x) g_D(y)}{g_D^2(A_{x, y})}\leq G^W_D(x, y)
 \leq c_2\dfrac{g_D(x) g_D(y)}{g_D^2(A_{x, y})}\log\left(1+|x-y|^{-2}\right),
  \end{equation}
where $c_k=c_k(d,  R_0, \Lambda_0, {\rm diam}(D)), k=1, 2.$

It is known that for a bounded Lipschitz domain $D,$ there exist $A_k, k=1,2$ such that for any $x\in D$ and $r>0$, one can find a point
$x_r$ in $D$ satisfying  $|x-x_r|\leq A_1(1\wedge r)$ and $\delta_D(x_r)\geq A_2 (1\wedge r).$
By the scale invariant parabolic BHP on Lipschitz domains in \cite[Theorem 1.6]{FGS} and a similar proof in \cite{Zh} (see also Remark 1.3 in \cite{Zh}), for fixed $T>0,$ there exist positive constants $c_k=c_k(d, R_0, \Lambda_0, {\rm diam}(D), T), k=3, 4$ such that for any $x, y\in D$ and $t\in (0, T),$
\begin{equation}\label{e:3.9'}
p^W_D(t, x, y)\geq c_3 \dfrac{g_D(x)g_D(y)}{g_D(x_{\sqrt t})g_D(y_{\sqrt t})} t^{-d/2} \exp(-c_4\dfrac{|x-y|^2}{t}).
 \end{equation}
 Hence, by \eqref{e:3.9'},
\begin{eqnarray}\label{e:3.10}
&&\int_0^{( {\rm diam} (D))^2} p^W_D(t, x, y) dt \nonumber\\
&\geq & c_3 \int_0^{({\rm diam} (D))^2} \dfrac{g_D(x)g_D(y)}{g_D(x_{\sqrt t})g_D(y_{\sqrt t})} t^{-d/2} \exp(-c_4\dfrac{|x-y|^2}{t})\,dt\nonumber \\
&\geq & c_5 \int_0^{r(x,y)^2}  \dfrac{g_D(x)g_D(y)}{g^2_D(A_{x, y})} t^{-d/2} \exp(-c_4\dfrac{|x-y|^2}{t})\,dt\nonumber \\
&=& c_5|x-y|^{2-d}\dfrac{g_D(x)g_D(y)}{g^2_D(A_{x, y})} \int_{|x-y|^2/r(x, y)^2}^\infty  s^{d/2-2} \exp(-c_4s) \,ds\nonumber \\
&\geq & c_5|x-y|^{2-d}\dfrac{g_D(x)g_D(y)}{g^2_D(A_{x, y})} \int_1^\infty  s^{d/2-2} \exp(-c_4s) \,ds \nonumber \\
&=& c_6|x-y|^{2-d}\dfrac{g_D(x)g_D(y)}{g^2_D(A_{x, y})},
\end{eqnarray}
where $c_k=c_k(d, R_0, \Lambda_0, {\rm diam}(D)), k=5, 6,$ the second inequality  holds as
$g_D(x_{\sqrt t})\vee g_D(y_{\sqrt t}) \leq c_7g_D(A_{x, y})$ for $t\leq r(x, y)^2$
 by the Carleson estimate of $W$ in Lipschitz domain (see e.g.  \cite[Lemma 4.4]{JK}).
Since $u(t)$ is  continuous and strictly positive with $ u(0+)=1$, there exists
$c_8=c_8( \phi, {\rm diam}(D))>0$
such that $u(t)\geq c_8$ for $t\in (0, ({\rm diam}(D))^2).$ We have by \eqref{e:3.10},
\begin{eqnarray*}
 U^{Z^D}(x, y)&\geq & \int_0^{({\rm diam} (D))^2} p^W_D(t, x, y) u(t)dt\geq c_8\cdot
 \int_0^{({\rm diam} (D))^2} p^W_D(t, x, y) dt\\
 & \geq &
  c_8c_6|x-y|^{2-d}\dfrac{g_D(x)g_D(y)}{g^2_D(A_{x, y})}.
 \end{eqnarray*}
This together with \eqref{e:3.1'} and \eqref{e:3.1''} proves the result.
 \qed

\smallskip

Next, we  give a result on Green functions of Brownian motion in Lipschitz domains with interior cone condition with common angle $\theta\in (\cos^{-1}(1/\sqrt d), \pi)$. Recall that
 for   $\theta\in (0, \pi)$ and $r>0,$ define
$$\Gamma_\theta(r):=\{x=(x_1, \ldots, x_d) \in \R^d: |x|<r \hbox{ and } x_d>|x|\cos\theta\}.$$
A domain $D$ satisfies the  interior cone
condition with common angle $\theta,$
 if
there is $a_0>0$ such that
for every point $z\in \partial D,$
there is a cone $\Gamma_\theta(z, a_0)\subset D$ with vertex at $z$ that is conjugate to $\Gamma_\theta(a_0);$
that is, $\Gamma_\theta(z, a_0)$ is the cone with vertex at $z$ that is obtained from $\Gamma_\theta(a_0)$ through parallel translation and rotation.

Let $D$ be a Lipschitz domain with characteristics $(R_0, \Lambda_0).$
 Without loss of generality, we  may and do  assume $R_0\in (0, 1).$
It is well known that
there exists   $ \kappa=\kappa(R_0, \Lambda_0)\in (0, 1/4)$ such that for $r\in (0, R_0)$
and $z\in\partial D,$
\begin{equation}\label{e:3.6}
\hbox{there exists } z_r\in D\cap \partial B(z, r)  \hbox{ with }
\kappa r\leq \delta_D(z_r)< r.
 \end{equation}
   In the following,  we always use $\kappa $ to denote the positive constant   in \eqref{e:3.6}.

\begin{lem}\label{L:2.5}
Let $D$ be a Lipschitz domain in $\R^d$ with characteristics $(R_0, \Lambda_0)$ that satisfies the interior cone condition with common angle $\theta\in (\cos^{-1}(1/\sqrt d), \pi)$.
Then there exist  constants $C=C(d, R_0, \Lambda_0, \theta )\in (0, 1)$ and $p=p(d, \theta )\in (0, 2)$  such that  for any $z_0\in\partial D$ and $r\in (0, R_0/4),$
$$
G^W_{D\cap B(z_0, 4r)}(x, y_0)\geq Cr^{2-d}(\delta_D(x)/r)^p, \quad x\in D\cap B(z_0, r/(2(1+\Lambda_0))),
$$
 where $y_0$ is a point in $D\cap   \partial B(z_0, 3r)$ with $3\kappa r\leq \delta_D(y_0)< 3r$.
 \end{lem}

\pf Under the hypothesis of the lemma,  by decreasing the value of $R_0$ if needed,  for every $z\in \partial D$, there is a cone $\Gamma_\theta (z, R_0)\subset D$ with vertex at $z$ that is conjugate  to $\Gamma_\theta (R_0)$.
 Let $z_0\in \partial D.$ For the simplicity of notation, for each $r>0,$ let $D_r(z_0):=D\cap B(z_0, r).$ By the scaling property of Brownian motion, for  $r>0,$
\begin{equation}\label{e:3.13'}
G^W_{D_{4r}(z_0)}(x, y)=r^{2-d}G^W_{r^{-1}D_{4r}(z_0)}(r^{-1}x, r^{-1}y),
\quad \hbox{for }x, y\in D_{4r}(z_0) .
\end{equation}

 Fix $r\in (0,  R_0/4).$   Let $x\in D_{r/(2(1+\Lambda_0))}(z_0)$.
  Set  $\tilde z_0:=r^{-1}z_0$ and   $\tilde x:=r^{-1}x.$
 Then $\tilde z_0 \in \partial (r^{-1}D)$ and $\tilde x\in (r^{-1}D)\cap  B(\tilde z_0, 1/ (2(1+\Lambda_0))$.
 Note that  $r^{-1}D$ is  a Lipschitz domain with characteristics $( r^{-1}R_0, \Lambda_0)$
 such that
   for every $z\in \partial ( r^{-1} D)$, there is a cone $\Gamma_\theta (z, r^{-1}  R_0)\subset r^{-1} D$ with vertex at $z$ that is conjugate  to $\Gamma_\theta (r^{-1}R_0)$.
 Denote by  an orthonormal coordinate system
$CS_{\tilde z_0}: y=(y_1, \cdots, y_{d-1}, y_d)=:(y^{(d-1)}, y_d)\in \R^{d-1}\times \R$ with its origin at $\tilde z_0$
and    a Lipschitz function
$\varphi_{\tilde z_0}: \R^{d-1}\rightarrow \R$ with Lipschitz constant $\Lambda_0$
 such that
\begin{eqnarray}\label{e:0.3n}
\hskip -0.3truein &&B(\tilde z_0, r^{-1}R_0)\cap r^{-1}D \nonumber\\
\hskip -0.3truein &= &    \big\{y=(y^{(d-1)}, y_d)\in B(0,r^{-1} R_0) \hbox{ in } CS_{\tilde z_0}: y_d>\varphi_{\tilde z_0}(y^{(d-1)}) \big\}.
\end{eqnarray}

 Let $(\tilde x^{(d-1)}, \tilde x_d)\in \R^{d-1}\times \R$ be the coordinate of $\tilde x$ in the  coordinate system $CS_{\tilde z_0}$ and
 let  $\tilde z_x=(\tilde x^{(d-1)}, \varphi_{\tilde z_0}(\tilde x^{(d-1)}))
 \in \partial (r^{-1}D)$.
Then
$$
|\tilde x-\tilde z_x|=|\tilde x_d-\varphi_{\tilde z_0}(\tilde x^{(d-1)})|\leq |\tilde x_d|+|\varphi_{\tilde z_0}(\tilde x^{(d-1)})|\leq \frac{1}{2(1+\Lambda_0)}+\frac{\Lambda_0}{2(1+\Lambda_0)}\leq \frac{1}{2}
$$
 and
 $|\tilde z_x-\tilde z_0|\leq |\tilde z_x-\tilde x|+|\tilde x-\tilde z_0| <  1.$
Let $\lambda_0:=\tan^{-1}(1/\Lambda_0) \wedge\theta.$
Observe that $r^{-1}R_0\geq 4$
 and $\tilde z_x + \Gamma_{\lambda_0}(1)  \subset r^{-1}D_{2r}(z_0)$.
 On the other hand, since $r^{-1}D$ satisfies the interior cone condition with common angle $\theta\in (\cos^{-1}(1/\sqrt d), \pi)$,
 there exists a cone $\Gamma_{\theta}(\tilde z_x, 1)\subset r^{-1}D_{2r}(z_0)$ with vertex at $\tilde z_x$ that is conjugate  to $\Gamma_\theta (1)$.
Note that
$
 (\tilde z_x + \Gamma_{\lambda_0}(1)) \cup \Gamma_{\theta}(\tilde z_x, 1)\subset r^{-1}D_{2r}(z_0)
$
 and by \eqref{e:0.3n},
$$r^{-1}D_{2r}(z_0)= \big\{y=(y^{(d-1)}, y_d)\in B(0, 2) \hbox{ in } CS_{\tilde z_0}: y_d>\varphi_{\tilde z_0}(y^{(d-1)}) \big\}.
$$
Thus there is a cone   $\Gamma^{(1)}_{\theta}(\tilde z_x, 1)$   with vertex at $\tilde z_x$ that is conjugate  to $\Gamma_{\theta} (1)$ so that
$
 (\tilde z_x + \Gamma_{\lambda_0}(1) ) \subset  \Gamma^{(1)}_{\theta}(\tilde z_x, 1)\subset r^{-1}D_{2r}(z_0).
$
 Take
 \begin{equation}\label{e:C4}
 \theta_0:=  \max \big\{ (\cos^{-1}(1/\sqrt d)+\theta)/2, \, \theta - (\lambda_0/2)\big\} .
 \end{equation}
 Then  $\theta_0\in (\cos^{-1}(1/\sqrt d), \theta)$.
Let $\Gamma_{\theta_0}(\tilde z_x, 1)$ be the cone with angle $\theta_0$ and vertex at $\tilde z_x$ that has the same axis as
$\Gamma^{(1)}_{\theta}(\tilde z_x, 1)$. Note that by \eqref{e:C4},
\begin{equation}\label{e:C5}
\tilde x \in \Gamma_{\lambda_0/2 }(\tilde z_x, 1)  \subset
\Gamma_{\theta_0}(\tilde z_x, 1) \subset \Gamma^{(1)}_{\theta}(\tilde z_x, 1) \subset  r^{-1}D_{2r}(z_0).
\end{equation}

For the simplicity of notation,  denote
  $\Gamma_{\theta_0}(\tilde z_x, 1)$
 by $U_x.$  Let $f_{U_x}$ be defined in \cite[Example 4.6.7]{Dav}, subject to the normalizations $\|f_{U_x}\|_\infty=1$ and $\lambda=0,$ which satisfies  $\Delta f_{U_x}=0$  on $U_x$ and $f_{U_x}=0$ on
 $\partial U_x\cap  \partial   \Gamma_{\theta_0} (\tilde z_x,  4).$
  By the argument in the last five lines containing (4.6.6) in \cite{Dav} and (3.1)-(3.3) in  \cite[Proposition 3.1]{AAC},
   there exist $c_1=c_1(d, \theta)>0$ and $p=p(d, \theta)\in (0, 2)$ such that
 $
 f_{U_x}(\tilde x)\geq c_1|\tilde x-\tilde z_x|^p \geq c_1(\delta_D(x)/r)^p.
 $

  Let  $\tilde y_0=r^{-1}y_0.$ Then $\tilde y_0 \in r^{-1}D\cap   \partial B(\tilde z_0, 3)$ with $3\kappa \leq \delta_D(\tilde y_0)< 3$.
   Note that $\tilde y\mapsto G^W_{r^{-1}D_{4r}(z_0)}(\tilde y, \tilde y_0)$ is harmonic in $r^{-1}D_{4r}(z_0)\setminus \{\tilde y_0\}$ and
    $
 \delta_{r^{-1}D_{4r}(z_0)}(\tilde y)  \geq \sin(\theta-\theta_0)$
  for $\tilde y\in \partial U_x\cap \partial B(\tilde z_x, 1)\subset r^{-1}D_{2r}(z_0).$
    Let $\tilde y_1\in r^{-1}D_{(3-\kappa)r}(z_0)\cap \partial B(\tilde y_0, 2\kappa).$
   By a standard Harnack chain argument,    there exist positive constants  $c_2=c_2(d, \theta,  R_0, \Lambda_0)$
   and $c_3=c_3(d, \theta,  R_0, \Lambda_0)$
    so that
 for  every $\tilde y\in \partial U_x\cap \partial B(\tilde z_x, 1),$
   $$
   G^W_{r^{-1}D_{4r}(z_0)}(\tilde y, \tilde y_0)\geq c_2 G^W_{r^{-1}D_{4r}(z_0)}(\tilde y_1, \tilde y_0)
   \geq  c_2G^W_{B(\tilde y_0, 3\kappa)}(\tilde y_1, \tilde y_0) \geq c_3,
   $$
   where the last inequality is due to the explicit Green function formula for Brownian motion in balls (see, e.g.,
    \cite[(27)-(28) in Section 2.3]{ChungZ}).
   Then by the maximum principle, we deduce that $G^W_{r^{-1}D_{4r}(z_0)}(\tilde y, \tilde y_0)\geq c_3f_{U_x}(\tilde y)$ for $\tilde y\in U_x.$
  Hence,
$$
G^W_{r^{-1}D_{4r}(z_0)}(r^{-1}x, r^{-1}y_0)\geq  c_3  f_{U_x}(r^{-1}x )\geq  c_1c_3(\delta_D(x)/r)^p
$$
for $x\in D_{r/(2(1+\Lambda_0))}(z_0).$
 This together with \eqref{e:3.13'} yields the desired estimate.
 \qed

 \begin{lem}\label{L:2.3'}
Let  $D$ be a Lipschitz domain with characteristics $(R_0, \Lambda_0)$ that satisfies the interior cone condition with common angle $\theta\in (\cos^{-1}(1/\sqrt d), \pi)$.
Then there exist $C=C(d, \phi, R_0,  \Lambda_0, \theta)>0$ and $p=p(d, \theta)\in (0, 2)$ such that for any $r\in (0, R_0/4)$ and any boundary point $z_0\in\partial D,$
$$
G_{D_{4r}(z_0)}(x, y_0)\geq Cr^{2-d}  (\delta_D(x) / r )^p $$
for $ x\in D\cap B(z_0, r/(2(1+\Lambda_0))),$
 where $D_{r}(z_0):=D\cap B(z_0, r)$ and $y_0$ is a point in $D\cap \partial B(z_0, 3r)$ with $3\kappa r\leq \delta_D(y_0)< 3r$.
\end{lem}

\pf Let $z_0\in \partial D$ and $r\in (0, R_0/4).$ Let $y_0$ be a point in $D\cap \partial B(z_0, 3r)$ with $3\kappa r\leq \delta_D(y_0)< 3r$.
It follows from  \cite{SV2} that the subordinate killed Brownian motion  $Z^{D_{4r}(z_0)}$ is a subprocess of $X^{D_{4r}(z_0)}$.
Hence,
$$
G_{D_{4r}(z_0)}(x, y_0)\geq U^{Z^{D_{4r}(z_0)}}(x, y_0), \quad x\in D_{r/(2(1+\Lambda_0))}(z_0).
$$

   Recall that  the potential density $u(\cdot)$ of the subordinator $S$  is
 continuous and strictly positive    that is bounded by $1$ with $u(0+)=1$.
 Note that $$p^W_{D_{4r}(z_0)}(t, x, y_0)=r^{-d} p^W_{r^{-1}D_{4r}(z_0)}(r^{-2}t, r^{-1}x, r^{-1}y_0), \quad x\in
  D_{r/(2(1+\Lambda_0))}(z_0).$$
We have
$$\begin{aligned}
U^{Z^{D_{4r}(z_0)}}(x, y_0)&=\int_0^\infty p^W_{D_{4r}(z_0)}(t, x, y_0) u(t)dt\\
&=\int_0^\infty r^{-d} p^W_{r^{-1}D_{4r}(z_0)}(r^{-2}t, r^{-1}x, r^{-1}y_0)u(t)\,dt\\
&\geq \int_0^{64r^2} r^{-d} p^W_{r^{-1}D_{4r}(z_0)}(r^{-2}t, r^{-1}x, r^{-1}y_0)u(t)\,dt\\
&= \int_0^{64} r^{2-d} p^W_{r^{-1}D_{4r}(z_0)}(s, r^{-1}x, r^{-1}y_0)u(r^2s)\,ds.
\end{aligned}$$
Note that $r^{-1}D_{4r}(z_0)$ is a Lipschitz domain  with   diameter less than $8$ and  its Lipschitz
 characteristics are determined by $(R_0, \Lambda_0).$
 Observe that $|r^{-1}x-r^{-1}y_0|>1$ for $x\in D_{r/(2(1+\Lambda_0))}(z_0).$
So it follows from \eqref{e:3.10}, \eqref{e:3.1'} and \eqref{e:3.1''}
 that there exists $c_1=c_1(d, R_0, \Lambda_0)>0$ such that for $x\in D_{r/(2(1+\Lambda_0))}(z_0),$
$$\int_0^{64} p^W_{r^{-1}D_{4r}(z_0)}(s, r^{-1}x, r^{-1}y_0)\,ds\geq c_1G^W_{r^{-1}D_{4r}(z_0)}(r^{-1}x, r^{-1}y_0).
$$
Since $u$ is
 continuous and strictly positive with $ u(0+)=1$,  there exists a positive constant $c_2=c_2(\phi)$ such that $u(s)\geq c_2$ for $s\in (0, 64).$
Hence we have for $x\in  D_{r/(2(1+\Lambda_0))}(z_0),$
\begin{eqnarray}\label{e:3.9}
 U^{Z^{D_{4r}(z_0)}}(x, y_0)
&\geq & \inf_{s\in (0, 64)}u(s) \cdot \int_0^{64} r^{2-d} p^W_{r^{-1}D_{4r}(z_0)}(s, r^{-1}x, r^{-1}y_0)\,ds \nonumber\\
&\geq & c_1c_2r^{2-d}G^W_{r^{-1}D_{4r}(z_0)}(r^{-1}x, r^{-1}y_0)\nonumber\\
&=& c_1 c_2 G^W_{D_{4r}(z_0)}(x, y_0),
\end{eqnarray}
 where in the last equality, we used the scaling property  of the Green function of $W.$
By \eqref{e:3.9} and Lemma \ref{L:2.5}, the desired conclusion is obtained.

 \qed

 The following lemma follows directly from   Theorem \ref{T:0.3}.

\begin{lem}\label{L:2.2'}
For each $r_0>0,$ there exists a constant $C=C(d, \phi,  r_0)>1$ such that
$
C^{-1}r^2\leq \E_x \left[ \tau_{B(x, r)} \right] \leq Cr^2
$
for any $x\in\R^d$ and $r\in (0, r_0).$
\end{lem}

  In the rest of this paper, we assume that there is a  density  $\mu (t)$ of the L\'evy measure $\mu (dt)$ of $S$
(i.e. $\mu(dt)=\mu(t)\,dt$) with respect to the Lebesgue measure on $[0, \infty)$
and  there is a constant $c>0$ so that
\begin{equation}\label{e:mu1}
\mu (t) \leq  c \mu (2t)  \   \hbox{ for } t\in (0, 8)  \quad \hbox{and} \quad
\mu(t)\leq c\mu(t+1) \     \hbox{ for } t\in (1, \infty).
\end{equation}
By \cite[Lemma 4.2]{RSV}, if the condition \eqref{e:mu1} holds, then there exists a positive constant $c_1$ such that
\begin{eqnarray}\label{e:3.19}
 &&j(r)\leq  c_1 j(2r) \  \hbox{ for }   r\in (0, 2) \nonumber\\
 &&   \hbox{and} \quad
 j(r) \leq c_1j(r+1)    \  \hbox{ for  }   r\in (1, \infty),
 \end{eqnarray}
where $j(r)$ is defined in \eqref{e:3.2}.
 Note that $j(r)$ is decreasing on $(0, \infty)$ by \eqref{e:3.2},
  we have
 \begin{eqnarray}\label{e:3.16a}
 &&j(2r)\leq j(r)\leq  c_1 j(2r) \  \hbox{  for }  r\in (0, 2) \nonumber\\
 && \hbox{and} \quad
 j(r+1)\leq j(r) \leq c_1j(r+1)  \  \hbox{ for  }   r\in (1, \infty).
 \end{eqnarray}
  Thus in view of \eqref{e:3.3}, condition \eqref{e:2.4} and hence Assumption
    {\bf (A3')}  holds  for $X$.
   Note that for this conclusion, we do not need the assumption
  that the subordinator $S$ has a positive drift.

 \smallskip

\begin{remark} \label{R:3.5} \rm
The Laplace exponent $\phi$ of a subordinator $S$ is said to be  a complete Bernstein function
if its L\'evy measure $\mu (dt)$ in \
 \eqref{e:1.12}
 has a completely monotone density with respect to the Lebesgue measure on $(0, \infty)$;
that is, if $\mu (dt) = \mu  (t) dt$ with $\mu \in C^\infty (0, \infty)$ and $(-1)^n \mu ^{(n)}(t) \geq 0$ on $(0, \infty)$ for every integer $n\geq 0$.
In this case, we say $S$ is a complete subordinator.
Suppose that $S$ is a complete subordinator having  a  positive drift $b>0$.
By \cite[Lemma 2.1]{KSV12},  there is a constant $c>0$ so that
\begin{equation}\label{e:3.12}
\mu (t)  \leq c \mu (t+1) \quad \hbox{for all } t>1.
\end{equation}
See \cite[Chapter 15]{SSV} for a list of examples of complete Bernstein functions.
See also examples of Laplace exponent $\phi$   of subordinator that satisfy the condition \eqref{e:mu1}  in Remark 1.7.
 \end{remark}

 \smallskip

The following result is obtained by \cite[Theorem 4.5]{RSV} and a similar argument in \cite[Proposition 2.2]{KSV2}.

\begin{prp}   \label{P:2.1}
There exists a positive constant $C=C(d, \phi)$ such that for any $r\in (0, 1], x_0\in\R^d$ and any function $h$ which is nonnegative and harmonic in $B(x_0, r)$ with respect to $X,$ we have
$$h(x)\leq Ch(y) \quad \hbox{for any } x,y\in B(x_0, r/2).$$
\end{prp}

 \smallskip

  For $d\geq 3,$ denote by $G(x, y)$  the Green function of $X$ in $\R^d.$
By  \cite[Theorem 3.1]{RSV},   for each $M>0,$  there exists a positive constant $c=c(d, \phi, M)$ such that
$ G(x, y)\leq c|x-y|^{2-d} $
for $ |x-y|\leq M .$
This in particular implies that for $d\geq 3,$ there exists a positive constant $c=c(d, \phi)$ such that for any $r\in (0, 1),$
\begin{equation}\label{e:3.0}
G_{B(0, r)}(x, y)\leq G(x, y) \leq c|x-y|^{2-d} \quad \hbox{for } x, y\in B(0, r).
\end{equation}

 We say an open set $D\subset \R^2$ is Greenian with respect to $X$ if the Green function $G_D(x, y)$ of $X$ in $D$ exists and is not identically infinite. For any  Greenian  (with respect to $X$)  open set $D$ in $\R^2$, and for any Borel subset $A\subset D,$ we define
\begin{eqnarray*}
 {\rm Cap}_D(A) &:=& \sup\big\{\mu(A): \mu \: \mbox{is a measure supported on} \: A    \nonumber \\
 && \hskip 0.8truein
 \mbox{with} \: \sup_{x\in D}\int_D G_D(x, y)\mu(dy)\leq 1\big\}.
\end{eqnarray*}
 The Dirichlet form $(\mathcal{E}, \mathcal{F})$ on $L^2(\R^d; dx)$ associated with $X$ is given by
$$
\mathcal{F}=W^{1, 2}(\R^d)=\Big\{u\in L^2(\R^d; dx): \dfrac{\partial u}{\partial x_i}\in L^2(\R^d; dx) \: \mbox{for every} \: 1\leq i\leq d \Big\}
$$
and for $u, v\in \mathcal{F},$
$$\mathcal{E}(u, v)=\dfrac{1}{2}\int_{\R^d} \nabla u(x)\cdot \nabla v(x)\,dx
+\int_{\R^d\times\R^d}(u(x)-u(y))(v(x)-v(y))j(|y-x|)\,dx\,dy.$$
 The following facts are known; see \cite{CF, FOT}.  Every function $u\in W^{1, 2}(\R^d)$ has an $ \mathcal{E}$-quasi-continuous version,
which is unique $ \mathcal{E}$-quasi-everywhere ($ \mathcal{E}$-q.e. in abbreviation) on $\R^d$.
We always represent $u\in W^{1, 2}(\R^d)$ by its $ \mathcal{E}$-quasi-continuous version.
For a  Greenian open set $D$ and $A\subset D$,
$$
{\rm Cap}_D(A)
=\inf\big\{\mathcal{E}(u, u): u\in W^{1, 2}(\R^d),  \  u\geq 1 \  \mathcal{E} \hbox{-q.e. on }   A
\hbox{ and }  u=0 \     \mathcal{E} \hbox{-q.e. on }   D^c   \big\}.
$$
We use ${\rm Cap}^W_D(\cdot)$ to denote the capacity measure of Brownian motion $W$ in $D.$
Let $\mathcal{E}^W$ be the Dirichlet form of Brownian motion $W.$ Since $\mathcal{E}^W\leq \mathcal{E},$  for any Greenian open set $D\subset \R^2,$
  \begin{equation}\label{e:3.16n}
  {\rm Cap}^W_D(A)\leq {\rm Cap}_D(A) \quad \mbox{for every Borel subset} \: A\subset D.
  \end{equation}

\begin{lem}\label{L:3.9}
Let $d=2.$ There exists $C=C(d)>0$ such that for any $r\in (0, 1)$ and $x\in B(0,  {3r}/{4})\setminus \{0\},$
$$G_{B(0, r)}(x, 0) \leq C \log(3r/|x|).$$
In particular,
$G_{B(0, r)}(x, 0)\leq C \log (3/a)$
for $a\in (0,  {1}/{4})$ and $x\in B(0,  {3r}/{4})\setminus B(0, ar).$
\end{lem}

\pf The proof is similar to \cite[Lemma 4.6]{CKSV1}  except that  we need to show  the constant $C$ is uniform for $r\in (0, 1).$
Fix $x\in  B(0,  {3r}/{4})\setminus \{0\}$ and let $a:=|x|/(3r).$ Then $a\in (0, 1/4).$ Since $\overline{B(0, ar)}= \overline{B(0, |x|/3)}$ is a compact subset of $B(0, r),$ there exists a capacitary measure $\mu_a$ for $\overline{B(0, ar)}$ such that ${\rm Cap}_{B(0, r)}(\overline{B(0, ar)})=\mu_a(\overline{B(0, ar)}).$
Note that $y\mapsto G_{B(0, r)}(x, y)$ is harmonic with respect to $X$ in $B(0, 2ar)=B(0, 2|x|/3).$ By the uniform Harnack inequality in Proposition \ref{P:2.1}, we have
$$\begin{aligned}
1&\geq \int_{\overline{B(0, ar)}} G_{B(0, r)}(x, y)\mu_a(dy)\\
&\geq \left(\inf_{y\in \overline{B(0, ar)}} G_{B(0, r)}(x, y)\right)\mu_a(\overline{B(0, ar)})\\
&\geq c_1G_{B(0, r)}(x, 0) {\rm Cap}_{B(0, r)}(\overline{B(0, ar)})\\
&\geq c_1G_{B(0, r)}(x, 0) {\rm Cap}^W_{B(0, r)}(\overline{B(0, ar)})\\
\end{aligned}$$
where the constant $c_1$ is independent of $r\in (0, 1).$ In the last inequality, we used \eqref{e:3.16n}.
Hence,
\begin{equation}\label{e:3.17n}
G_{B(0, r)}(x, 0)\leq  \dfrac{c_1^{-1}}{{\rm Cap}^W_{B(0, r)}(\overline{B(0, ar)})}=\dfrac{c_1^{-1}}{{\rm Cap}^W_{B(0, r)}(\overline{B(0, |x|/3)})}.
\end{equation}
By \cite[Lemma 4.5]{CKSV1}, there exists $c_2>0$ such that for any $a\in (0, 1/4),$
\begin{equation}\label{e:3.18n}
{\rm Cap}^W_{B(0, 1)}(\overline{B(0, a)})\geq \dfrac{c_2}{\log(1/a)}.
\end{equation}
By the scaling property of $W$,
 we have $ G^W_{B(0, 1)}(x, y)=r^{d-2}G^W_{B(0, r)}(rx, ry)$ for $x, y\in B(0, 1).$
 Hence, by the definition of ${\rm Cap}^W_{B(0, 1)}(\overline{B(0, a)}),$
\begin{equation}\label{e:3.19n}
{\rm Cap}^W_{B(0, 1)}(\overline{B(0, a)})
 = {\rm Cap}^W_{B(0, r)}(\overline{B(0, ar)}).
\end{equation}
The conclusion now follows from \eqref{e:3.17n}-\eqref{e:3.19n}.
\qed

 \begin{lem}\label{L:3.10}
Let $D$ be a Lipschitz domain with characteristics $(R_0, \Lambda_0).$
There exists a positive constant $\rho_0=\rho_0(d, \phi, R_0, \Lambda_0)\in (0, 1)$ such that for any $ x\in D$ with $\delta_D(x)\leq R_0/2,$
$$\P_x(X_{\tau_{B(x, 2\delta_D(x))\cap D}}\in D^c)\geq \rho_0.
$$
\end{lem}

\pf The proof is similar to \cite[Lemma 4.1]{CKSV} or \cite[Lemma 5.1]{KSV2}.
Recall that the potential density function $u$ of $S_t$ is  strictly positive and continuous   on $[0, \infty)$ with $u(0+)=1$. Hence, for each $t>0$,  there exists a positive constant $c_t>0$ such that $\inf_{s\in (0, t)}u(s)\geq c_t.$
For the simplicity of notation, let $D_x:=B(x, 2\delta_D(x))\cap D.$

Let $Z^{D_x}_t:=W^{D_x}(S_t)$.
 We will use $\zeta$ to denote the life time of the process $Z^{D_x}_t.$
 By  \cite[Corollary 4.2(ii)]{SV2}, we have  for $ x\in D$ with $\delta_D(x)\leq R_0/2,$
\begin{eqnarray}\label{e:3.26}
&&\P_x(X_{\tau_{B(x, 2\delta_D(x))\cap D}}\in D^c) \nonumber\\
&\geq & \P_x(X_{\tau_{D_x}}\in \partial D\cap B(x, 2\delta_D(x)))\nonumber\\
&\geq &  \P_x(Z^{D_x}_{\zeta-}\in \partial D\cap B(x, 2\delta_D(x)))\nonumber\\
&= &  \E_x [u(\tau^W_{D_x}); W_{\tau^W_{D_x}}\in \partial D\cap B(x, 2\delta_D(x))]\nonumber\\
&\geq &  \E_x [u(\tau^W_{D_x});  \tau^W_{D_x}\leq t \  \hbox{ and }  \  W_{\tau^W_{D_x}}\in \partial D\cap B(x, 2\delta_D(x))]\nonumber\\
&\geq &  \inf_{s\in (0, t)} u(s)  \Big( \P_x(W_{\tau^W_{D_x}}\in \partial D\cap B(x, 2\delta_D(x))-\P_x(\tau^W_{D_x}> t)  \Big) .
\end{eqnarray}
Since $D$ is a Lipschitz domain, by \cite[(4.3)]{CKSV},  there exists $c_1=c_1(d, R_0, \Lambda_0)\in (0, 1)$ such that
\begin{equation}\label{e:3.27n}
\P_x(W_{\tau^W_{D_x}}\in \partial D\cap B(x, 2\delta_D(x))\geq c_1.
\end{equation}
On the other hand, there exists $c_2=c_2(d)>0$ such that
\begin{equation}\label{e:3.28n}
\P_x \big(\tau^W_{D_x}> t \big)\leq \dfrac{\E_x\tau^W_{D_x}}{t}\leq \dfrac{\E_x\tau^W_{B(x, 2\delta_D(x))}}{t}\leq c_2\dfrac{\delta_D(x)^2}{t}\leq c_2\dfrac{R_0^2}{t}.
\end{equation}
We choose a large enough constant  $t_0=t_0(R_0, \Lambda_0)>0$ such that $c_1-c_2( R_0^2/t_0) \geq \frac{c_1}{2}.$
Hence, by \eqref{e:3.26}-\eqref{e:3.28n},
$\P_x(X_{\tau_{B(x, 2\delta_D(x))\cap D}}\in D^c)\geq {c_{t_0}c_1}/{2}.$
  This proves the lemma  by taking $\rho_0:= {c_{t_0}c_1}/{2}.$
\qed

\smallskip

\begin{lem}\label{L:2.4}
Suppose $D$ is a Lipschitz domain with characteristics $(R_0, \Lambda_0)$ and satisfies the interior cone condition with common angle $\theta\in (\cos^{-1}(1/\sqrt d), \pi)$.
There exists $C=C(d, \phi, R_0, \Lambda_0, \theta)>0$   such that for any $r\in (0, R_0/4)$ and $z_0\in \partial D,$  we have
$$\E_x \tau_{D\cap B(z_0, r/(2(1+\Lambda_0)))}\leq
Cr^d G_{D\cap B(z_0, 4r)}(x, y^\ast)$$ for  $x\in D\cap B(z_0,  r/(2(1+\Lambda_0)) ), $
where  $y^\ast\in D\cap \partial B(z_0, 3r)$ with $3\kappa r< \delta_D(y^\ast)< 3r.$
\end{lem}

\pf This proof mainly uses the idea of the``box" method  developed by Bass and Burdzy \cite{BB2} and then adapted analytically in Aikawa \cite{Ai}.

Fix $z_0\in \partial D.$ For the simplicity of notation, we let $D_r(z_0):=D\cap B(z_0, r).$
Recall that the subordinate Brownian motion has L\'evy system $(N(x, dy), t)$ with $N(x, dy)= j(|y-x|)dy$.
 By  \eqref{e:3.16a}  and a similar argument as that for Lemma \ref{L:2.2}, there exists $c_1\geq 1 $ such that for any $r\in (0, R_0/4)$ and $x\in D_{r/(2(1+\Lambda_0))}(z_0),$
\begin{eqnarray}\label{e:3.7}
c_1^{-1}\E_x \tau_{D_{r/(2(1+\Lambda_0))}(z_0)}  r^d j(r)
&\leq& \P_x(X_{\tau_{D_{r/(2(1+\Lambda_0))}(z_0)}}\in   B(y^\ast, \kappa r))\nonumber\\
&\leq & c_1\E_x \tau_{D_{r/(2(1+\Lambda_0))}(z_0)} r^d j(r).
\end{eqnarray}
Hence, it suffices to prove that there exists $C=C(d, \phi, R_0, \Lambda_0, \theta)>0$   such that for any $r\in (0, R_0/4)$ and $x\in D_{r/(2(1+\Lambda_0))}(z_0),$
\begin{equation}\label{e:2.1}
\dfrac{\P_x(X_{\tau_{D_{r/(2(1+\Lambda_0))}(z_0)}}\in   B(y^\ast, \kappa r))}{r^{d+2} j(r)}
\leq
C\dfrac{G_{D_{4r}(z_0)}(x, y^\ast)}{r^{2-d}}.
\end{equation}

  It follows from \eqref{e:3.0} and Lemma \ref{L:3.9} that for $d\geq 2$,
there exists $c_2$ such that for any $r\in (0, R_0/4),$
$
G_{D_{4r}(z_0)}(x, y^\ast)\leq  G_{B(x, 6r)}(x, y^\ast)\leq c_2 r^{2-d}$
for  $x\in  D_{r/(2(1+\Lambda_0))}(z_0).$
Hence, we can take $M_1$ such that for any $r\in (0, R_0/4),$
$M_1r^{d-2}G_{D_{4r}(z_0)}(x, y^\ast)\leq 1$  for $ x\in  D_{r/(2(1+\Lambda_0))}(z_0).$
Define for each $j\geq 0,$
$$W_j:=\{x\in D_{r/(2(1+\Lambda_0))}(z_0): 2^{-(j+1)}\leq M_1r^{d-2}G_{D_{4r}(z_0)}(x, y^\ast)< 2^{-j}\}.$$
Then
$$D_{r/(2(1+\Lambda_0))}(z_0)=\cup_{j\geq 0} W_j.$$
Define for each $j\geq 0,$
$$
U_j:=\cup_{k\geq j} W_k
 = \big\{  x\in D_{r/(2(1+\Lambda_0))}(z_0):   M_1r^{d-2}G_{D_{4r}(z_0)}(x, y^\ast)< 2^{-j} \big\}
$$
and
$$
 J_j:=\cup_{i=0}^j W_i
 = \big\{x\in D_{r/(2(1+\Lambda_0))}(z_0):   M_1r^{d-2}G_{D_{4r}(z_0)}(x, y^\ast)\geq  2^{-(j+1)} \big\} .
 $$
Let $\omega_0(x):=\P_x(X_{\tau_{D_{r/(2(1+\Lambda_0))}(z_0)}}\in  B(y^\ast, \kappa r))$ and
$$
d_j:=\sup_{x\in J_j}\dfrac{\omega_0(x)/ (r^{d+2} j(r))}{G_{D_{4r}(z_0)}(x, y^\ast)/r^{2-d}}=
\sup_{x\in J_j}\dfrac{\omega_0(x)}{ r^{2d} j(r) G_{D_{4r}(z_0)}(x, y^\ast)}.
$$
To prove \eqref{e:2.1}, it is sufficient to show that there exists $C=C(d, \phi, R_0, \Lambda_0, \theta)$ independent of $r$ such that
\begin{equation}\label{e:2.1''}
\sup_{j \geq 0}d_j\leq C<\infty.
\end{equation}

Note that by Lemma \ref{L:2.2'}, there exists a positive constant $c_3$ such that $\E_x \tau_{D_{r/(2(1+\Lambda_0))}(z_0)}\leq \E_x\tau_{B(x, 2r)}\leq c_3r^2$ for any $r\in (0, 1).$
Hence, by  \eqref{e:3.7},
\begin{equation}\label{e:2.4'}
\omega_0(x)=\P_x(X_{\tau_{D_{r/(2(1+\Lambda_0))}(z_0)}}\in   B(y^\ast, \kappa r))\leq c_1c_3  r^{d+2} j(r)
\end{equation}
for $x\in D_{r/(2(1+\Lambda_0))}(z_0).$
Note that for $x\in J_0,$ $r^{d-2}G_{D_{4r}(z_0)}(x, y^\ast)\geq \frac{1}{2M_1}.$
Thus we have  $d_0<2M_1c_1c_3.$

Now suppose $x\in W_j$ for $ j\geq 1.$
Let  $A:=B(y^\ast, \kappa r)$ and $\tau_j:=\tau_{U_j}.$
We have
\begin{equation}\label{e:2.3}
w_0(x)
 = \P_x \big( X_{\tau_j}\in J_{j-1} \hbox{ and } X_{\tau_{D_{r/(2(1+\Lambda_0))}(z_0)}}\in   A \big)
    +\P_x \left( X_{\tau_j}\in A \right).
 \end{equation}
For the first item, by the strong Markov property of $X$, it follows that for $x\in W_j,$
\begin{eqnarray}\label{e:2.7}
&& \P_x \big( X_{\tau_j}\in J_{j-1} \hbox{ and } X_{\tau_{D_{r/(2(1+\Lambda_0))}(z_0)}}\in  A \big) \nonumber \\
&=&\E_x  \big[ \P_{X_{\tau_j}}( X_{\tau_{D_{r/(2(1+\Lambda_0))}(z_0)}}\in   A );     X_{\tau_j}\in J_{j-1} \big] \nonumber\\
&\leq & d_{j-1}\E_x \big[  r^{2d} j(r) G_{D_{4r}(z_0)}(X_{\tau_j}, y^\ast); X_{\tau_j}\in J_{j-1} \big] \nonumber\\
&\leq & d_{j-1} r^{2d} j(r) \E_x   \big[  G_{D_{4r}(z_0)}  (X_{\tau_j},y^\ast)  \big] \nonumber\\
&\leq & d_{j-1} r^{2d} j(r) G_{D_{4r}(z_0)}(x,y^\ast ) .
\end{eqnarray}

In the following, we estimate the second item in \eqref{e:2.3}. Let $\delta_{U_j}:=\sup_{y\in U_j}\delta_D(y),$
 which is no larger than $r$ as $U_j \subset D_r (z_0)$.
 We assert that  there exists a positive constant $c_4$  such that for $r\in (0, R_0/4)$ and $j\geq 1,$
\begin{equation}\label{e:2.8}
\P_x(X_{\tau_j}\in A)\leq c_4(\delta_{U_j})^2 r^d j(r), \quad x\in U_j.
\end{equation}
Define  $F_j:=\{y\in D: \delta_D(y)\leq
\delta_{U_j}  \}.$
 By the definition of $\delta_{U_j},$ we have $U_j\subset F_j.$
 For simplicity,  let  $B_x:=B(x, 2\delta_{U_j})$ and  $C_x:=F_j\cap B_x.$  Define
$p_0(x, A):=\P_x(X_{\tau_{C_x}}\in A) $ and
$p_{k+1}(x, A):=\E_x[p_k(X_{\tau_{C_x}}, A); X_{\tau_{C_x}}\in F_j] $
for $k=0, 1, \cdots.$
Then $p_k(x, A)$ is the $\P_x$ probability of the event that the process $X$ goes to $A$
after exactly $k$ jumps from one set $C_y$ to another.
By   a similar argument of  \cite[(3.44)-(3.48) in Lemma 7]{BB}, we assert that for $x\in F_j,$
\begin{equation}\label{e:2.17}
\P_x(X_{\tau_{F_j}}\in A)=\sum_{k=0}^\infty p_k(x, A).
\end{equation}
 Indeed, by the strong Markov property of $X,$
$$\P_x(X_{\tau_{F_j}}\in A)=\P_x(X_{\tau_{C_x}}\in A)+\E_x [\P_{X_{\tau_{C_x}}}(X_{\tau_{F_j}}\in A); X_{\tau_{C_x}}\in F_j].$$
Define
$r_0(x, A):=\P_x(X_{\tau_{F_j}}\in A) $ and
$r_{k+1}(x, A):=\E_x[r_k(X_{\tau_{C_x}}, A); X_{\tau_{C_x}}\in F_j] $
for $ k=0, 1, \cdots.$
Then $r_k(x, A)$ is the $\P_x$ probability of the event that the process $X$ goes to $A$
after more than $k$ jumps from one set $C_y$ to another.
By the induction argument,
\begin{equation}\label{e:3.38}
\P_x(X_{\tau_{F_j}}\in A)=\sum_{i=0}^k p_i(x, A)+r_{k+1}(x, A).
\end{equation}
 By Lemma \ref{L:3.10}, there exists $\rho_0=\rho_0(d, \phi, R_0, \Lambda_0)\in (0, 1)$  such that  for any $j\geq 1,$
\begin{equation}\label{e:2.14}
 \P_x(X_{\tau_{C_x}}\in F_j)\leq \P_x(X_{\tau_{B(x, 2\delta_D(x))\cap D}}\in D)\leq 1-\rho_0, \quad {\rm for} \quad x\in F_j.
 \end{equation}
  We have by \eqref{e:2.14} and an induction argument,
$$r_{k+1}(x, A)\leq (1-\rho_0)^{k+1}\rightarrow 0  \quad {\rm as} \quad k\rightarrow\infty.$$
 This together with \eqref{e:3.38} establishes \eqref{e:2.17}.

Note that by  \eqref{e:3.16a}  and Lemma \ref{L:2.2'}, there exist positive constants $c_k, k=5, 6$ such that
$$\begin{aligned}
\sup_{x\in F_j}p_0(x, A)&= \sup_{x\in F_j}\P_x(X_{\tau_{C_x}}\in A)\\
&\leq
 \sup_{x\in F_j}\P_x(X_{\tau_{B_x}}\in A)\leq c_5\sup_{x\in F_j}\E_x \tau_{B_x}
  r^d j(r)\\
  &   \leq c_6(\delta_{U_j})^2  r^d j(r).
 \end{aligned}$$
  Let $k\geq 0.$  Suppose $\sup_{x\in F_j}p_k(x, A)\leq c_6(1-\rho_0)^k(\delta_{U_j})^2  r^d j(r),$
  then by \eqref{e:2.14},
$$\begin{aligned}
\sup_{x\in F_j}p_{k+1}(x, A)&=\sup_{x\in F_j}\E_x[p_k(X_{\tau_{C_x}}, A); X_{\tau_{C_x}}\in F_j]\\
&\leq c_6(1-\rho_0)^{k}(\delta_{U_j})^2 r^d j(r)  \sup_{x\in F_j}\P_x(X_{\tau_{C_x}}\in F_j)\\
&\leq c_6(1-\rho_0)^{k+1}(\delta_{U_j})^2  r^d j(r).
\end{aligned}$$
Hence, by the induction and \eqref{e:2.17}, we have
\begin{eqnarray*}
&&\sup_{x\in F_j}\P_x(X_{\tau_{F_j}}\in A)\leq \sum_{k=0}^\infty\sup_{x\in F_j}p_{k}(x, A)\\
&\leq& \sum_{k=0}^\infty c_6(1-\rho_0)^{k}(\delta_{U_j})^2  r^d j(r)
 \leq \dfrac{c_6}{\rho_0}(\delta_{U_j})^2  r^d j(r).
\end{eqnarray*}
Note that $U_j\subset F_j,$ we obtain \eqref{e:2.8}.

 By Lemma \ref{L:2.3'} and the definition of $U_j$,  there exist $ p=p(d, \theta)\in (0, 2)$ and a positive constant $c_7=c_7(d, \phi, R_0, \Lambda_0, \theta)$ such that for $j\geq 1$ and any $u\in U_j,$
\begin{equation}\label{e:2.5}
(\delta_D(u)/r)^p\leq c_7G_{D_{4r}(z_0)}(u, y^\ast)r^{d-2} \leq c_7M_1^{-1}2^{-j}.
\end{equation}
Hence,  it follows from \eqref{e:2.8} and \eqref{e:2.5} that for $x\in W_j,$
$$
\begin{aligned}
\dfrac{\P_x(X_{\tau_j}\in B(y^\ast, \kappa r))}{ r^{d+2}  j(r)}
\leq c_4(\delta_{U_j}/r)^2
\leq c_4 (c_7M_1^{-1}2^{-j})^{2/p}.
\end{aligned}
$$
On the other hand, for $x\in W_j,$
$G_{D_{4r}(z_0)}(x, y^\ast)r^{d-2}\geq  M_1^{-1}2^{-(j+1)}.$
Thus, for $x\in W_j,$
\begin{equation}\label{e:2.35}
\P_x(X_{\tau_j}\in B(y^\ast, \kappa r))
\leq 2c_4c_7^{2/p}M_1^{-(2-p)/p} 2^{-j(2-p)/p} r^{2d} j(r)
G_{D_{4r}(z_0)}(x,y^\ast).
\end{equation}
For each $j\geq 1,$ let  $b_j:=2c_4c_7^{2/p}M_1^{-(2-p)/p} 2^{-j(2-p)/p}.$
Hence by \eqref{e:2.3}, \eqref{e:2.7}  and \eqref{e:2.35}, for each $j\geq 1,$
$$w_0(x)\leq d_{j-1}  r^{2d} j(r) G_{D_{4r}(z_0)}(x,y^\ast)+b_j  r^{2d} j(r) G_{D_{4r}(z_0)}(x,y^\ast), \quad x\in W_j.$$
That is,    $ d_j\leq d_{j-1}+b_j$  for $j\geq 1.$
Hence
$$
\sup_{i\geq 0} d_i\leq d(0)+\sum_{j=1}^\infty b_j<2M_1c_1c_3+2c_4c_7^{2/p}M_1^{-(2-p)/p}\sum_{j=1}^\infty2^{-j(2-p)/p}<\infty.
$$
This  yields that \eqref{e:2.1''} holds. The proof is complete.
\qed

\begin{lem}\label{L:2.8}
Let $D$ be a Lipschitz domain in $\R^d$ with characteristics $(R_0, \Lambda_0)$   satisfying an   interior cone condition with common angle $\theta\in (\cos^{-1}(1/\sqrt d), \pi)$.
There exists $C=C(d, \phi,  R_0, \Lambda_0, \theta)>0$  such that for any $z_0\in \partial D, r\in (0, R_0/2)$ and $x\in D_{r/(8(1+\Lambda_0))}(z_0),$
$$
\P_x \big( X_{\tau_{D_{2r}(z_0)\setminus \overline{B(y_0, \kappa r/2)}}}\in \overline{B(y_0, \kappa r/2)} \, \big)
\geq C r^{-2} \E_x  \big[ \tau_{D_{r/2}(z_0)}\big],
$$
where $D_r(z_0):=D\cap B(z_0, r)$ and  $y_0 \in D\cap \partial B(z_0, r)$ with $\kappa r< \delta_D(y_0)< r.$
\end{lem}

\pf Let $F:=D_{2r}(z_0)\setminus \overline{B(y_0, \kappa r/2)}.$ Since $G_{D_{2r}(z_0)}(\cdot, y_0)$ is harmonic in $D_{2r}(z_0)\setminus \{y_0\},$ we have for $x\in  D_{r/(8(1+\Lambda_0))}(z_0) ,$
\begin{eqnarray}\label{e:3.16}
&& G_{D_{2r}(z_0)}(x, y_0)
= \E_x \big[ G_{D_{2r}(z_0)}(X_{\tau_F}, y_0); X_{\tau_F}\in \overline{B(y_0, \kappa r/2)} \, \big] \nonumber\\
&=& \E_x  \big[ G_{D_{2r}(z_0)}(X_{\tau_F}, y_0); X_{\tau_F}\in B(y_0, \kappa r/4) \big]\nonumber\\
&&\quad+\E_x  \big[ G_{D_{2r}(z_0)}(X_{\tau_F}, y_0);
  \kappa r/4 \leq | X_{\tau_F} -y_0| \leq \kappa r/2
 \big]\nonumber\\
&\leq & \int_{B(y_0, \kappa r/4)}\int_F G_{F}(x, y) j(|y-z|)   G_{D_{2r}(z_0)}(z, y_0)\,dy\,dz \nonumber\\
&&\hskip 0.1truein +\sup_{z\in \overline{B(y_0, \kappa r/2)}\setminus B(y_0, \kappa r/4)}  G_{B(z, 4r)}(z, y_0)\,
\P_x \big(   \kappa r/4 \leq | X_{\tau_F} -y_0| \leq \kappa r/2 \big)\nonumber\\
&\leq & c_1 j(r) \E_x\tau_F\int_{B(y_0, \kappa r/4)}G_{D_{2r}(z_0)}(y_0, z)\,dz \nonumber\\
&&\hskip 0.1truein+ c_2r^{2-d}\P_x \big(   \kappa r/4 \leq | X_{\tau_F} -y_0| \leq \kappa r/2 \big),
\end{eqnarray}
  where in the last inequality,   Lemma \ref{L:3.9}, \eqref{e:3.16a} and \eqref{e:3.0}  are used and $c_k=c_k(d, \phi, R_0, \Lambda_0)$ for $ k=1,2.$

 Note that by Lemma \ref{L:2.2'}, there exists a positive constant $c_3=c_3(d, \phi)>1$ such that
\begin{equation}\label{e:3.42}
\int_{B(y_0, \kappa r/4)}G_{D_{2r}(z_0)}(y_0, z)\,dz\leq \E_{y_0} \tau_{B(y_0, 4r)}\leq c_3r^2.
\end{equation}
Moreover, it follows from \eqref{e:3.16a} that there exists a positive constant $c_4=c_4(d, \phi, R_0, \theta)>1$ such that for any $r\in (0, 1),$
\begin{equation}\label{e:3.43}
\P_x(X_{\tau_F}\in B(y_0, \kappa r/4))
\geq c_4  r^d j(r)  \E_x \left[ \tau_F \right].
\end{equation}
Hence, by \eqref{e:3.42} and \eqref{e:3.43}, the first item on the right side of \eqref{e:3.16} satisfies
\begin{equation}\label{e:3.39'}
 c_1 j(r) \E_x\tau_F\int_{B(y_0, \kappa r/4)}G_{D_{2r}(z_0)}(y_0, z)\,dz
\leq  c_1c_3 c_4^{-1}r^{2-d}\P_x(X_{\tau_F}\in B(y_0, \kappa r/4)).
\end{equation}
Let $c_5:=(c_1c_3 c_4)^{-1}\vee c_2.$ We have by \eqref{e:3.16} and  \eqref{e:3.39'},
$$
G_{D_{2r}(z_0)}(x, y_0) \leq c_5r^{2-d}\P_x(X_{\tau_{F}}\in \overline{B(y_0, \kappa r/2)}).
$$
That is,
\begin{equation}\label{e:3.20}
\P_x \big( X_{\tau_{D_{2r}(z_0)\setminus \overline{B(y_0, \kappa r/2)}}}\in \overline{B(y_0, \kappa r/2)} \big)
\geq c_5^{-1} r^{d-2} G_{D_{2r}(z_0)} (x, y_0)
\end{equation}
for $x\in D_{r/(8(1+\Lambda_0))}(z_0).$
Hence, combining \eqref{e:3.20} with Lemma \ref{L:2.4}, there exists $c_6=c_6(d, \phi,  R_0, \Lambda_0, \theta)$ such that
$$\P_x \big( X_{\tau_{D_{2r}(z_0)\setminus \overline{B(y_0, \kappa r/2)}}}\in \overline{B(y_0, \kappa r/2)} \big)
\geq c_6 r^{-2}\E_x \tau_{D_{r/(4(1+\Lambda_0))}(z_0)}$$
for $x\in D_{r/(8(1+\Lambda_0))}(z_0).$
By Theorem \ref{T:2.8}, there exists $c_7=c_7(d, \Lambda_0)>0$ such that $\E_x \tau_{D_{r/(4(1+\Lambda_0))}(z_0)}\geq c_7\E_x \tau_{D_{r/2}(z_0)}$ for
$x\in D_{r/(8(1+\Lambda_0))}(z_0).$
Hence, the desired conclusion follows.
\qed

\smallskip

\begin{lem}\label{L:3.7}
Suppose $D$ is a Lipschitz domain with characteristics $(R_0, \Lambda_0).$
There exists a constant $C=C(d, \phi,  R_0, \Lambda_0)>0$  such that  for any $z_0\in\partial D, r\in (0, R_0/2)$ and any nonnegative harmonic function $h$
on $D\cap B(z_0, 2r)$  vanishing continuously on  $D^c\cap B(z_0, 2r),$
$$h(x)\leq Ch(x_0)  \quad \hbox{for } x\in D\cap B(z_0, r),$$
where $x_0$ is a point in $D\cap B(z_0, r)$ with $\kappa r\leq\delta_D(x_0)< r.$
\end{lem}

 \pf    Let $z_0\in \partial D.$ Define $\rho_{z_0}(x):=x_d-\varphi_{z_0}(\tilde x),$ $(\tilde x, x_d)\in \R^{d-1}\times \R$ is the coordinates of $x$ in the orthonormal coordinate system
$CS_{z_0}$ with its origin at $z_0$ and the local Lipschitz function  $\varphi_{z_0}$ with Lipschitz constant $\Lambda_0$  in
$CS_{z_0}.$
By a very similar argument in   Lemma \ref{L:3.10}  with $\rho_{z_0}(x)$ in place of $\delta_D(x)$ ,
   there exists a constant $c_1=c_1(d, \phi, R_0, \Lambda_0)>0$ such that for any $z_0\in\partial D$ and $\rho_{z_0}(x)\in (0, R_0/2),$
   \begin{equation}\label{e:3.44}
   \P_x(\tau_{D\cap B(x, 2\rho_{z_0}(x))}\in D^c)\geq c_1,
   \end{equation}
   Through  a very similar argument as that for  \cite[Theorem 5.3]{KSV2} by using  \eqref{e:3.44} and  the  condition \eqref{e:mu1} in place of \cite[Lemma 5.1]{KSV2}  and the condition  $\mu(t)\leq c\mu(2t)$ for $t\in (0, 1]$ there,   one gets the desired estimate.
      \qed

\begin{thm}\label{T4}
Suppose $D$ is a Lipschitz domain with characteristics $(R_0, \Lambda_0)$ and satisfies the interior cone condition with common angle $\theta\in (\cos^{-1}(1/\sqrt d), \pi)$.
Then there exists $C=C(d, \phi, R_0,  \Lambda_0, \theta)>0$  such that for any $z_0\in\partial D, r\in (0, R_0/2)$ and any nonnegative harmonic function $h$
on $D\cap B(z_0, 2r)$  vanishing continuously on the  points of $D^c\cap B(z_0, 2r),$
\begin{equation}\label{e:3.40'}
\dfrac{h(x)}{h(y)}\leq C \dfrac{\E_x  \left[ \tau_{D_{2r}(z_0)} \right] }{\E_y  \left[ \tau_{D_{2r}(z_0)}\right]}
\quad \hbox{for }  x,y\in D\cap B(z_0, r),
\end{equation}
where $D_r(z_0):=D\cap B(z_0, r).$
Consequently, there exists $C=C(d, \phi, R_0,  \Lambda_0, \theta)>0$  such that for any $z_0\in\partial D, r\in (0, R_0/2)$ and any two nonnegative harmonic functions $h_k, k=1,2$
on $D\cap B(z_0, 2r)$  vanishing continuously on $D^c\cap B(z_0, 2r),$
\begin{equation}\label{e:1.3b}
\dfrac{h_1(x)}{h_1(y)}\leq C \dfrac{h_2(x)}{h_2(y)}  \quad   \hbox{for } x,y\in D\cap B(z_0, r).
\end{equation}
\end{thm}

\pf Let $z_0\in \partial D$ and $r\in (0, R_0/2).$ Let $h$ be a nonnegative harmonic function
on $D\cap B(z_0, 2r)$ and vanishes continuously on the  points of $D^c\cap B(z_0, 2r).$
For the simplicity of notation, let $B_r(z_0):=B(z_0, r).$
 By the uniform Harnack inequality Proposition \ref{P:2.1} and a standard Harnack chain argument, it only suffices to prove \eqref{e:3.40'} for $x, y\in D\cap B(z_0, r/(8(1+\Lambda_0))  ).$

 Let $x_0 \in D\cap  \partial B(z_0, r)$
with $\kappa r< \delta_D(x_0)< r.$
Let $\{U_n; n\geq 1\}$ be an increasing sequence of relatively compact subsets of $D_{2r}(z_0)$ so that
$\overline U_n\subset U_{n+1}$ and $\cup_{n=1}^\infty U_n = D_{2r}(z_0)$. Let $F:=D_{2r}(z_0)\setminus \overline{B(x_0, \kappa r/2)}.$
Then for  $x\in D_{r/(8(1+\Lambda_0)) }(z_0)\subset F$,
\begin{eqnarray*}
h(x)&=& \lim_{n\to \infty} \E_x [ h(X_{\tau_{U_n\cap F}}) ]\\
&\geq &  \lim_{n\to \infty}\E_x [h(X_{\tau_{U_n\cap F}}); X_{\tau_{U_n\cap F}}\in  \overline {B(x_0, \kappa r/2)}  ]\\
&\geq & \inf_{z\in \overline{B(x_0, \kappa r/2)}}h(z)\P_x (X_{\tau_{ F}}\in \overline{B(x_0, \kappa r/2)}).\\
 \end{eqnarray*}
Hence, by the Harnack principle Proposition \ref{P:2.1} and Lemma \ref{L:2.8},
\begin{equation}\label{e:3.32a}\begin{aligned}
h(x)
&\geq c_1h(x_0)\P_x \big(X_{\tau_{D_{2r}(z_0)\setminus \overline{B(x_0, \kappa r/2)}}}\in \overline{B(x_0, \kappa r/2)} \big)\\
&\geq c_2h(x_0) r^{-2}  \E_x  \left[ \tau_{D_{r/2}(z_0)} \right],
\end{aligned}\end{equation}
where $c_k=c_k(d, \phi, R_0, \Lambda_0, \theta)>0, k=1,2.$

Since $h$ vanishes continuously on the  points of $D^c\cap B(z_0, 2r),$ by  Remark \ref{R:1.4}(ii) and a similar argument as that for \cite[Lemma 4.2]{CKSV}, $h$ is regular harmonic in $D_{r/2}(z_0).$   Thus, for  $x\in D_{r/(8(1+\Lambda_0)) }(z_0),$
\begin{equation}\label{e:3.23}\begin{aligned}
 h(x)
&=\E_x h(X_{\tau_{D_{r/2}(z_0)}})\\
&=\E_x  \left[h(X_{\tau_{D_{r/2}(z_0)}}); X_{\tau_{D_{r/2}(z_0)}}\in D_{r}(z_0) \right]\\
&\quad+\E_x \left[ h(X_{\tau_{D_{r/2}(z_0)}}); X_{\tau_{D_{r/2}(z_0)}}\in  B^c_{r}(z_0) \right].
\end{aligned}\end{equation}
 By Lemmas \ref{L:3.7} and   \ref{L:2.10},  there exists $c_3=c_3(d, \phi, R_0, \Lambda_0)$ such that
 \begin{equation}\label{e:3.33}\begin{aligned}
  &\E_x [h(X_{\tau_{D_{r/2}(z_0)}}); X_{\tau_{D_{r/2}(z_0)}}\in D_{r}(z_0)]\\
 \leq &\sup_{z\in D_{r}(z_0)}h(z)\cdot\P_x(X_{\tau_{D_{r/2}(z_0)}}\in  D_{r}(z_0))\\
\leq &c_3h(x_0) r^{-2}  \E_x  \left[ \tau_{D_{r/2}(z_0)}\right]  .
 \end{aligned}\end{equation}
For the second item in \eqref{e:3.23}, if follows from  \eqref{e:3.16a} that there exists
 $c_4=c_4(d, \phi)$
such that
\begin{eqnarray*}
&&\E_x  \big[ h(X_{\tau_{D_{r/2}(z_0)}}); X_{\tau_{D_{r/2}(z_0)}}\in B^c_{r}(z_0) \big]\\
&=& \int_{B_{r}^c(z_0)}\int_{D_{r/2}(z_0)} G_{D_{r/2}(z_0)}(x, y)  j(|y-z|) h(z)\,dy\,dz\\
&\leq & c_4   \E_x \left[ \tau_{D_{r/2}(z_0)} \right] \int_{ B^c_{r}(z_0)}j(|z-z_0|) h(z)\,dz.
\end{eqnarray*}
Let $y_0 \in D\cap \partial B(z_0, r/2)$ with $\kappa r/2< \delta_D(y_0)< r/2.$
It follows from Lemma \ref{L:2.2'} and \eqref{e:3.16a} that
\begin{eqnarray*}
h(y_0) &\geq &  \E_{y_0} \big[ h(X_{\tau_{B(y_0, \kappa r/2)}}); X_{\tau_{B(y_0, \kappa r/2)}}\in B^c_{r}(z_0) \big]\\
&= & \int_{B_{r}^c(z_0)}\int_{B(y_0, \kappa r/2)} G_{B(y_0, \kappa r/2)}(y_0, y) j(|y-z|)  h(z)\,dy\,dz\\
&\geq&  c_5\E_{y_0}\left[ \tau_{B(y_0, \kappa r/2)} \right]  \int_{ B^c_{r}(z_0)}j(|z-z_0|) h(z)\,dz\\
&\geq & c_6 r^2   \int_{ B^c_{r}(z_0)}j(|z-z_0|) h(z)\,dz,
\end{eqnarray*}
 where $c_k=c_k(d, \phi, R_0, \Lambda_0)>0, k=5,6.$ By the uniform Harnack principle Proposition \ref{P:2.1}, there exists $c_7=c_7(d, \phi)>0$ such that $h(y_0)\leq c_7 h(x_0).$
Hence,
\begin{equation}\label{e:3.25}
\E_x  \left[ h(X_{\tau_{D_{r/2}(z_0)}}); X_{\tau_{D_{r/2}(z_0)}}\in B^c_{r}(z_0) \right]
 \leq  \dfrac{c_4c_7}{c_6}h(x_0) r^{-2}  \E_x \left[ \tau_{D_{r/2}(z_0)}\right].
\end{equation}
Consequently, by \eqref{e:3.23}-\eqref{e:3.25}, for $x\in D_{r/(8(1+\Lambda_0) ) }(z_0),$
\begin{equation}\label{e:3.5}
h(x)\leq (c_3+\dfrac{c_4c_7}{c_6})h(x_0)r^{-2} \E_x  \big[\tau_{D_{r/2}(z_0)} \big].
\end{equation}
Then by \eqref{e:3.5} and \eqref{e:3.32a},
$$
c_2h(x_0) \dfrac{\E_x  \big[ \tau_{D_{r/2}(z_0)} \big] }{r^2}\leq h(x)\leq (c_3+\dfrac{c_4c_7}{c_6})h(x_0) \dfrac{\E_x  \big[ \tau_{D_{r/2}(z_0)} \big] }{r^2}
$$
for $x\in  D_{r/(8(1+\Lambda_0))}(z_0).$
Note that by Theorem \ref{T:2.8},  there exists a positive constant $c_8=c_8(d)$ such that $c_8\E_x{\tau_{D_{2r}(z_0)}}\leq \E_x{\tau_{D_{r/2}(z_0)}}\leq \E_x{\tau_{D_{2r}(z_0)}}.$ Hence, the desired conclusion is obtained.
\qed

\begin{cor}\label{C1}
 Suppose that  $D$ is a $C^1$ domain in $\R^d$ with $d\geq 2$. Then the scale invariant BHP holds for $X$ in $D.$
\end{cor}

\pf  For each $z\in\partial D,$ let $T_z$ be the tangent plane of $D$ at the boundary point $z.$
Denote by $CS_z: y=(y_1, \cdots, y_{d-1}, y_d)=:(\wt {y}, y_d)\in \R^{d-1}\times \R$ an orthonormal coordinate system with its origin at $z$ such that
$$
T_z=\big\{ y=(\wt {y}, y_d)\in \R^{d-1}\times \R \hbox{ in } CS_z: y_d=0  \big\}.
$$
Note that $D$ is a $C^1$ domain with uniform modulo of continuity.  Then there exist $R_0>0,$ a uniform continuity function $f$ in  $[0, \infty)$ with $f(0)=0$ and a $C^1$ function $\varphi_z: \R^{d-1}\rightarrow \R$ satisfying $\varphi_z(\tilde 0)=\nabla \varphi_z(\tilde 0)=0$ and $|\nabla \varphi_z(\tilde x)-\nabla \varphi_z(\tilde y)|\leq f(|\tilde x-\tilde y|)$  such that
$$
B( z, R_0)\cap D= \left\{y=(\wt {y}, y_d)\in B(0,R_0) \hbox{ in } CS_z: y_d>\varphi_z(\wt {y}) \right\}.
$$
Note that $\nabla \varphi_z(\tilde 0)=0,$ then  in the coordinate system $CS_z,$
$|\nabla \varphi_z(\tilde y)|=|\nabla \varphi_z(\tilde y)-\nabla \varphi_z(\tilde 0)|\leq f(|\tilde y|)$
when  $|\tilde y|<R_0.$
Hence, there exists a small enough constant $\ee_0$ such that in the coordinate system $CS_z,$
$|\nabla \varphi_z(\tilde y)|< 1/\sqrt{d-1}$  when  $ |\tilde y|<\ee_0.$
 Consequently,  $D$ satisfies the interior cone condition with common angle $\theta\in (\cos^{-1}(1/\sqrt d), \pi)$.
The desired conclusion is obtained by Theorem \ref{T4}.
\qed

\smallskip

   In the remaining part of this section,   we will prove the last assert in Theorem \ref{T:1.10};    that is,
       we   show that the scale invariant BHP    for $X$
     fails on a truncated cone with the angle  $\theta\in (0, \cos^{-1}(1/\sqrt d)]$ in $\R^d$ with $d\geq 2.$
     Recall that for $\theta\in (0, \pi),$
      $\Gamma_\theta$ is the vertical truncated circular cone in $\R^d$ with vertex at $0$ and  with  angle $\theta$
     defined in \eqref{e:cone} with $r=1$.
  The truncated cone $\Gamma_\theta$ is a  bounded Lipschitz domain. Denote its Lipschitz  characteristics by $(R_0, \Lambda_0)$.
 Fix some  $x_0\in \Gamma_\theta$ and
   define  $r(x, y)=\delta_{\Gamma_\theta} (x)\vee\delta_{\Gamma_\theta} (y)\vee|x-y|.$
   We let $A_{x, y}=x_0$ if $r(x, y)\geq R_0/32$ and when $r:=r(x, y)\leq R_0/32,$ $A_{x, y}$ is a point in ${\Gamma_\theta}$
   so that $B(A_{x, y}, \kappa_0 r)\subset \Gamma_\theta\cap B(x, 3r)\cap B(y, 3r)$ with $\kappa_0:= {1}/{(2\sqrt{1+\Lambda_0^2})}.$

  For each open set $B,$ denote by $\tau^W_B$ the first  time of Brownian motion $W$  exiting $B.$

\begin{lem}\label{L:4.4'}
 Let $\Gamma_\theta$ be a truncated cone with  the  angle  $\theta\in (0, \cos^{-1}(1/\sqrt d)]$ in $\R^d$ with $d\geq 2.$
 Let $x=(0, \cdots, 0, a)$ in $\R^d$ with  $a\in (0, 1/8).$ There are positive constants $C_k=C_k(d, \theta), k=1, 2$ such that
 $$
 \E_x  \left[ \tau^W_{\Gamma_\theta} \right]  \geq  C_1  |x|^2\log (1/(2|x|) \quad \mbox{when} \quad \theta=\cos^{-1}(1/\sqrt d);
 $$
and
 $$\E_x  \left[ \tau^W_{\Gamma_\theta} \right]  \geq  C_2 |x|^2 \quad \mbox{when} \quad \theta\in (0, \cos^{-1}(1/\sqrt d)).$$
 \end{lem}

 \pf     For $r\in (0, 1),$   recall that
  $
  \Gamma_\theta(r):=\{x\in \Gamma_\theta: |x|< r\}.
  $
 Let $x_0$ be a fixed point in $\Gamma_\theta\setminus \Gamma_\theta(3/4).$
 Let  $g_{\Gamma_\theta}(x):=G^W_{\Gamma_\theta}(x_0, x)\wedge 1.$
 By \eqref{e:3.1'}  and \eqref{e:3.1''}, there exists $c_1=c_1(d, \theta)$ such that for any $x=(0, \cdots, 0, a)$ in $\R^d$ with $a\in (0, 1/8),$
\begin{eqnarray*}
 \E_x  [ \tau^W_{\Gamma_\theta}]
  &=&\int_{\Gamma_\theta} G^W_{\Gamma_\theta}(x, y)\,dy
 \geq  \int_{\Gamma_{\theta/2}(1/2)\setminus \Gamma_{\theta/2}(2|x|)} G^W_{\Gamma_\theta}(x, y)\,dy\\
 &\geq & c_1\int_{\Gamma_{\theta/2}(1/2)\setminus \Gamma_{\theta/2}(2|x|)} \dfrac{g_{\Gamma_\theta}(x)}{g_{\Gamma_\theta}(A_{x, y})}|x-y|^{2-d}\dfrac{g_{\Gamma_\theta}(y)}{g_{\Gamma_\theta}(A_{x, y})}\,dy.
 \end{eqnarray*}
 Note that for $y\in \Gamma_{\theta/2}(1/2)\setminus \Gamma_{\theta/2}(2|x|),$ we can take $A_{x, y}=y.$ It follows that
\begin{equation}\label{e:3.54}
 \E_x   \left[ \tau^W_{\Gamma_\theta} \right]
   \geq c_1\int_{\Gamma_{\theta/2}(1/2)\setminus \Gamma_{\theta/2}(2|x|)} \dfrac{g_{\Gamma_\theta}(x)}{g_{\Gamma_\theta}(y)}|x-y|^{2-d}\,dy.
\end{equation}

Let $\phi_0$ be the positive eigenfunction in $L^2(\Gamma_\theta)$  corresponding to the first positive eigenvalue $\lambda_1$ of $-\Delta$ in $\Gamma_\theta$ with zero Dirichlet boundary condition, normalized to have $\int_{\Gamma_\theta} \phi_0^2(x)\,dx=1$.
It follows from Example 4.6.7  in \cite{Dav} that there exists $q=q(d, \theta)>0$ such that
$\phi_0(x)$ decays at rate $|x|^q$ as $x\rightarrow 0$ along the axis of the cone $\Gamma_\theta.$
Moreover, by virtue of (4.6.6) in \cite{Dav} and (3.1)-(3.3) in \cite{AAC},  $q=2$ if $\theta=\cos^{-1}(1/\sqrt d),$ and $q>2$  if $\theta\in (0, \cos^{-1}(1/\sqrt d)).$

It is known that $\phi_0(y)\asymp  g_{\Gamma_\theta}(y)$ for $y\in \Gamma_\theta,$  where the comparison constant depends only on $(d, \theta)$ (see  \cite[Lemma 3.2]{AAC}).
  Hence,  by the uniform Harnack inequality of $W$ and the standard chain argument,
\begin{equation}\label{e:3.28}
\begin{aligned}
& g_{\Gamma_\theta}(y) \asymp |y|^2 \quad \hbox{if }  \theta=\cos^{-1}(1/\sqrt d);\\
& \mbox{and} \quad g_{\Gamma_\theta}(y) \asymp |y|^\gamma  \quad \hbox{if }    \theta\in (0, \cos^{-1}(1/\sqrt d))
 \end{aligned} \end{equation}
  for any $y\in \Gamma_{\theta/2}(1/2),$ where $\gamma=\gamma(d,\theta)\in (2, \infty)$ and the comparison constants depend only on $(d, \theta).$
  In particular, for any $x=(0, \cdots, 0, a)$ in $\R^d$ with $a\in (0, 1/8),$
  \begin{equation}\label{e:3.56}
  \begin{aligned}
  &g_{\Gamma_\theta}(x) \asymp |x|^2 \quad \hbox{if }  \theta=\cos^{-1}(1/\sqrt d); \\
 & \mbox{and} \quad g_{\Gamma_\theta}(x) \asymp |x|^\gamma  \quad \hbox{if }    \theta\in (0, \cos^{-1}(1/\sqrt d)).
 \end{aligned} \end{equation}
Note that $|y|\leq |x-y|+|x|\leq 2|x-y|$ for   $y\in \Gamma_{\theta/2}(1/2)\setminus \Gamma_{\theta/2}(2|x|).$
 If $\theta=\cos^{-1}(1/\sqrt d),$
 then by \eqref{e:3.28}, there exists $c_2=c_2(d, \theta)$ such that $g_{\Gamma_\theta}(y)\leq c_2|y|^2\leq  4c_2|x-y|^2$ for  $y\in \Gamma_{\theta/2}(1/2)\setminus \Gamma_{\theta/2}(2|x|).$  Thus  by this together with \eqref{e:3.54} and \eqref{e:3.56}, there exist $c_k=c_k(d, \theta)$ with $k=3, 4$ such that when $\theta=\cos^{-1}(1/\sqrt d),$
 $$
 \E_x   \left[ \tau^W_{\Gamma_\theta} \right]
     \geq  c_3\int_{\Gamma_{\theta/2}(1/2)\setminus \Gamma_{\theta/2}(2|x|)} \dfrac{|x|^2}{|x-y|^2}|x-y|^{2-d}\,dy
  \geq c_4|x|^2\log(|x|^{-1}).
  $$
 If $\theta\in (0, \cos^{-1}(1/\sqrt d)),$ then by \eqref{e:3.28}, there exists $c_5=c_5(d, \theta)$ such that
 $g_{\Gamma_\theta}(y)\leq c_5|y|^\gamma\leq  c_5 2^\gamma|x-y|^\gamma$ for  $y\in \Gamma_{\theta/2}(1/2)\setminus \Gamma_{\theta/2}(2|x|).$
Hence, by this together with \eqref{e:3.54} and \eqref{e:3.56}, there exist $c_k=c_k(d, \theta)$ with $k=6, 7$ such that when $\theta\in (0, \cos^{-1}(1/\sqrt d)),$
 $$
 \E_x  \left[ \tau^W_{\Gamma_\theta} \right]
   \geq  c_6\int_{\Gamma_{\theta/2}(1/2)\setminus \Gamma_{\theta/2}(2|x|)} \dfrac{|x|^\gamma}{|x-y|^\gamma}|x-y|^{2-d}\,dy
   \geq c_7|x|^2.
 $$
 The proof is complete.
   \qed

 \begin{lem}\label{L:3.16}
Let $D$ be an open set in $\R^d.$
Suppose there exists a cone $V$ with vertex $z\in \partial D$ such that $V\cap B(z, r)\subseteq D^c$ for some $r>0.$
Then $z$ is regular with  respect to  $X$ for $D^c.$
\end{lem}

\pf Without loss of generality, we assume $z=0.$ For each $t>0,$ we have
\begin{equation}\label{e:3.47}
\begin{aligned}
 \P_0 (\tau_D\leq t) &\geq \P_0(X_t\in D^c)\geq \P_0(X_t\in V\cap B(0, r))\\
 &\geq \P_0(X_t\in V)-\P_0(X_t\notin B(0, r)).
 \end{aligned}\end{equation}
For each $\lambda>0,$ define $Y^{(\lambda)}_t:=\lambda X_{\lambda^{-2}t}.$
Then it follows from \eqref{e:3.47} and $t^{-1/2}V=V$ that
\begin{equation}\label{e:3.47'}
\P_0 (\tau_D\leq t)\geq \P_0(Y^{(t^{-1/2})}_1\in V)-\P_0(Y^{(t^{-1/2})}_1\notin B(0, r/\sqrt{t})).
\end{equation}
Let   $\phi^{(\lambda)}(r):=\lambda^{-2}\phi(\lambda^2 r)$ for $r>0.$
Note that the characteristic function of $Y^{(\lambda)}_t$ satisfies for $x, \xi\in\R^d$,
\begin{equation}\label{e:3.48}
\E_x [e^{i\xi\cdot (Y^{(\lambda)}_t-Y^{(\lambda)}_0)}]=\E_{x/\lambda} [e^{i (\lambda\xi)\cdot (X_{\lambda^{-2}t}-X_0)}]\\
=e^{-t\phi^{(\lambda)}(|\xi|^2)} .
\end{equation}
 By \eqref{e:3.1}, $\phi^{(\lambda)}(r)=r +\lambda^{-2}\int_0^\infty (1-e^{-\lambda^2 r s})\,\mu(ds)$ for $r\geq 0.$
   Note that  for positive $r$ and  $s$,
    $\lambda^{-2} (1-e^{-\lambda^2 r s})\leq \lambda^{-2}(1\wedge (\lambda^2 rs))\leq 1\wedge (rs)$  when   $\lambda>1.$ Hence, by the dominated convergence theorem,
\begin{equation}\label{e:3.51}
\lim_{\lambda \to \infty} \phi^{(\lambda)}(r) =  r   \quad \mbox{for } r\geq  0.
\end{equation}
 Then $Y^{(\lambda)}_1$ converges to $W_1$ in law as $\lambda\rightarrow\infty.$
 So by letting $t\rightarrow 0$ in \eqref{e:3.47'}, we have
$
\P_0 (\tau_D=0)\geq \P_0(W_1\in V)>0.
$
By  Blumenthal's zero-one law,  $\P_z(\tau_D=0)=1.$ That is,  $z$ is regular for $D^c.$
\qed

\begin{thm}\label{T5}
 Let $\Gamma_\theta$ be a  truncated cone in $\R^d$ with  the angle $\theta\in (0, \cos^{-1}(1/\sqrt d)]$ in $\R^d$ with $d\geq 2.$ Then the scale invariant BHP fails to hold  for $X$  in $\Gamma_\theta$ near the vertex $0.$
 \end{thm}

 \pf  We  use the contradiction method to prove this result.
 Recall the definition of $\Gamma_\theta (r)$ from \eqref{e:cone} and that $\Gamma_\theta :=\Gamma_\theta(1).$
Let $r\in (0, 1/4)$ and $y_0=( 0, \cdots, 0, 3r).$
Define  $h(y):=\P_y(  X_{\tau_{\Gamma_\theta(2r)}} \in   \Gamma^c_\theta(4r) )$
  and $g(y):=\int_{B( y_0,  (3r\sin\theta\wedge  r)/2 ))} G_{\Gamma_\theta(4r)}(y, z)\,dz.$ Then $h$ and $g$  are
   non-negative   regular harmonic functions   in $\Gamma_\theta(2r)$ with respect to $X$
   and vanish on $B(0, 2r) \setminus \Gamma_\theta(2r).$
 Note that by
   \eqref{e:3.16a} and Lemma \ref{L:2.2},
   $\P_y(  X_{\tau_{\Gamma_\theta(2r)}} \in   \Gamma_\theta(4r)^c )
   \asymp \E_y\tau_{\Gamma_\theta(2r)}\cdot
    \int_{\Gamma^c_\theta(4r)} j(|z|)\,dz $ for   $ y\in \Gamma_\theta(r),$
where the comparison constant is independent of  $r\in (0, 1/4)$.
 Suppose the scale invariant BHP holds for $X$ in $\Gamma_\theta,$ then  there exists $c_1=c_1(d, \theta)>1$ such that for any
  $r\in (0, 1/4)$
and $y_1, y_2\in \Gamma_\theta(r),$
\begin{equation}\label{e:3.27'}
c_1^{-1}\dfrac{g(y_1)}{g(y_2)}\leq \dfrac{\E_{y_1} \tau_{\Gamma_\theta(2r)}}{\E_{y_2} \tau_{\Gamma_\theta(2r)}}\leq c_1\dfrac{g(y_1)}{g(y_2)}.
\end{equation}
It follows from  Theorem \ref{T:2.8} that there exists $c_2=c_2(d)>1$ such that for any $r\in (0, 1/4),$
$c_2^{-1}\E_x \tau_{\Gamma_\theta(4r)}\leq \E_x \tau_{\Gamma_\theta(2r)}\leq \E_x \tau_{\Gamma_\theta(4r)}$
for $x\in \Gamma_\theta(r).$
Hence, \eqref{e:3.27'} is equivalent that there exists $c_3=c_3(d, \theta)>1$ such that for any $r\in (0, 1/4)$
and $y_1, y_2\in \Gamma_\theta(r),$
\begin{equation}\label{e:3.27}
c_3^{-1}\dfrac{g(y_1)}{g(y_2)}\leq \dfrac{\E_{y_1} \tau_{\Gamma_\theta(4r)}}{\E_{y_2} \tau_{\Gamma_\theta(4r)}}\leq c_3\dfrac{g(y_1)}{g(y_2)}.
\end{equation}

 For each $\lambda>0,$ define $Y^{(\lambda)}_t:=\lambda X_{\lambda^{-2}t}$ and  $\phi^{(\lambda)}(s):=\lambda^{-2}\phi(\lambda^2 s)$ for $s>0.$ By \eqref{e:3.48}, $Y^{(\lambda)}_t$ is a subordinate Brownian motion with the characteristic exponent function  $\phi^{(\lambda)}(|\xi|^2).$
Let $\lambda>0.$
  Denote by $\E^{(\lambda)}$ and $G^{(\lambda)}$ the expectation and the Green function of the process $Y^{(\lambda)}.$
  For each open set $B,$ denote by $\tau^{(\lambda)}_B$ the first exit time of $Y^{(\lambda)}$ from $B.$
  Note that for $r>0,$
  \begin{equation}\label{e:3.57}
  p_{\Gamma_\theta(r)}(t,x , y)=r^{-d} p^{(1/r)}_{r^{-1}\Gamma_\theta(r)}(r^{-2} t, r^{-1}x, r^{-1}y), \quad t>0, x, y\in \Gamma_\theta(r),
  \end{equation}
  where $p_{\Gamma_\theta(r)}(t,x , y)$ and $p^{(1/r)}_{r^{-1}\Gamma_\theta(r)}(t, x, y)$ are the transition density functions of the part process of $X$ in $\Gamma_\theta(r)$ and the part process of $Y^{(1/r)}$ in $r^{-1}\Gamma_\theta(r).$
By integrating with respect to the time $t$ in \eqref{e:3.57}, we have the following scaling change formula
\begin{equation}\label{e:3.39}
G_{\Gamma_\theta(r)}(x, y)=r^{2-d}G^{(1/r)}_{r^{-1}\Gamma_\theta(r)} (r^{-1}x, r^{-1}y), \quad x,y \in \Gamma_\theta(r),
\end{equation}
and thus
\begin{equation}\label{e:3.40}
\E_x  [ \tau_{\Gamma_\theta(r)}] =r^2 \E^{(1/r)}_{r^{-1}x} \big [ \tau^{(1/r)}_{r^{-1}\Gamma_\theta(r)} \big], \quad x\in \Gamma_\theta(r).
\end{equation}
 Note that for any $r\in (0, 1),$ $r^{-1}\Gamma_\theta(4r)=\Gamma_\theta(4).$
 Hence,  for $y\in \Gamma_\theta(4r),$
\begin{equation}\label{e:3.68}
\begin{aligned}
g(y)
&=r^{2-d} \int_{B( y_0,  (3r\sin\theta\wedge  r)/2) )} G^{(1/r)}_{\Gamma_\theta(4)}(r^{-1}y, r^{-1}z)\,dz\\
&=r^2 \int_{B(r^{-1} y_0,  (3\sin\theta\wedge  1)/2))} G^{(1/r)}_{\Gamma_\theta(4)}(r^{-1}y, z)\,dz.
\end{aligned}
\end{equation}
Let $\tilde y_0:=r^{-1}y_0=(0, \cdots, 0, 3)$ and $g^{(r)}(y):=\int_{B(\tilde y_0,  (3\sin\theta\wedge  1)/2) )}
G^{(1/r)}_{\Gamma_\theta(4)}( y, z)\,dz$ for $y\in \Gamma_\theta(4).$
 Thus by \eqref{e:3.40} and \eqref{e:3.68},  \eqref{e:3.27} is equivalent  that for  any $r\in (0, 1/4)$
and $y_1, y_2\in \Gamma_\theta(1),$
\begin{equation}\label{e:c1.5}
c_3^{-1} \dfrac{g^{(r)}(y_1)}{g^{(r)}(y_2)}
\leq \dfrac{\E^{(1/r)}_{y_1} \big[  \tau^{(1/r)}_{\Gamma_\theta(4)} \big]}{\E^{(1/r)}_{y_2} \big[  \tau^{(1/r)}_{\Gamma_\theta(4)} \big]}
\leq c_3 \dfrac{g^{(r)}(y_1)}{g^{(r)}(y_2)}.
 \end{equation}

 Note that  $\lim_{\lambda \to \infty} \phi^{(\lambda)}(s) = s$ by \eqref{e:3.51}.
By \cite[Theorem 1.2 and Proposition 2.1]{CS2},
for each $f\in L^2(\Gamma_\theta(4), dy),$ $G^{(r^{-1})}_{\Gamma_\theta(4)}f$  converges strongly to
$G^{W}_{\Gamma_\theta(4)}f$ in $ L^2(\Gamma_\theta(4), dy)$ as $r\rightarrow 0.$
  Take $f={\mathbbm 1}_{\Gamma_\theta(4)}(\cdot)$ and $f={\mathbbm 1}_{B(\tilde y_0, (3\sin\theta\wedge 1)/2)}(\cdot)$, respectively,
There exists a subsequence $\{r_n\}\subset (0, \infty)$ such that $\lim_{n\to \infty} r_n = 0$  and  for
 almost every  $y$ in $\Gamma_\theta(4),$
$\lim_{n\to \infty} \E^{(1/r_n)}_{y} \big[  \tau^{(1/r_n)}_{\Gamma_\theta(4)} \big] = \E^W_y  \big[ \tau^W_{\Gamma_\theta(4)} \big]$
and
$ \lim_{n\to \infty} g^{(r_n)}(y) =  \int_{B(\tilde y_0,  (3\sin\theta\wedge  1)/2) )}G^W_{\Gamma_\theta(4)}(y, z)\,dz =: g^W(y) .$
 Hence, by taking $r_n\rightarrow 0$ in \eqref{e:c1.5}, we have  for almost every $y_1, y_2$ in $\Gamma_\theta(1),$
  \begin{equation}\label{e:0.7}
   c_3^{-1} \dfrac{g^W(y_1)}{g^W(y_2)}
\leq \dfrac{\E^W_{y_1} \big[  \tau^W_{\Gamma_\theta(4)} \big]}{\E^W_{y_2} \left[  \tau^W_{\Gamma_\theta(4)} \right]}
\leq c_3\dfrac{g^W(y_1)}{g^W(y_2)}.
\end{equation}
Note that for each $y\in \Gamma_\theta(1),$ $z\mapsto G^W_{\Gamma_\theta(4)}(y, z)$ is harmonic with respect to $W$ in $B(\tilde y_0, 3\sin\theta\wedge 1).$ By the uniform Harnack principle of $W,$ there exists $c_4=c_4(d)>1$ such that for any $y\in \Gamma_\theta(1),$
$$c_4^{-1}G^W_{\Gamma_\theta(4)}(y, \tilde y_0)\leq g^W(y)=\int_{B(\tilde y_0,  (3\sin\theta\wedge  1)/2) )}G^W_{\Gamma_\theta(4)}(y, z)\,dz\leq  c_4G^W_{\Gamma_\theta(4)}(y, \tilde y_0).$$
Hence by \eqref{e:0.7}, there  exists $c_5=c_5(d, \theta)>1$ such that for  almost every  $y_1, y_2$ in $\Gamma_\theta(1),$
 \begin{equation}\label{e:3.69}
 c_5^{-1} \dfrac{G^W_{\Gamma_\theta(4)}(y_1, \tilde y_0)}{G^W_{\Gamma_\theta(4)}(y_2, \tilde y_0)}
\leq \dfrac{\E^W_{y_1} \left[  \tau^W_{\Gamma_\theta(4)} \right]}{\E^W_{y_2} \left[  \tau^W_{\Gamma_\theta(4)} \right]}
\leq c_5 \dfrac{G^W_{\Gamma_\theta(4)}(y_1, \tilde y_0)}{G^W_{\Gamma_\theta(4)}(y_2, \tilde y_0)}.
\end{equation}
It is known that  $G^W_{\Gamma_\theta(4)}(y, \tilde y_0) $  and $\E^W_{y} \left[  \tau^W_{\Gamma_\theta(4)} \right]$  are continuous for  $y\in \Gamma_\theta(1).$
Hence, \eqref{e:3.69} holds for each $y_1$ and $y_2$ in $\Gamma_\theta(1).$ Thus,
 for
$\tilde x=(0, \cdots, 0, a)$ with $a<1,$
\begin{equation}\label{e:c1.6}
c_5^{-1} \dfrac{G^{W}_{\Gamma_\theta(4)}(\tilde x, \tilde y_0)}{G^{W}_{\Gamma_\theta(4)}(\tilde x_0, \tilde y_0)}
\leq \dfrac{\E^{W}_{\tilde x}  \left[ \tau^W_{\Gamma_\theta(4)} \right]}{\E^{W}_{\tilde x_0}  \left[ \tau^W_{\Gamma_\theta(4)} \right]}
\leq c_5\dfrac{ G^{W}_{\Gamma_\theta(4)}(\tilde x, \tilde y_0)}{G^{W}_{\Gamma_\theta(4)}(\tilde x_0, \tilde y_0)}.
\end{equation}
where   $\tilde x_0=(0, \cdots, 0, 1/2).$
Note that there is $c_6=c_6(d)>1$ such that
$$c_6\geq \E^{W}_{\tilde x_0} \tau^{W}_{B(\tilde x_0, 4)}\geq \E^{W}_{\tilde x_0}  \tau^{W}_{\Gamma_\theta(4)}
\geq \E^{W}_{\tilde x_0} \tau^{W}_{B(\tilde x_0, 1/4\sin\theta)}\geq c_6^{-1}.$$
By \cite[Lemma 6.7]{ChungZ}, there is $c_7=c_7(d, \theta)>1$ so that
$c_7^{-1}\leq G^{W}_{\Gamma_\theta(4)}(\tilde x_0, \tilde y_0)\leq c_7.$
Hence, \eqref{e:c1.6} is equivalent that there is $c_8=c_8(d, \theta)>1$ such that for
$\tilde x=(0, \cdots, 0, a)$ with $a<1,$
\begin{equation}\label{e:c1.6'}
c_8^{-1} G^{W}_{\Gamma_\theta(4)}(\tilde x, \tilde y_0)
\leq \E^{W}_{\tilde x}  \left[ \tau^W_{\Gamma_\theta(4)} \right]
\leq c_8 G^{W}_{\Gamma_\theta(4)}(\tilde x, \tilde y_0).
\end{equation}
However, by  Lemma \ref{L:4.4'} and \eqref{e:3.56}, \eqref{e:c1.6'} does not hold for $\Delta$  if $\theta\in (0, \cos^{-1}(1/\sqrt d)].$
  This contradiction shows that  the scale invariant  BHP fails to hold for $X$  on  the truncated cone $\Gamma_\theta$ with
   $\theta\in (0, \cos^{-1}(1/\sqrt d)]$.
\qed

\begin{remark}\rm By \cite[Theorem 3.5  and Example 5.5]{BKK2}  and \eqref{e:3.19},
 non-scale invariant boundary Harnack principle holds on   open sets for
any subordinate  Brownian motion  with Gaussian component that satisfies  condition \eqref{e:mu1}.
\end{remark}

\section{Appendix}

In the Appendix,  we provide a detailed proof for Pruitt's result Theorem \ref{T:0.3}  for any  L\'evy process in $\R^d$. Let $Y$ be a  L\'evy process in $\R^d$  with L\'evy triple $(A, b, \nu).$
 Recall that the function $\Phi (r)$ is defined on $(0, \infty)$ by \eqref{e:Phi}.

\begin{lem}\label{L:1.3}
Let $Y$ be a L\'evy process in $\R^d$ starting from $Y_0=0$ and $M_t:=\sup_{0\leq s\leq t}|Y_s|.$
There exists $c=c(d)>1$ such that for any $r>0,$
\begin{equation}\label{e:0.8}
\P(M_t\geq r)\leq ct \Phi(r) \quad \hbox{and} \quad \P(M_t\leq r)\leq \dfrac{c}{t\Phi(r)}.
\end{equation}
\end{lem}

\pf  Let $(A, b, \nu)$ be the L\'evy triple of $Y$ so the   L\'evy-Khintchine  $\psi$ of $Y$
is given by \eqref{e:2.25}.
Fix $r>0.$ Define for $u\in \R^d$,
\begin{equation}\label{e:0.9a}
\begin{aligned}
\psi_1(u)=& -\dfrac{1}{2} (u, Au)+i(\tilde b_r, u)
+\int_{|z|\leq r} \left(e^{i(u, z)}-1-   i(u, z)\right)\nu(dz),
\end{aligned}\end{equation}
\begin{equation}\label{e:0.10a}
 \psi_2(u)=\int_{|z|>  r} (e^{i(u, z)}-1)\nu(dz).
\end{equation}
Note that  $\psi(u)=\psi_1(u)+\psi_2(u).$ Let $Y^{(1)}$ and $Y^{(2)}$ be the two independent L\'evy  processes
with L\'evy characteristics $\psi_1$ and $\psi_2$, respectively.
 In other words, $Y^{(1)}$ is a L\'evy  process  with L\'evy triple
$(A, \tilde b_r,  {\mathbbm 1}_{\{|z|\leq r\}} \nu (dz) )$ and $Y^{(2)}$ is a (spatial) Poisson point process
with    intensity measure ${\mathbbm 1}_{\{|z|>r\}} \nu (dz)$ that is independent of $Y^{(1)}.$
Define
$K(r):=r^{-2}\left(\sum_{i=1}^d a_{ii}+\int_{|z|\leq r} |z|^2\nu(dz)\right),$
$L(r):=r^{-1}|\tilde b_r|$ and $G(r):=\nu\{z: |z|>r\}.$
Then $\Phi(r)=G(r)+K(r)+L(r).$

By differentiation of the characteristic function,
$\E Y^{(1)}_t=-it \psi'_1(0)=t\tilde b_r$ and
$\E|Y^{(1)}_t- \E Y^{(1)}_t|^2=t r^2K(r).$
For the compound Poisson process $Y^{(2)}_t,$
$\P(Y^{(2)}_s\neq 0\: \mbox{for some}\: s\leq t)=1-e^{-tG(r)}\leq tG(r).$
Fix $t>0.$ Note that $|\E Y^{(1)}_t|= r t L(r).$ If $|\E Y^{(1)}_t|\geq r/2,$ we have
 $$\P(M_t\geq r)\leq 1 \leq \frac{2 |\E Y^{(1)}_t|}r
 \leq 2tL(r)\leq 2t \Phi(r).$$
On the other hand, if $|\E Y^{(1)}_t|< r/2,$ then
$|\E Y^{(1)}_s| = \frac{s}{t} |\E Y^{(1)}_t| <r/2$ for all $s\in [0, t],$ so
by  Doob's martingale inequality,
 $$\begin{aligned}
 \P(M_t\geq r)&\leq \P(Y^{(2)}_s\neq 0\: \mbox{for some}\: s\leq t)+\P(\sup_{0\leq s\leq t}|Y^{(1)}_s|\geq r)\\
 &\leq tG(r)+\P(\sup_{0\leq s\leq t}|Y^{(1)}_s-\E Y^{(1)}_s|\geq r/2)\\
 &\leq tG(r)+\dfrac{\E|Y^{(1)}_t- \E Y^{(1)}_t|^2}{(r/2)^2}\\
 &\leq tG(r)+\dfrac{tr^2 K(r)}{(r/2)^2}\leq 4t \Phi(r).\\
 \end{aligned}$$
 That is, the first inequality in \eqref{e:0.8} holds.

 We will prove the second inequality in \eqref{e:0.8}   by dividing into    three cases according to which one is the largest
 among $\{K(2r), G(2r), L(2r)\}$.
In the following, with a little abuse of notation, we  let $\psi_1$ and $\psi_2$ be the characteristics defined in \eqref{e:0.9a} and \eqref{e:0.10a} with $2r$ in place of $r,$ let $Y^{(1)}$ and $Y^{(2)}$ be the corresponding L\'evy processes with characteristics $\psi_1$ and $\psi_2$.

  (i) Suppose $K(2r)\geq  \max\{ G(2r), L(2r)\}.$ In this case,  $K(2r)\geq  \frac{1}{3}\Phi(2r).$
 For any one-dimensional random variable $Z$ in $\R$ with characteristic function $g$,
 as  is done in  \cite[p.955]{P}, for any $a>0$,
$$\frac{\cos 1}2 \,  a^2 \P(a|Z|\leq 1)\leq \int_0^a \big( \int_0^a  |g(s) | ds \big) dt \leq a  \int_0^a  |g(s) | ds.$$
 Thus    with  $c_1 = \frac{2}{\cos 1}$,
 \begin{equation}\label{e:0.9}
 \P(|Z|\leq r)\leq c_1  r \int_0^{1/r}|g(s)|\,ds \quad \hbox{for every } r>0.
 \end{equation}

Fix $t>0,$ write  $Y_t=(Y^{(1)}_t, \cdots, Y^{(d)}_t).$ For each $i\geq 1,$ let $g_{Y^{(i)}_t}$ be the characteristic function of $Y^{(i)}_t.$
Denote by $g_{Y_t}$ the characteristic functions of $Y_t.$
For each  constant $\omega$ in $\R,$ let ${\bf \tilde \omega_i}=(0, \cdots, 0, \overset{(i)}{\omega}, 0, \cdots, 0)\in\R^d.$
There exists $c_2\in (0, 1)$ such that for each $\omega\in (0, (2r)^{-1}),$
\begin{eqnarray}\label{e:0.10}
|g_{Y^{(i)}_t}(\omega)|&=&
 |g_{Y_t}({\bf \tilde \omega_i})|=\exp(t {\rm Re} \psi({\bf \tilde\omega_i})) \nonumber \\
& \leq &\exp\big(t\big[-\dfrac{1}{2} ({\bf \tilde\omega_i}, A{\bf \tilde \omega_i})+\int_{\{z\in\R^d: |z|\leq 2r\}} (\cos(\omega z_i)-1)\nu(dz)\big]\big)  \nonumber \\
& \leq &\exp\big(t\omega^2\big[-\dfrac{1}{2} a_{ii}-c_2\int_{\{z\in\R^d: |z|\leq 2r\}} z_i^2\nu(dz)\big]\big)
\end{eqnarray}
Note that
$\P(|Y_t|\leq 2r)\leq \inf_{1\leq i\leq d}\P(|Y^{(i)}_t|\leq 2r).$
Hence, by \eqref{e:0.9} and \eqref{e:0.10},
\begin{eqnarray*}
&& \P(|Y_t|\leq 2r)
 \leq \inf_{1\leq i\leq d} \P(|Y^{(i)}_t|\leq 2r)\\
&\leq & 2c_1r\int_0^{(2r)^{-1}}\inf_{1\leq i\leq d}\exp\big(t\omega^2\big[-\dfrac{1}{2} a_{ii}-c_2\int_{\{z\in\R^d: |z|\leq 2r\}} z_i^2\nu(dz)\big]\big)\,d\omega\\
&=& 2c_1r\int_0^{(2r)^{-1}}\exp\big(t\omega^2\big[-\dfrac{1}{2} \sup_{1\leq i\leq d}a_{ii}-c_2\sup_{1\leq i\leq d}\int_{\{z\in\R^d: |z|\leq 2r\}} z_i^2\nu(dz)\big]\big)\,d\omega\\
&\leq & 2c_1r\int_0^{(2r)^{-1}}\exp\big(t\omega^2d^{-1}\big[-\dfrac{1}{2} \sum_{i=1}^d a_{ii}-c_2\sum_{i=1}^d\int_{\{z\in\R^d: |z|\leq 2r\}} \sum_{i=1}^d z_i^2\nu(dz)\big]\big)\,d\omega\\
&\leq & 2c_1r\int_0^\infty\exp\big(-(c_2\wedge\frac{1}{2})t\omega^2d^{-1} 4r^2 K(2r)\big)\,d\omega
= \dfrac{c_3}{(tK(2r))^{1/2}}.
\end{eqnarray*}
Therefore,
$$
\P(M_t\leq r) \leq \P (|Y_{t/2}|\leq r)\P(|Y_t-Y_{t/2}|\leq 2r)
 \leq \dfrac{2c_3^2}{tK(2r)}\leq \dfrac{ 6 c_3^2}{t\Phi(2r)}.
 $$

\smallskip

(ii) Suppose $G(2r)\geq   \max\{ L(2r), K(2r)\}.$ Then $G(2r)\geq  \frac{1}{3}\Phi(2r)$ and so
$$
\P(M_t\leq r)
\leq  \P(Y^{(2)}_s= 0\: \mbox{for any}\: 0\leq s\leq t)
=e^{-tG(2r)}\leq (tG(2r))^{-1}\leq \frac{3}{t\Phi(2r)}.
$$

\smallskip

(iii) Suppose $L(2r)\geq   \max\{ G(2r), K(2r)\}.$ If $tL(2r)\leq 2,$ then
$\P(M_t\leq r)\leq 1\leq \frac{2}{tL(2r)}.$
When $tL(2r)> 2,$ $|\E Y^{(1)}_t|=2r t L(2r)\geq 4r.$ In this case,
\begin{eqnarray*}
\P(M_t\leq r)&\leq & \P \big(Y^{(2)}_s= 0\: \mbox{for any}\: 0\leq s\leq t
 \hbox{ and } \sup_{0\leq s\leq t}|Y^{(1)}_s|\leq r  \big)\\
&\leq&  \P \Big(|Y^{(1)}_t|\leq  r  \Big)
\, \leq \, \P \Big( |Y^{(1)}_t-\E Y^{(1)}_t|\geq  3 r t L(2r)/2  \Big)\\
&\leq&  \dfrac{\E \big[ |Y^{(1)}_t-\E Y^{(1)}_t|^2\big] }{ (3 rt  L(2r)/2)^2  }
\leq  \dfrac{t (2r)^2 K(2r)}{9t^2 r^2 L(2r)^2/4}   \leq \dfrac{2}{ tL(2r)}\leq   \dfrac{ 6} {t\Phi(2r)} .
\end{eqnarray*}
\qed

\medskip

\noindent {\bf Proof of Theorem \ref{T:0.3}.}
Note that  for any $r>0$,
\begin{equation}\label{e:1.8}
\E_x\tau_{B(x, r)}=\E_0 \tau_{B(0, r)}=\int_0^\infty \P_0(\tau_{B(0, r)}>t)\,dt=\int_0^\infty\P(M_t\leq r) \,dt.
\end{equation}
Since
$
\{M_{2t}\leq r\}\subset  \Big\{M_t\leq r \hbox{ and } \sup_{s\in [0, t]}|Y_{t+s}-Y_t|\leq 2r\Big\},
$
we have by the independent increments property of $Y$,  Lemmas \ref{L:1.1} and \ref{L:1.3}
that
\begin{eqnarray}\label{e:1.9}
\P (M_{2t} \leq r) &\leq&
\P(M_t\leq r) \, \P \Big( \sup_{s\in [0, t]}|Y_{t+s}-Y_t|\leq 2r \Big) =
 \P(M_t\leq r)\P(M_t\leq 2r) \nonumber \\
&\leq& \frac{c_1^2}{t^2\Phi (r) \Phi (2r)} \leq \frac{c_2}{t^2\Phi(r)^2}.
\end{eqnarray}
Hence,
\begin{eqnarray*}
&&\E_x \tau_{B(x, r)} = \frac12\int_0^\infty \P(M_{2t}\leq r)dt\\
&\leq & \frac12\int_0^{1/\Phi(r)}1 \,dt+\frac12\int_{1/\Phi(r)}^\infty \frac{c_2}{t^2\Phi(r)^2}\,dt
\leq \frac12(1+c_2)\dfrac{1}{\Phi(r)}.
\end{eqnarray*}

 For the lower bound,
 by the first inequality in \eqref{e:0.8}, there exists $c_3=c_3(d)>1$ such that
$\P(M_t\geq r)\leq c_3t \Phi(r)$ for  $ r>0.$
For each  $r>0$, take $t_0 >0$ so that  $c_3t_0 \Phi(r)=1/2$.
It then follows from \eqref{e:1.8} that
$$
\E_x \tau_{B(x, r)} \geq \int_0^{t_0} \P (M_{t} <r) dt \geq t_0/2 = \frac{1}{4c_3 \Phi (r)}.
$$
 This establishes Theorem \ref{T:0.3}.
  \qed

\medskip

{\bf Acknowledgement.} The authors are grateful to Pat Fitzsimmons for the discussion of  conservative and dissipative decomposition of an excessive measure in \cite{FM}.

\end{document}